\documentclass{amsart}
\usepackage{latexsym,amscd,amssymb,epsfig}

\theoremstyle{plain}
  \newtheorem{theorem}{Theorem}[section]
  \newtheorem{proposition}[theorem]{Proposition}
  \newtheorem{lemma}[theorem]{Lemma}
  \newtheorem{corollary}[theorem]{Corollary}
  
  \newtheorem{conjecture}[theorem]{Conjecture}
  \newtheorem{exercise}[theorem]{Exercise}
\theoremstyle{definition}
  \newtheorem{definition}[theorem]{Definition}
  \newtheorem{example}[theorem]{Example}
  \newtheorem{examples}[theorem]{Examples}
  \newtheorem{question}[theorem]{Question}
\theoremstyle{remark}
  \newtheorem{remark}[theorem]{Remark}

\numberwithin{equation}{section}
\numberwithin{figure}{section}

\DeclareMathOperator{\Line}{Line}

\def\Sym{\mathfrak{S}}
\def\peak{\mathrm{peak}}

\def\LL{\mathcal{L}}
\def\F{\mathcal{F}}
\def\PS{\mathrm{PS}}

\def\nw{\,{\nwarrow}}

\def\Hom{\mathrm{Hom}}

\def\Vertices{\mathrm{Vertices}}
\def\Span{\mathrm{Span}}
\def\Des{\mathrm{Des}}
\def\des{\mathm{des}}

\def\N{\mathcal{N}}

\def\B{\mathcal{B}}

\def\des{\mathrm{des}}
\def\wdes{\mathrm{wdes}}

\def\Z{\mathbb{Z}}
\def\R{\mathbb{R}}

\def\C{\mathbb{C}}
\def\A{\mathrm{Path}}
\def\cycle{\mathrm{Cycle}}
\def\Kite{\mathrm{Kite}}

\title{Faces of generalized permutohedra}

\author{Alexander Postnikov}
\address{Massachusetts Institute of Technology}
\email{apost at math.mit.edu}

\author{Victor Reiner}
\address{University of  Minnesota}
\email{reiner at math.umn.edu}

\author{Lauren Williams}
\address{Harvard University}
\email{lauren at math.harvard.edu}

\date{April 28, 2007}

\thanks{A.P.\ was supported in part by NSF CAREER Award DMS-0504629.  
V.R.\ was supported in part by NSF Grant DMS-0601010.
L.W.\ was supported in part by an NSF Postdoctoral Fellowship.}

\begin{document}

\begin{abstract}
The aim of the paper is to calculate face numbers of simple 
generalized permutohedra, and
study their $f$-, $h$- and $\gamma$-vectors.
These polytopes include permutohedra, associahedra,
graph-associahedra, simple graphic zonotopes, nestohedra, and
other interesting polytopes.  

We give several explicit formulas for $h$-vectors and $\gamma$-vectors
involving descent statistics.  This includes a combinatorial interpretation 
for $\gamma$-vectors of a large class of generalized
permutohedra which are flag simple polytopes, and confirms for them
Gal's conjecture on nonnegativity of $\gamma$-vectors.

We calculate explicit generating functions and formulae 
for $h$-polynomials of various families
of graph-associahedra, including those
corresponding to all Dynkin diagrams of finite and affine types.
We also 
discuss relations with Narayana numbers 
and with Simon Newcomb's problem.

We give (and conjecture) upper and lower bounds for $f$-, $h$-, and $\gamma$-vectors
within several classes of generalized permutohedra.

An appendix discusses the equivalence of various notions of deformations of simple polytopes.
\end{abstract}

\maketitle

\tableofcontents

\section{Introduction}
\label{sec:intro}

Generalized permutohedra are a very well-behaved class of convex polytopes
studied in~\cite{Post}, as generalizations of the classical permutohedra,
associahedra, cyclohedra, etc.  That work explored their wonderful properties
from the point of view of valuations such as volumes, mixed volumes, and number
of lattice points.  This paper focuses on their further good behavior with
respect to face enumeration in the case when they are simple polytopes.

Simple generalized permutohedra include as an important subclass (the duals of)
the {\it nested set complexes\/} considered by DeConcini and Procesi in their
work on {\it wonderful compactifications\/} of hyperplane arrangements; see
\cite{DP, FeichtnerSturmfels}.  In particular, when the arrangement comes from
a Coxeter system, one obtains interesting flag simple polytopes 
studied by Davis, Januszkiewicz, and Scott \cite{DJS}.  These polytopes can be
combinatorially described in terms of the corresponding Coxeter graph.
Carr and Devadoss \cite{CarrDevadoss} studied these polytopes for
arbitrary graphs and called them {\it graph-associahedra.} 

We mention here two other recent papers in which generalized permutohedra have
appeared.  
Morton, Pachter, Shiu, Sturmfels, and Wienand \cite{MPSSW} considered generalized 
permutohedra from the point of view of {\it rank tests\/} on ordinal data in statistics.   The
normal fans of generalized permutohedra are what they called {\it submodular rank
tests}.  Agnarsson and Morris \cite{AgnarssonMorris} investigated closely the
$1$-skeleton (vertices and edges) in the special case where generalized
permutohedra are Minkowski sums of standard simplices.

\medskip

Let us formulate several results of the present paper.  A few definitions are required.
A {\it connected building set\/} $\B$ on $[n]:=\{1,\dots,n\}$
is a collection of nonempty subsets in $[n]$
such that 
\begin{enumerate}
\item 
if $I,J\in \B$ and $I\cap J\ne \emptyset$, then $I\cup J\in \B$,
\item
$\B$ contains all singletons $\{i\}$ and the whole set $[n]$,
\end{enumerate}
see Definition~\ref{def:building-set-definition}.
An interesting subclass of {\it graphical building sets\/} $\B(G)$ 
comes from connected graphs $G$ on $[n]$.  The building set $\B(G)$ contains all 
nonempty subsets of vertices $I\subseteq[n]$ such that the induced graph $G|_I$ is connected.

The {\it nestohedron\/} $P_\B$ is defined (see Definition~\ref{def:nestohedron})
as the Minkowski sum 
$$
P_\B = \sum_{I\in \B} \Delta_I
$$
of the coordinate simplices $\Delta_I:=\mathrm{ConvexHull}(e_i\mid i\in I)$, 
where the $e_i$ are the endpoints of the coordinate vectors in $\R^n$.
According to \cite[Theorem~7.4]{Post} and \cite[Theorem~3.14]{FeichtnerSturmfels}
(see Theorem~\ref{th:nestohedra_dual_to_nested_complexes} below),
the nestohedron $P_\B$ is a simple polytope which is dual to a simplicial nested set complex.
For a graphical building set $\B(G)$, the nestohedron $P_{\B(G)}$ is called
the {\it graph-associahedron.}  In the case when $G$ is the $n$-path, 
$P_{\B(G)}$ is the usual associahedron; and in the case when $G=K_n$ is the complete graph,
$P_{\B(G)}$ is the usual permutohedron.

Recall that the $f$-vector and the $h$-vector of a simple $d$-dimensional polytope
$P$
are $(f_0,f_1,\dots,f_d)$ and $(h_0,h_1,\dots,h_d)$, where $f_i$ is the number
of
$i$-dimensional faces of $P$ and $\sum h_i \, (t+1)^i = \sum f_i \,t^i$.
It is known that the $h$-vector of a simple polytope is 
positive and symmetric.
Since the $h$-vector is symmetric, one can define
another vector called the $\gamma$-vector 
$(\gamma_1,\gamma_2,\dots,\gamma_{\lfloor d/2\rfloor})$ 
by the relation
$$
\sum_{i=0}^d h_i\, t^i = \sum_{i=0}^{\lfloor \frac{d}{2} \rfloor}  \gamma_i \, t^i (1+t)^{d-2i}.
$$

A simplicial complex $\Delta$ is called a {\it flag complex\/} 
(or a {\it clique complex\/}) if its simplices are cliques (i.e., subsets of vertices
with complete induced subgraphs) of some graph (1-skeleton of $\Delta$).
Say that a simple polytope is {\it flag\/} if its dual
simplicial complex is flag.

Gal conjectured~\cite{Gal} that the $\gamma$-vector has nonnegative entries for any flag simple polytope.

Let us that say a connected building set $\B$ is {\it chordal\/} if, for any of
the sets $I=\{i_1<\dots<i_r\}$ in $\B$, all subsets $\{i_s,i_{s+1},\dots,i_r\}$ also belong to $\B$;
see Definition~\ref{def:chordal_building}.
By Proposition~\ref{prop:graphical_chordal}, graphical chordal building sets $\B(G)$ are exactly 
building sets coming from chordal graphs.
By Proposition~\ref{prop:chordal_are_flag}, all nestohedra $P_\B$ for chordal building sets
are flag simple polytopes.  So Gal's conjecture applies to this class of {\it chordal nestohedra,}  which include graph-associahedra for chordal graphs
and, in particular, for trees.

For a building set $\B$ on $[n]$, define (see Definition~\ref{def:characterization})
the set $\Sym_n(\B)$ of {\it $\B$-permutations\/} as the set of permutations $w$ of size $n$
such that, for any $i=1,\dots,n$, there exists $I\in \B$ such that
$I\subseteq \{w(1),\dots,w(i)\}$, and $I$ contains both $w(i)$ and $\max\{w(1),w(2),\dots,w(i)\}$.
It turns out that $\B$-permutations are in bijection with vertices of the nestohedron
$P_\B$; see Proposition~\ref{prop:B-trees-B-permutations}.

Let $\des(w)=\#\{i\mid w(i)>w(i+1)\}$ denote the number of descents in a permutation $w$.
Let $\widehat{\Sym}_n$ be the subset of permutations $w$ of size $n$ without two consecutive
descents and without final descent, i.e., there is no $i\in [n-1]$ such that $w(i)>w(i+1)>w(i+2)$,  
assuming that $w(n+1)=0$.

\begin{theorem}  
{\rm (Corollary~\ref{cor:h_des_chordal} and Theorem~\ref{th:gamma-chordal})}
Let $\B$ be a connected chordal building set on $[n]$.  
Then the $h$-vector of the nestohedron $P_\B$ is given by
$$
\sum_i h_i \, t^i = \sum _{w\in \Sym_n(\B)} t^{\des(w)},
$$
and the $\gamma$-vector of the nestohedron $P_\B$ is given by
$$
\sum_i \gamma_i \, t^i = \sum _{w\in \Sym_n(\B)\cap \widehat{\Sym}_n} t^{\des(w)}.
$$
\end{theorem}

This result shows that Gal's conjecture is true for chordal nestohedra.

\medskip

The paper is structured as follows.

Section~\ref{sec:fpolynomial} reviews polytopes, cones, fans, and gives 
basic terminology of face enumeration for polytopes ({\it $f$-vectors}),
simple polytopes ({\it $h$-vectors}), and {\it flag\/} simple polytopes 
({\it $\gamma$-vectors\/}). 

Section~\ref{sec:cones-fans-preposets} reviews the definition of generalized
permutohedra, and recasts this definition equivalently in terms of their normal
fans.
It then sets up the dictionary between preposets, and cones and fans coming from the braid
arrangement.  In particular, one finds that each vertex in a generalized
associahedron has associated to it a {\it poset\/} that describes its normal
cone.  This is used to characterize when the polytope is simple, namely when the
associated posets have Hasse diagrams which are trees.
In Section~\ref{sec:simples} this leads to a combinatorial formula formula for
the $h$-vector in terms of descent statistics on these tree-posets.

The remainder of the paper deals with subclasses of simple generalized
permutohedra.  Section~\ref{sec:zonotopes} dispenses quickly with the very
restrictive class of simple zonotopal generalized permutohedra, namely the 
{\it simple graphic zonotopes}.

Section~\ref{sec:nested} then moves on to the interesting class of
{\it nestohedra\/} $P_\B$ coming from a building set $\B$, where the
posets associated to each vertex are rooted trees.  These include graph-associahedra. 
Section~\ref{sec:FlagNestohedra} characterizes the flag nestohedra.

Section~\ref{sec:B-trees-permutations} discusses
$\B$-trees and $\B$-permutations.
These trees and permutations are in bijection with each other
and with vertices of the nestohedron $P_\B$.
The $h$-polynomial of a nestohedron is the descent-generating
function for $\B$-trees.
Then Section~\ref{sec:chordal} introduces the class
of chordal building sets and shows that $h$-polynomials of their
nestohedra are descent-generating functions for $\B$-permutations.

Section~\ref{sec:three-graph-ass-examples} illustrates
these formulas for $h$-polynomials by several examples: 
the classical permutohedron and associahedron, the cyclohedron,
the stellohedron (the graph-associahedron for the star graph), 
and the Stanley-Pitman polytope.

Section~\ref{sec:gamma_nestohedra} gives a combinatorial
formula for the $\gamma$-vector of all chordal nestohedra
as a descent-generating function (or peak-generating function) 
for a subset of $\B$-permutations. 
This result implies Gal's nonnegativity conjecture 
for this class of polytopes.
The warm-up example here is the classical permutohedron, and the section concludes
with the examples of the associahedron and cyclohedron.

Section~\ref{sec:one-branching} calculates the generating functions for
$f$-polynomials of the graph-associahedra for all trees with one branching point
and discuss a relation with Simon Newcomb's problem.
Section~\ref{sec:almost-paths} deals with graphs that are formed by a path with
two small fixed graphs attached to the ends.  It turns out that the $h$-vectors
of graph-associahedra for such  path-like graphs can be expressed 
in terms of $h$-vectors of classical associahedra.
The section includes explicit formulas for graph-associahedra for 
the Dynkin diagrams of all finite and affine Coxeter groups.

Section~\ref{sec:flossing} gives some bounds and
monotonicity conjectures for face numbers of generalized permutohedra.  

The paper ends with an Appendix which clarifies the equivalence between various
kinds of deformations of a simple polytope.

\medskip
\noindent
{\sc Acknowledgments:}  The authors thank
Federico Ardila, 
Richard Ehrenborg, 
Ira Gessel, 
Sangwook Kim, 
Jason Morton,
Margaret Readdy,
Anne Shiu, 
Richard Stanley, 
John Stembridge,
Bernd Sturmfels,  
Oliver Wienand,
and Andrei Zelevinsky 
for helpful conversations.

\section{Face numbers}
\label{sec:fpolynomial}

This section recalls some standard definitions from the theory of convex
polytopes and formulate Gal's extension of the Charney-Davis conjecture.

\subsection{Polytopes, cones, and fans}
\label{ssec:cones_fans}

A {\it convex polytope\/} $P$ is the convex hull of a finite collection of
points in $\R^n$.  The {\it dimension\/} of a polytope (or any other subset in
$\R^n$) is the dimension of its affine span.

A {\it polyhedral cone\/} in $\R^n$ is a subset defined by a conjunction of
weak inequalities of the form $\lambda(x) \geq 0$ for linear forms $\lambda
\in (\R^n)^*$.
A {\it face\/} of a polyhedral cone is a subset of the cone given by replacing
some of the inequalities $\lambda(x)\geq 0$ by the equalities $\lambda(x)=0$.


Two polyhedral cones $\sigma_1, \sigma_2$ {\it intersect properly\/} if
their intersection is a face of each.
A {\it complete fan of cones\/} $\F$ in $\R^n$ is a collection
of distinct nonempty polyhedral cones covering $\R^n$ such that
(1) every nonempty face of a cone in $\F$ is also a cone in
$\F$,  and
(2) any two cones in $\F$ intersect properly.
Cones in a fan $\F$ are also called {\it faces\/} of~$\F$.

Note
that fans can be alternatively defined only in terms of their top dimensional
faces, as collections of distinct pairwise properly
intersecting $n$-dimensional cones
covering $\R^n$.







A {\it face\/} $F$ of a convex polytope $P$ is the set of points in $P$ where
some linear functional $\lambda\in (\R^n)^*$ achieves its maximum on $P$, i.e.,
$$
F=\{x\in P\mid \lambda(x) = \max\{\lambda(y)\mid y\in P\}\}.
$$
Faces that consist of a
single point are called {\it vertices\/} and 1-dimensional faces are called
{\it edges\/} of $P$.

Given any convex polytope $P$ in $\R^n$ and a face $F$ of $P$, the {\it normal
cone to $P$ at $F$}, denoted $\N_F(P)$, is the subset of linear functionals
$\lambda\in(\R^n)^*$ whose maximum on $P$ is achieved
on all of the points in the face
$F$, i.e.,
$$
\N_F(P) := \{\lambda \in (\R^n)^* \mid \lambda(x) =
\max\{\lambda(y)\mid y\in P\}
\textrm{ for all } x\in F\}.
$$
Then $\N_F(P)$ is a polyhedral cone in $(\R^n)^*$, and the collection of all
such cones $\N_F(P)$ as one ranges through all faces $F$ of $P$ gives a
complete fan in $(\R^n)^*$ called the {\it normal fan\/} $\N(P)$.  A fan of the
form $\N(P)$ for some polytope $P$ is called a {\it polytopal fan.}

The combinatorial structure of faces of $P$ can be encoded by the {\it lattice
of faces\/} of $P$ ordered via inclusion.  This structure is also encoded by
the normal fan $\N(P)$.  Indeed, the map $F \mapsto \N_F(P)$ is an
inclusion-reversing bijection between the faces of $P$ and the faces of
$\N(P)$.

A cone is called {\it pointed\/} if it contains no lines (1-dimensional linear
subspaces), or equivalently, if it can be defined by a conjunction of
inequalities $\lambda_i(x) \geq  0$ in which the $\lambda_i$ span $(\R^n)^*$.
A fan is called {\it pointed\/} if all its faces are pointed.

If the polytope $P\subset \R^n$ is full-dimensional, that is $\dim P = n$, then
the normal fan $\N(P)$ is pointed.  For polytopes $P$ of lower dimension $d$,
define the $(n-d)$-dimensional subspace $P^\perp \subset (\R^n)^*$
of linear functionals which are constant on $P$.  Then all cones in the normal
fan $\N(P)$ contain the subspace $P^\perp$.  Thus $\N(P)$ can be reduced to a
pointed fan in the space $(\R^n)^*/P^\perp$.

A polytope $P$ is called {\it simple\/}  if any vertex of $P$ is incident to
exactly $d = \dim P$ edges.  A cone is called {\it simplicial\/} if it can be
given by a conjunctions of linear inequalities $\lambda_i(x)\geq 0$ and linear
equations $\mu_j(x) = 0$ where the covectors $\lambda_i$ and $\mu_j$ form a
basis in $(\R^n)^*$.  A fan is called {\it simplicial\/} if all its faces are
simplicial.  Clearly, simplicial cones and fans are pointed.
A convex polytope $P\subset \R^n$ is simple if and only if its (reduced)
normal fan $\N(P)/P^\perp$ is simplicial.


The {\it dual simplicial complex\/} $\Delta_P$ of a simple polytope $P$ is the
simplicial complex obtained by intersecting the (reduced) normal fan
$\N(P)/P^\perp$ with the unit sphere.  Note that $i$-simplices of $\Delta_P$
correspond to faces of $P$ of codimension $i+1$.




\subsection{$f$-vectors and $h$-vectors}

For a $d$-dimensional polytope $P$, the {\it face
number\/} $f_i(P)$ is the number of $i$-dimensional faces of $P$.  The
vector $(f_0(P),\dots,f_d(P))$ is called the {\it $f$-vector\/},
and the polynomial $f_P(t) = \sum_{i=0}^d f_i(P)\, t^i$ is called
the {\it $f$-polynomial\/} of $P$.

Similarly, for a $d$-dimensional fan $\F$, $f_i(\F)$ is the number
of $i$-dimensional faces of $\F$, and
$f_\F(t) = \sum_{i=0}^d f_i(\F)\,t^i$.
Note that face numbers of a polytope $P$ and its (reduced) normal cone
$\F = \N(P)/P^\perp$ are related as
$f_i(P) = f_{d-i}(\F)$, or equivalently, $f_P(t) = t^d\,f_\F(t^{-1})$.

We will most often deal with the case where $P$ is a simple polytope,
or equivalently, when $\F$ is a simplicial fan.
In these situations, there is a more compact encoding of the face numbers $f_i(P)$
or $f_i(\F)$ by smaller nonnegative integers.  One defines the {\it
$h$-vector} $(h_0(P),\dots,h_d(P))$ and {\it $h$-polynomial}
$h_P(t) = \sum_{i=0}^d h_i(P)\, t^i$ uniquely by the relation
\begin{equation}
\label{eq:f-h-relation}
f_P(t) = h_P(t+1),
\quad \textrm{or equivalently,} \quad
f_j(P) = \sum_{i} \binom{i}{j} h_i(P),\  j=0,\dots,d.
\end{equation}
For a simplicial fan $\F$, the $h$-vector $(h_0(\F),\dots,h_d(\F)$
and the  $h$-polynomial $h_\F(t) = \sum_{i=0}^d h_i(\F)\,t^i$
are defined by the relation
$t^d \,f_\F(t^{-1}) = h_\F(t+1)$,
or equivalently,
$f_j(\F) = \sum_{i} \binom{i}{d-j} h_i(\F)$, for $j=0,\dots,d$.
Thus the $h$-vector of a simple polytope
coincides with the $h$-vector of its normal fan.

The nonnegativity of $h_i(P)$ for a simple polytope $P$ comes from its
well-known combinatorial interpretation \cite[\S8.2]{Ziegler} in terms of the
$1$-skeleton of the simple polytope $P$.
Let us extend this interpretation to arbitrary complete simplicial fans.

For a simplicial fan $\F$ in $\R^d$, construct the graph $G_\F$ with
vertices corresponding to $d$-dimensional cones and edges corresponding to
$(d-1)$-dimensional cones of $\F$, where two vertices of $G_\F$ are
connected by an edge whenever the corresponding cones share a
$(d-1)$-dimensional face.  Pick a vector $g\in \R^d$ that does not belong to
any $(d-1)$-dimensional face of $\F$ and orient edges of $G_\F$, as
follows.  Orient an edge $\{\sigma_1,\sigma_2\}$ corresponding to two cones
$\sigma_1$ and $\sigma_2$
in $\F$ as $(\sigma_1,\sigma_2)$ if the vector $g$ points from
$\sigma_1$ to $\sigma_2$ (in a small neighborhood of the common face of
these cones).

\begin{proposition}
\label{prop:h_combin_int}
For a simplicial fan $\F$,
the $i$th entry $h_i(\F)$ of its $h$-vector equals
the number of vertices with outdegree $i$
in the oriented graph $G_\F$.
These numbers satisfy the Dehn-Sommerville symmetry:
$h_i(\F) = h_{d-i}(\F)$.
\end{proposition}

\begin{corollary}
\label{cor:h_combin_int}
{\rm (cf.~\cite[\S8.2]{Ziegler})} \
Let $P\in \R^n$ be a simple polytope.  Pick a generic linear form
$\lambda\in(\R^n)^*$.  Let $G_P$ be the $1$-skeleton of $P$ with edges
directed so that $\lambda$ increases on each edge.
Then $h_i(P)$ is the number of vertices in $G_P$ of outdegree $i$.
\end{corollary}

\begin{proof}[Proof of Proposition~\ref{prop:h_combin_int}]
The graph $G_\F$ has a unique vertex of outdegree 0.  Indeed, this is the
vertex corresponding to the cone in $\F$ containing the vector $g$.  For
any face $F$ of $\F$ (of an arbitrary dimension), let $G_\F(F)$ be the
induced subgraph on the set of $d$-dimensional cones of
$\F$ containing $F$ as a face.  Then
$G_\F(F) \simeq G_{\F'}$, where $\F'$ is the link
of the face $F$ in the fan
$\F$,  which is also a simplicial fan of smaller dimension.  Thus the
subgraph $G_\F(F)$ also contains a unique vertex of outdegree $0$ (in this
subgraph).

There is a surjective map $\phi: F \mapsto \sigma$ from all faces of
$\F$ to vertices of $G_\F$ (i.e., $d$-dimensional faces of $\F$)
that sends a face $F$ to the vertex $\sigma$ of outdegree $0$ in the subgraph
$G_\F(F)$.  Now, for a vertex $\sigma$ of $G_\F$ of outdegree $i$, the
preimage $\phi^{-1}(\sigma)$ contains exactly $\binom{d-i}{d-j}$
faces of dimension $j$.  Indeed, $\phi^{-1}(\sigma)$ is formed by
taking all possible intersections of $\sigma$ with some subset
of its $(d-1)$-dimensional faces $\{F_1,\dots,F_{d-i}\}$
on which the vector $g$ is directed towards the interior of $\sigma$;
there are exactly $d-i$ such faces because $\sigma$ has indegree $i$
in $G_\F$.  Thus a face of dimension $j$ in $\phi^{-1}(\sigma)$
has the form $F_{i_1}\cap \cdots \cap F_{i_{d-j}}$ for a $(d-j)$-element subset
$\{i_1,\dots,i_{d-j}\}\subseteq [d-i]$.

Let $\tilde h_i$ be the number of vertices of $G_\F$ of outdegree $i$.
Counting $j$-dimensional faces in preimages $\phi^{-1}(\sigma)$
one obtains the relation
$f_j(\F) =\sum_i \binom{d-i}{d-j}\,\tilde h_i$.
Comparing this with the definition of $h_i(\F)$, one deduces
that $h_i(\F) = \tilde h_{d-i}$.

Note that the numbers $h_i(\F)$ do not depend upon the choice of the vector $g$.
It follows that the numbers $\tilde h_i$ of vertices with given outdegrees also
do not depend on $g$.  Replacing the vector $g$ with $-g$ reverses the
orientation of all edges in the $d$-regular graph $G_\F$, implying the
the symmetry $\tilde h_i = \tilde h_{d-i}$.  \end{proof}

%

The Dehn-Sommerville symmetry means that $h$-polynomials are
palindromic polynomials: $t^d\, h_\F(\frac{1}{t}) = h_\F(t)$.
In this sense the $h$-vector encoding is more compact, since it is determined
by roughly half of its entries, namely $h_0,h_1,\ldots,h_{\lfloor \frac{d}{2} \rfloor}$.

Whenever possible, we will try to either give further explicit combinatorial
interpretations or generating functions for the $f$- and $h$-polynomials of
simple generalized permutohedra.

\subsection{Flag simple polytopes and $\gamma$-vectors}
\label{ssec:flag_simple_gamma}

%

A simplicial complex $\Delta$ is called a {\it flag simplicial complex\/} or
{\it clique complex\/} if it has the following property:  a collection $C$ of
vertices of $\Delta$ forms a simplex in $\Delta$ if and only if there is an
edge in the $1$-skeleton of $\Delta$ between any two vertices in $C$.
Thus flag simplicial complexes can be uniquely recovered from their
$1$-skeleta.

Let us say that a simple polytope $P$ is a {\it flag polytope\/}
if its dual simplicial complex $\Delta_P$ is a flag simplicial complex.




We next discuss $\gamma$-vectors of flag simple polytopes, as introduced by Gal \cite{Gal}
and independently in a slightly different context by Br\"anden \cite{Braenden, Braenden2};
see also the discussion in \cite[\S 1D]{Stembridge}.
A conjecture of Charney and Davis \cite{CharneyDavis}
led Gal \cite{Gal} to define the following equivalent encoding of the
$f$-vector or $h$-vector of a simple polytope $P$, in terms of smaller
integers, which are conjecturally nonnegative when $P$ is {\it flag.}
Every palindromic polynomial $h(t)$ of degree $d$ has a unique
expansion in terms of centered binomials $t^i (1+t)^{d-2i}$ for $0 \leq i \leq d/2$,
and so one can define the entries $\gamma_i=\gamma_i(P)$ of the {\it $\gamma$-vector}
$(\gamma_0, \gamma_1, \dots , \gamma_{\lfloor \frac{d}{2} \rfloor})$ and the $\gamma$-polynomial
$\gamma_P(t):=\sum_{i=0}^{\lfloor \frac{d}{2} \rfloor} \gamma_i t^i$  uniquely
by
$$
h_P(t) = \sum_{i=0}^{\lfloor \frac{d}{2} \rfloor}  \gamma_i \, t^i (1+t)^{d-2i}
       = (1+t)^d \gamma_P\left(\frac{t}{(1+t)^2}\right).
$$


\begin{conjecture}
\label{conj:Gal-conjecture}
{\rm Gal~\cite{Gal}} \
The $\gamma$-vector has nonnegative entries for any flag simple polytope.
More generally, the nonnegativity of the $\gamma$-vector holds
for every flag simplicial homology sphere.
\end{conjecture}


Thus we will try to give explicit combinatorial
interpretations, where possible, for the $\gamma$-vectors
of flag simple generalized permutohedra.
As will be seen in Section~\ref{ssec:when-is-flag}, any graph-associahedron
is a flag simple polytope.

\begin{remark} \
\label{Boolean-classes-remark}
Section~\ref{sec:gamma_nestohedra} later will employ a
a certain combinatorial approach to Gal's conjecture and $\gamma$-vector
nonnegativity that goes back to work of Shapiro, Woan, and Getu \cite{ShapiroWoanGetu},
also used by Foata and Strehl, and more recently by Br\"anden; see
\cite{Braenden2} for a thorough discussion.

Suppose $P$ is a simple polytope and one has a combinatorial formula for the $h$-polynomial
$h_P(t) = \sum_{a\in A} t^{f(a)}$, where $f(a)$ is some
statistic on the set $A$.  Suppose further that
one has a partition of $A$ into {\it $f$-symmetric Boolean classes},
i.e.\ such that the $f$-generating
function for each class is $t^r(1+t)^{2n-r}$ for some $r$.
Let $\widehat A \subset A$ be
the set of representatives of the classes where $f(a)$
takes its minimal value.  Then the $\gamma$-polynomial equals
$\gamma_P(t) = \sum_{a\in \widehat A} t^{f(a)}$.

Call $f(a)$ a ``generalized descent-statistic.''
Additionally, define 
$$
\peak(a)=\min\{f(b) \mid a \text{ and
}b\text{ in the same class}\}+1
$$ 
and call it a ``generalized
peak statistic.''  The reason for this terminology will become
apparent in Section~\ref{sec:gamma_nestohedra}.
\end{remark}

\section{Generalized permutohedra and the cone-preposet dictionary}
\label{sec:cones-fans-preposets}

This section reviews the definition of generalized permutohedra
from~\cite{Post}.  It then records some observations about the relation between
cones and fans coming from the braid arrangement and preposets.  (Normal fans
of generalized permutohedra are examples of such fans.) This leads to a
characterization for when generalized permutohedra are simple, an
interpretation for their $h$-vector in this situation, and a corollary about
when the associated toric variety is smooth.

The material in this section and in the Appendix
(Section~\ref{sec:deformations}) has substantial overlap with recent work on
{\it rank tests of non-parametric statistics}~\cite{MPSSW}.  We have tried to
indicate in places the corresponding terminology used in that paper.

\subsection{Generalized permutohedra}

Recall that a {\it usual permutohedron} in $\R^n$ is the
convex hull of $n!$ points obtained by permuting the
coordinates of any vector $(a_1,\dots,a_n)$ with strictly
increasing coordinates $a_1< \cdots < a_n$.
So the vertices of a usual permutohedron can be
labelled $v_w=(a_{w^{-1}(1)},\dots,a_{w^{-1}(n)})$ by the permutations $w$
in the symmetric group $\Sym_n$.
The edges of this permutohedron are $[v_w, v_{ws_i}]$,
where $s_i =(i,i+1)$ is an adjacent transposition.
Then, for any $w \in \Sym_n$ and any $s_i$, one has
\begin{equation}
\label{eq:genperm-defn-eqn}
v_w - v_{ws_i} = k_{w,i} (e_{w(i)}-e_{w(i+1)})
\end{equation}
where the $k_{w,i}$ are some {\it strictly positive} real scalars, and
$e_1,\ldots,e_n$ are the standard basis vectors in $\R^n$.

Note that a usual permutohedron in $\R^n$ has dimension
$d=n-1$,  because it is contained in an affine hyperplane
where the sum of coordinates $x_1+\cdots + x_n$ is constant.

\begin{definition}
\label{def:generalized_permutohedron}
\cite[Definition 6.1]{Post} \
A {\it generalized permutohedron} $P$ is the convex hull of $n!$ points $v_w$
in $\R^n$ labelled by the permutations $w$ in the symmetric group $\Sym_n$,
such that for any $w \in \Sym_n$ and any adjacent transposition $s_i$, one
still has equation~\eqref{eq:genperm-defn-eqn}, but with $k_{w,i}$ assumed only
to be nonnegative, that is, $k_{w,i}$ can vanish.
\end{definition}

The Appendix shows that all $n!$ points $v_w$ in a generalized
permutohedron $P$ are actually vertices of $P$ (possibly with repetitions); see
Theorem~\ref{th:deformation}.  Thus a generalized permutohedron $P$ comes
naturally equipped with the surjective map $\Psi_P:\Sym_n\to \Vertices(P)$
given by $\Psi_P:w\mapsto v_w$, for $w\in\Sym_m$.

Definition~\ref{def:generalized_permutohedron} says that a generalized
permutohedron is obtained by moving the vertices of the usual permutohedron in
such a way that directions of edges are preserved, but some edges (and higher
dimensional faces) may degenerate.  In the Appendix such deformations of a
simple polytope are shown to be equivalent to various other notions of
deformation; see
Proposition~\ref{prop:permutohedron-summand-characterization} below
and the more general Theorem~\ref{th:deformation}.

\subsection{Braid arrangement}

Let $x_1,\dots,x_n$ be the usual coordinates in $\R^n$.
Let $\R^n/(1,\dots,1)\R\simeq \R^{n-1}$ denote the quotient space
modulo the 1-dimensional subspace generated by the vector $(1,\dots,1)$.
The {\it braid
arrangement\/} is the arrangement of hyperplanes $\{x_i-x_j = 0\}_{1\leq i < j
\leq n}$ in the space $\R^n/(1,\dots,1)\R$.
These hyperplanes subdivide the space into the polyhedral cones
$$
C_w :=
\{x_{w(1)}\leq x_{w(2)} \leq \cdots \leq x_{w(n)}\}
$$
labelled by permutations $w\in\Sym_n$, called {\it Weyl chambers\/} (of type $A$).
The Weyl chambers and their lower dimensional faces form a complete simplicial
fan, sometimes called the {\it braid arrangement fan.}

Note that a usual permutohedron $P$ has dimension $d=n-1$, so one can reduce
its normal fan modulo the 1-dimensional subspace $P^\perp = (1,\dots,1)\R$.
The braid arrangement fan is exactly the (reduced) normal fan $\N(P)/P^\perp$
for a usual permutohedron $P\subset \R^n$.
Indeed, the (reduced) normal cone $\N_{v_w}(P)/P^\perp$ of
$P$ at  vertex $v_w$ is exactly the Weyl chamber $C_w$.
(Here one identifies $\R^n$ with $(\R^n)^*$ via the standard inner product.)

Recall that the {\it Minkowski sum\/} $P+Q$ of two polytopes $P, Q\subset \R^n$
is the polytope $P+Q:=\{x+y \mid x\in P, y\in Q \}$.  Say that $P$ is a {\it
Minkowski summand\/} of $R$, if there is a polytope $Q$ such that $P+Q=R$.  Say
that a fan $\F$ is {\it refined\/} by a fan $\F'$
if any cone in $\F$ is a union of cones in $\F'$.
The following proposition is a special case of
Theorem~\ref{th:deformation}.

\begin{proposition}
\label{prop:permutohedron-summand-characterization}
A polytope $P$ in $\R^n$ is a generalized permutohedron if and only if
its normal fan (reduced by $(1,\dots,1)\R$) is refined by the braid
arrangement fan.

Also, generalized permutohedra are exactly the polytopes
arising as Minkowski summands of usual permutohedra.
\end{proposition}

This proposition shows that generalized permutohedra
lead to the study of cones given by some inequalities of the form
$x_i-x_j\geq 0$ and fans formed by such cones.
Such cones are naturally related to posets and preposets.

\subsection{Preposets, equivalence relations, and posets}

Recall that a {\it binary relation} $R$ on a set $X$ is a subset of
$R\subseteq X \times X$.  A {\it preposet\/} is a reflexive and transitive
binary relation $R$, that is $(x,x)\in R$ for all $x\in X$, and  whenever
$(x,y), (y,z) \in R$ one has $(x,z)\in R$.  In this case we will often use the
notation $x \preceq_R y$ instead of $(x,y) \in R$.
Let us also write $x\prec_R y$ whenever $x\preceq_R y$ and $x\neq y$.

An {\it equivalence relation} $\equiv$ is the special case of a preposet whose binary relation is
also symmetric.  Every preposet $Q$ gives rise to an equivalence relation
$\equiv_Q$ defined by $x \equiv_Q y$ if and only if both $x \preceq_Q y$ and $y \preceq_Q x$.
A {\it poset} is the special case of a preposet $Q$ whose associated
equivalence relation $\equiv_Q$ is the trivial partition, having
only singleton equivalence classes.

Every preposet $Q$ gives rise to the poset $Q/{\equiv_Q}$ on the equivalence
classes $X/{\equiv_Q}$.  A preposet $Q$ is uniquely determined by $\equiv_Q$
and $Q/{\equiv_Q}$, that is, a preposet is just an equivalence relation together
with a poset structure on the equivalence classes.

A {\it preposet} $Q$ on $X$ is {\it connected} if the undirected graph having
vertices $X$ and edges $\{x,y\}$ for all $x \preceq_Q y$ is connected.

A {\it covering relation\/} $x\lessdot_Q y$ in a poset $Q$ is a pair of elements
$x\prec_Q y$ such that there is no $z$ such that $x\prec_Q z\prec_Q y$.  The
{\it Hasse diagram\/} of a poset $Q$ on $X$ is the directed graph on $X$ with
edges $(x,y)$ for covering relations $x\lessdot_Q y$.

Let us say that a poset
$Q$ is a {\it tree-poset} if its Hasse diagram is a spanning tree on $X$.
Thus tree-posets correspond to directed trees on the vertex set $X$.

A {\it linear extension\/} of a poset $Q$ on $X$ is a linear ordering
$(y_1,\dots,y_n)$ of all elements in $X$ such that $y_1\prec_Q y_2\prec_Q \cdots
\prec_Q y_n$.  Let $\LL(Q)$ denote the set of all linear extensions of $Q$.

The {\it union\/} $R_1\cup R_2$ of two binary relations $R_1,R_2$ on $X$ is
just their union as two subsets of $X\times X$.  Given any reflexive binary
relation $Q$, denote by $\overline{Q}$ the preposet which is its transitive
closure.  Note that if $Q_1$ and $Q_2$ are two preposets on the same set $X$,
then the binary relation $Q_1 \cup Q_2$ is not necessarily a preposet.
However, we can obtain a preposet by taking its
transitive closure $\overline{Q_1 \cup Q_2}$.

Let $R\subseteq Q$ denote containment of binary relations on the same set,
meaning containment as subsets of $X\times X$.  Also let $R^{op}$ be the {\it
opposite\/} binary relation, that is $(x,y)\in R^{op}$ if and only if $(y,x)\in
R$.

For two preposets $P$ and $Q$ on the same set, let us say
that $Q$ is a {\it contraction\/} of $P$ if there is a binary relation
$R \subseteq P$ such that $Q = \overline {P \cup R^{op}}$.
In other words, the equivalence classes of $\equiv_Q$ are obtained
by merging some equivalence classes of $\equiv_P$
{\it along relations\/}
in $P$ and the poset structure on equivalence classes of $\equiv_Q$
is induced from the poset structure on classes of $\equiv_P$.

For example, the preposet $1<\{2,3\}<4$ (where $\{2,3\}$ is an equivalence
class) is a contraction of the poset $P = (1<3,\, 2<3,\, 1<4,\, 2<4)$.
However, the preposet $(\{1,2\}<3,\ \{1,2\}<4)$ is not a contraction of $P$
because $1$ and $2$ are incomparable in $P$.

\begin{definition}
\label{def:poset_fan}
Let us say that two preposets $Q_1$ and $Q_2$ on the same set
{\it intersect properly\/}
if the preposet $\overline{Q_1 \cup Q_2}$ is both a contraction of
$Q_1$ and of $Q_2$.


A {\it complete fan of posets\footnote{In \cite{MPSSW}, this is called a {\it
convex rank test.}} on $X$} is a collection of distinct posets on $X$ which
pairwise intersect properly, and whose linear extensions (disjointly) cover all
linear orders on $X$.
\end{definition}

Compare Definition~\ref{def:poset_fan}
to the definitions of properly intersecting cones and complete fan of cones;
see Section~\ref{ssec:cones_fans}.
This connection will be
elucidated in Proposition~\ref{prop:cone-preposet-dictionary}.

\begin{example}
\label{ex:intersections}
The two posets $P_1:= 1<2$ and $P_2:= 2<1$ on the set $\{1,2\}$ intersect
properly.  Here $\overline{P_1 \cup P_2}$ is equal to
$\{1<2, 2<1\}$.  These $P_1$ and $P_2$ form a complete fan of posets.

However,
the two posets $Q_1:=2<3$ and $Q_2:=1<2<3$ on the set $\{1,2,3\}$ do not
intersect properly.  In this case $\overline{Q_1 \cup Q_2} = Q_2$,
which is not a contraction of $Q_1$.
\end{example}

\subsection{The dictionary}

Let us say that a {\it braid cone\/} is a polyhedral cone in the space
$\R^n/(1,\dots,1)\R\simeq \R^{n-1}$ given by a conjunction of inequalities
of the form $x_i - x_j \geq 0$.
In other words, braid cones are polyhedral cones formed by
unions of Weyl chambers or their lower dimensional faces.

There is an obvious bijection between preposets and braid cones.
For a preposet $Q$ on the set $[n]$, let $\sigma_Q$ be the braid
cone in the space $\R^n/(1,\dots,1)\R$ defined by the conjunction of the
inequalities $x_i \leq x_j$ for all $i \preceq_Q j$.  Conversely, one can always
reconstruct the preposet $Q$ from the cone $\sigma_Q$ by saying that
$i\preceq_Q j$ whenever $x_i \leq x_j$ for all points in $\sigma_Q$.


\begin{proposition}
\label{prop:cone-preposet-dictionary}
Let the cones $\sigma, \sigma'$ correspond to the preposets $Q, Q'$
under the above bijection.  Then
\begin{enumerate}
\item[(1)] The preposet $\overline{Q \cup Q'}$ corresponds to the cone $\sigma \cap \sigma'$.
\item[(2)] The preposet $Q$ is a contraction of $Q'$ if and only if the cone
$\sigma$ is a face $\sigma'$.
\item[(3)] The preposets $Q, Q'$ intersect properly if and only if
the cones $\sigma, \sigma'$ do.
\item[(4)] $Q$ is a poset if and only if
$\sigma$ is a full-dimensional cone, i.e., $\dim\sigma = n-1$.
\item[(5)] The equivalence relation $\equiv_Q$ corresponds to the linear
span $\Span(\sigma)$ of $\sigma$.
\item[(6)] The poset $Q/{\equiv_Q}$
corresponds to a full-dimensional cone inside $\Span(\sigma_Q)$.
\item[(7)] The preposet $Q$ is connected if and only if the cone $\sigma$
is pointed.

\item[(8)] If $Q$ is a poset, then the minimal set of inequalities
describing the cone $\sigma$ is $\{x_i \leq x_j \mid i \lessdot_Q j\}$.
(These inequalities associated with covering relations in $Q$ are exactly
the facet inequalities for $\sigma$.)

\item[(9)] $Q$ is a tree-poset if and only if $\sigma$ is a
full-dimensional simplicial cone.
\item[(10)] For $w\in \Sym_n$, the cone $\sigma$ contains the
Weyl chamber $C_w$ if and only if $Q$ is a poset and $w$
is its linear extension, that is $w(1)\prec_Q w(2)\prec_Q \cdots \prec_Q w(n)$.
\end{enumerate}
\end{proposition}

\begin{proof}
(1)  The cone $\sigma\cap \sigma'$ is given by conjunction of all inequalities
for $\sigma$ and $\sigma'$.  The corresponding preposet is obtained by adding
all inequalities that follow from these, i.e.,
by taking the transitive closure of $Q\cup Q'$.

(2) Faces of $\sigma'$ are obtained by replacing some inequalities
$x_i \leq x_j$ defining $\sigma'$ with equalities $x_i = x_j$, or
equivalently, by adding the opposite inequalities $x_i \geq x_j$.

(3) follows from (1) and (2).

(4) $\sigma$ is full-dimensional if its defining relations do not include any
equalities $x_i = x_j$, that is $\equiv_Q$ has only singleton
equivalence classes.

(5)  The cone associated with the equivalence relation $\equiv_Q$
is given by the equations $x_i = x_j$ for $i\equiv_Q j$, which is exactly
$\Span(\sigma)$.

(6) Follows from (4) and (5).

(7) The maximal subspace contained in the half-space $\{x_i \leq x_j\}$ is
given by $x_i  = x_j$.  Thus the maximal subspace contained in
the cone $\sigma$ is given by the conjunction of equations
$x_i = x_j$ for $i\leq_Q j$.
If $Q$ is disconnected then this subspace has a positive dimension.
If $Q$ is connected then this subspace is given by $x_1 = \cdots  = x_n$,
which is just the origin in the space $\R^n/(1,\dots,1)\R$.

(8) The inequalities for the covering relations $i\lessdot_Q j$ imply all other
inequalities for $\sigma$ and they cannot be reduced to a smaller
set of inequalities.

(9)  By (4) and (7) full-dimensional pointed cones correspond
to connected posets.  These cones will be simplicial if they are given by
exactly $n-1$ inequalities.  By~(8) this means that the corresponding poset
should have exactly $n-1$ covering relations, i.e., it is a tree-poset.

(10)  Follows from (4) and definitions.
\end{proof}

According to Proposition~\ref{prop:cone-preposet-dictionary},
a full-dimensional braid cone $\sigma$
associated with a poset $Q$ can be described in three different
ways (via all relations in $Q$, via covering relations in $Q$,
and via linear extensions $\LL(Q)$ of $Q$) as
$$
\sigma = \{x_i\leq x_j \mid i\preceq_{Q}j\}
=\{x_i\leq x_j \mid i\lessdot_{Q}j\}
=\bigcup_{w\in\LL(Q)} C_w.
$$

Let $\F$ be a family of $d$-cones in $\R^d$ which intersect properly.
Since they have disjoint interiors, they will correspond to a complete fan if
and only if their closures cover $\R^d$, or equivalently, their spherical
volumes sum to the volume of the full $(d-1)$-sphere.

A braid cone corresponding to a poset $Q$ is the union
of the Weyl chambers $C_w$ for all linear extensions $w \in \LL(Q)$,
and every Weyl chamber has the same spherical volume
($\frac{1}{n!}$ of the sphere)
due to the transitive Weyl group action.
Therefore, a collection of properly intersecting posets
$\{Q_1,\ldots,Q_t\}$ on $[n]$ correspond to a complete fan on braid
cones if and only if
$$
\bigcup_{i=1}^t \LL(Q_i) = \Sym_n
\textrm{ (disjoint union), or equivalently, if and only if }
\sum_{i=1}^t |\LL(Q_i)| = n!,
$$
cf.~ Definition~\ref{def:poset_fan}.

\begin{corollary}
\label{cor:fans=fans}
A complete fan of braid cones
(resp., pointed braid cones, simplicial braid cones) in $\R^n/(1,\dots,1)\R$
corresponds to
a complete fan of posets (resp., connected posets, tree-posets)
on $[n]$.
\end{corollary}


Using Proposition~\ref{prop:permutohedron-summand-characterization}, we can
relate Proposition~\ref{prop:cone-preposet-dictionary} and
Corollary~\ref{cor:fans=fans} to generalized permutohedra.
Indeed, normal cones of a generalized permutohedron
(reduced modulo $(1,\dots,1)\R$) are braid cones.

For a generalized permutohedron $P$,
define the {\it vertex poset\/} $Q_v$ at a vertex $v\in\Vertices(P)$
as the poset on $[n]$ associated with the
normal cone $\N_v(P)/(1,\dots,1)\R$ at
the vertex $v$, as above.

\begin{corollary}
\label{cor:fans_posets}
For a generalized permutohedron (resp., simple generalized
permutohedron) $P$, the collection
of vertex posets $\{Q_v\mid v\in \Vertices(P)\}$ is a
complete fan of posets
(resp., tree-posets).
\end{corollary}

Thus normal fans of generalized permutohedra correspond to
certain complete fans of posets, which we call {\it polytopal.}
In~\cite{MPSSW}, the authors call such fans {\it submodular
rank tests,} since they are in bijection with the faces of the cone of
submodular functions.  That cone is precisely the deformation cone we discuss
in the Appendix.

\begin{example}
\label{ex:simplicial-non-polytopal-example}
In \cite{MPSSW}, the authors modify an example of Studen\'y \cite{Studeny}
to exhibit a non-polytopal complete fan of posets.  They also kindly
provided us with the following further nonpolytopal example, having
$16$ posets $Q_v$, all of them tree-posets:
$(1,2<3<4)$ (which means that $1<3$ and $2<3$),
$(1,2<4<3)$,
$(3,4<1<2)$,
$(3,4<2<1)$,
$(1<4<2,3)$,
$(4<1<2,3)$,
$(2<3<1,4)$,
$(3<2<1,4)$,
$(1<3<2<4)$,
$(1<3<4<2)$,
$(3<1<2<4)$,
$(3<1<4<2)$,
$(2<4<1<3)$,
$(2<4<3<1)$,
$(4<2<1<3)$,
$(4<2<3<1)$.
This gives a complete fan of simplicial cones, but does {\it not}
correspond to a (simple) generalized permutohedron.
\end{example}

Recall that $\Psi_P:\Sym_n\to\Vertices(P)$ is the surjective map
$\Psi_P:w\mapsto v_w$; see Definition~\ref{def:generalized_permutohedron}.
The previous discussion immediately implies the following corollary.

\begin{corollary}
\label{cor:Map}
Let $P$ be a generalized permutohedron in $\R^n$,
and $v\in\Vertices(P)$ be its vertex.
For $w\in\Sym_n$, one has $\Psi_P(w) = v$ whenever the normal cone $\N_v(P)$
contains the Weyl chamber $C_w$.
The preimage $\Psi_P^{-1}(v)\subseteq\Sym_n$ of a vertex
$v\in\Vertices(P)$ is
the set $\LL(Q_v)$ of all linear extensions of the vertex poset $Q_v$.
\end{corollary}

We remark on the significance of this cone-preposet dictionary for toric varieties associated to
generalized permutohedra or their normal fans; see Fulton \cite{Fulton} for
further background.

A complete fan $\F$ of polyhedral cones in $\R^d$ whose cones are rational with
respect to $\Z^d$ gives rise to a toric variety $X_\F$, which is normal,
complete and of complex dimension $d$.

This toric variety is {\it projective} if and only if $\F$ is the normal fan
$\N(P)$ for some polytope $P$, in which case one also denotes $X_\F$ by
$X_P$.

The toric variety $X_\F$ is {\it quasi-smooth} or {\it orbifold} if and only if
$\F$ is a complete fan of simplicial cones;  in the projective case, where
$\F=\N(P)$, this corresponds to $P$ being a simple polytope.

In this situation, the $h$-numbers of $\F$ (or of $P$) have the auxiliary geometric meaning
as the (singular cohomology) Betti numbers $h_i = \dim H^i(X_{\F},\C)$.
The symmetry $h_i=h_{d-i}$ reflects Poincar\'e duality for this quasi-smooth variety.

The toric variety $X_\F$ is {\it smooth}
exactly when every top-dimensional cone of $\F$ is
not only simplicial but {\it unimodular}, that is, the primitive vectors on its extreme rays
form a $\Z$-basis for $\Z^d$.  Equivalently, the facet inequalities $\ell_1,\ldots,\ell_d$
can be chosen to form a $\Z$-basis for $(\Z^d)^*=\Hom(\Z^d,\Z)$ inside $(\R^d)^*$.
One has $X_\F$ both smooth and projective if and only if  $\F=\N(P)$ for
a {\it Delzant} polytope $P$, that is, one which is simple and has every vertex normal
cone unimodular.

\begin{corollary} {\rm (cf.~ \cite[\S5]{Zel})} \
A complete fan $\F$ of posets gives rise to a complete toric variety $X_\F$,
which will be projective if and only if $\F$ is associated
with the normal fan $\N(P)$ for a generalized permutohedron.

A complete fan $\F$ of tree-posets gives rise to a (smooth, not just orbifold)
toric variety $X_\F$, which will be projective if and only $\F$ is associated
with the normal fan $\N(P)$ of a simple generalized permutohedron.  In other
words, simple generalized permutohedra are always  Delzant.
\end{corollary}

\begin{proof}
All the assertions should be clear from the above discussion, except
for the last one about simple generalized permutohedra being Delzant.
However, a tree-poset $Q$ corresponds to a set of functionals $x_i-x_j$ for
the edges $\{i,j\}$ of a tree, which are well-known to give a $\Z$-basis for
$(\Z^d)^*$, cf.~ \cite[Proposition 7.10]{Post}.
\end{proof}

\section{Simple generalized permutohedra}
\label{sec:simples}

\subsection{Descents of tree-posets and $h$-vectors}
\label{ssec:descents_of_tree_posets}

  The goal of this section is to combinatorially interpret the $h$-vector of any
simple generalized permutohedron.

\begin{definition}
\label{def:des_tree_poset}
Given a poset $Q$ on $[n]$, define the {\it descent set\/}
$\Des(Q)$ to be the set of ordered pairs $(i,j)$ for which $i \lessdot_Q j$
is a covering relation in $Q$ with $i >_\Z j$, and define the statistic
{\it number of descents\/} $\des(Q):=|\Des(Q)|$.
\end{definition}

\begin{theorem}
\label{th:shelling-and-descents}
Let $P$ be a simple generalized permutohedron, with
vertex posets $\{Q_v\}_{v \in \Vertices(P)}$.  Then one has
the following expression for its $h$-polynomial:
\begin{equation}
\label{eq:h-polynomial-is-descent-polynomial}
h_P(t) = \sum_{v \in \Vertices(P)} t^{\des(Q_v)}.
\end{equation}

More generally, for a complete fan
$\F = \{Q_v\}$ of tree-posets
(see Definition~\ref{def:poset_fan}), one also has
$h_\F(t) = \sum_{v} t^{\des(Q_v)}$.
\end{theorem}


\begin{proof}
(cf.\ proof of Proposition 7.10 in \cite{Post})
Let us prove the more general claim about fans of tree-posets,
that is, simplicial fans coarsening the braid arrangement fan.

Pick a generic vector $g=(g_1,\dots,g_n)\in\R^n$
such that $g_1<\cdots < g_n$ and construct the directed graph $G_\F$,
as in Proposition~\ref{prop:h_combin_int}.
Let $\sigma = \{x_i \leq x_j\mid i\lessdot_{Q_v} j\}$ be the cone of $\F$
associated with poset $Q_v$,
see Proposition~\ref{prop:cone-preposet-dictionary}(8).
Let $\sigma'$ be an adjacent cone separated from $\sigma$
by the facet $x_i = x_j$, $i\lessdot_{Q_v} j$.
The vector $g$ points from $\sigma$ to $\sigma'$ if and only if
$g_i >_{\R} g_j$, or equivalently, $i>_{\Z}j$.
Thus the outdegree of $\sigma$ in the graph $G_\F$ is exactly
the descent number $\des(Q)$.  The claim now follows
from Proposition~\ref{prop:h_combin_int}.
\end{proof}

%
%

For a usual permutohedron $P$ in $\R^n$, the vertex posets $Q_v$ are just
all linear orders on $[n]$.  So its $h$-polynomial $h_P(t)$ is the
classical {\it Eulerian polynomial\/}\footnote{Note that a more
standard convention is to call  $t A_n(t)$ the Eulerian polynomial.}
\begin{equation}
\label{eq:usual-Eulerian-polynomial}
A_n(t):=\sum_{w \in \Sym_n} t^{\des(w)},
\end{equation}
where $\des(w):=\#\{i\mid w(i)>w(i+1)\}$ is the usual
{\it descent number\/} of a permutation $w$.

Any element $w$ in the Weyl group $\Sym_n$ sends a complete fan
$\F=\{Q_i\}$ of tree-posets to another such complete fan
$w\F=\{wQ_i\}$, by relabelling all of the posets.  Since $w\F$ is an
isomorphic simplicial complex, with the same $h$-vector, this leads to
a curious corollary.

\begin{definition}
Given a tree-poset $Q$ on $[n]$, define its {\it generalized Eulerian polynomial}
$$
A_Q(t):=\sum_{w \in \Sym_n} t^{\des(wQ)}.
$$
Note that $A_Q$ depends upon $Q$ only as an {\it unlabelled} poset.
\end{definition}

When $Q$ is a linear order, $A_Q(t)$ is the
usual Eulerian polynomial $A_n(t)$.

\begin{corollary}
\label{cor:h-as-average}
The $h$-polynomial $h_P(t)$ of a simple generalized permutohedron $P$ is the
``average'' of the generalized Eulerian polynomials of its
vertex tree-posets $Q_v$:
$$
h_P(t) = \frac{1}{n!} \sum_{v \in \Vertices(P)} A_{Q_v}(t).
$$
\end{corollary}

See Example~\ref{ex:zonotopal-example} below for an illustration of
Theorem~\ref{th:shelling-and-descents} and Corollary~\ref{cor:h-as-average}.

\subsection{Bounds on the $h$-vector and monotonicity}

 It is natural to ask for upper and lower bounds on the $h$-vectors of
simple generalized permutohedra.  Some of these follow immediately from an $h$-vector
monotonicity result of Stanley that applies to complete simplicial
fans.  To state it, we recall a definition from that paper.

\begin{definition}
Say that a simplicial complex $\Delta'$ is a {\it geometric subdivision\/} of a
simplicial complex $\Delta$ if they have geometric realizations which are
topological spaces on the same underlying set, and every face of $\Delta'$ is
contained in a single face of $\Delta$.
\end{definition}

\begin{theorem}
\label{th:subdivision-monotonicity}
{\rm (see \cite[Theorem 4.1]{Stanley-localh})} \
If $\Delta'$ is a geometric subdivision of a Cohen-Macaulay simplicial complex
$\Delta$, then the h-vector of $\Delta'$ is componentwise weakly larger than
that of $\Delta$.

In particular this holds when $\Delta, \Delta'$ come from two complete
simplicial fans and $\Delta'$ refines $\Delta$, e.g., the normal fans of  two
simple polytopes $P, P'$ in which $P$ is a Minkowski summand of $P$.
\end{theorem}

\begin{corollary}
\label{cor:upperbound}
A simple generalized permutohedron $P$ in $\R^n$
has $h$-polynomial coefficientwise
smaller than that of the permutohedron,
namely the Eulerian polynomial $A_n(t)$.
\end{corollary}

\begin{proof}
Proposition~\ref{prop:permutohedron-summand-characterization} tells us that
the normal fan of $P$ is refined by that of the permutohedron, so the
above theorem applies.
\end{proof}

\begin{remark}
  Does the permutohedron also provide an upper bound for
the $f$-vectors, flag $f$- and flag $h$-vectors, generalized $h$-vectors, and $cd$-indices of generalized
permutohedra also in the non-simple case?  Is there also a monotonicity result for
these other forms of face and flag number data when one has two generalized
permutohedra $P, P'$ in which $P$ is a Minkowski summand of $P'$?

  The answer is ``Yes'' for $f$-vectors and flag $f$-vectors, which clearly
increase under subdivision.  The answer is also ``Yes'' for
generalized $h$-vectors, which Stanley also showed \cite[Corollary 7.11]{Stanley-localh}
can only increase under geometric subdivisions of rational convex polytopes.
But for flag $h$-vectors and $cd$-indices, this is not so clear.
\end{remark}

Later on (Example~\ref{ex:simplex-example},
Section~\ref{sec:Stanley-Pitman-example}, and
Section~\ref{sec:flossing}) there will be more to say about lower bounds
for $h$-vectors of simple generalized permutohedra within various classes.


\section{The case of zonotopes}
\label{sec:zonotopes}

  This section illustrates some of the foregoing results in
the case where the simple generalized permutohedron is a zonotope;
see also \cite[\S8.6]{Post}.
Zonotopal generalized permutohedra are exactly the graphic zonotopes,
and the simple zonotopes among them correspond to a very restrictive class
of graphs that are easily dealt with.

  A {\it zonotope} is a convex polytope $Z$ which is the Minkowski sum of one-dimensional polytopes
(line segments), or equivalently, a polytope whose normal fan $\N(Z)$ coincides with chambers and cones of a
hyperplane arrangement.  Under this equivalence, the line segments which are the Minkowski summands of
$Z$ lie in the directions of the normal vectors to the hyperplanes in the arrangement.
Given a graph $G=(V,E)$ without loops or multiple edges, on node set $V=[n]$ and with edge set $E$,
define the associated {\it graphic zonotope} $Z_G$ to be the
Minkowski sum of line segments in the
directions $\{e_i - e_j\}_{ij \in E}$.

Proposition~\ref{prop:permutohedron-summand-characterization} then immediately implies the following.

\begin{proposition}
The zonotopal generalized permutohedra are exactly the graphic zonotopes $Z_G$.
\end{proposition}

 Simple zonotopes are very special among all zonotopes, and simple graphic
zonotopes have been observed \cite[Remark 5.2]{Kim} to correspond to a very
restrictive class of graphic zonotopes, namely those whose biconnected
components are all complete graphs.

%

Recall that for a graph $G=(V,E)$, there is an equivalence relation on $E$
defined by saying $e \sim e'$ if there is some circuit (i.e., cycle which is
minimal with respect to inclusion of edges) of $G$ containing both $e, e'$.
The $\sim$-equivalence classes are then called {\it biconnected components} of
$G$.

\begin{proposition}
\label{prop:simple-graphic-zonotopes}
{\rm \cite[Remark 5.2]{Kim}} \
The graphic zonotope $Z_G$ corresponding to a graph $G=(V,E)$ is a simple polytope
if and only if
every biconnected component of $G$ is the set of edges of a complete
subgraph some subset of the vertices $V$.

In this case, if $V_1,\ldots,V_r \subseteq V$ are the node sets for these
complete subgraphs, then $Z_G$ is isomorphic to the Cartesian product of
usual permutohedra of dimensions $|V_j|-1$ for $j=1,2,\dots,r$.
\end{proposition}

Let us give another description for this class of graphs.  For a graph
$F$ with $n$ edges $e_1,\dots,e_n$,  the {\it line graph\/} $\Line(F)$ of
$F$ is the graph on the vertex set $[n]$ where $\{i,j\}$ is an edge in
$\Line(F)$ if and
only if the edges $e_i$ and $e_j$ of $F$ have a common vertex.
The following claim is left an exercise for the reader.

\begin{exercise}
For a graph $G$, all biconnected components of $G$ are edge sets
of complete graphs if and only if $G$ is isomorphic to the line graph
$\Line(F)$ of some forest $F$.
Biconnected components of $\Line(F)$ correspond to non-leaf vertices of $F$.
\end{exercise}

For the sake of completeness, included here is a proof
of Proposition~\ref{prop:simple-graphic-zonotopes}.

\begin{proof}[Proof of Proposition~\ref{prop:simple-graphic-zonotopes}]
If the biconnected components of $G$ induce subgraphs isomorphic to graphs
$G_1,\dots,G_r$ then one can easily check that $Z_G$ is the Cartesian product
of the zonotopes $Z_{G_i}$.  Since a Cartesian product of polytopes is simple
if and only if each factor is simple, this reduces to the case where $r=1$.
Also note that when $r=1$ and $G$ is a complete graph, then $Z_G$ is the
permutohedron, which is well-known to be simple.

For the reverse implication, assume $G$ is biconnected but not a complete
graph, and it will suffice, by
Proposition~\ref{prop:cone-preposet-dictionary}(9), to construct a vertex $v$
of $Z_G$ whose poset $Q_v$ is not a tree-poset.  One uses the fact
\cite{GreeneZaslavsky} that a vertex $v$ in the graphic zonotope $Z_G$
corresponds to an acyclic orientation of $G$, and the associated poset $Q_v$ on
$V$ is simply the transitive closure of this orientation.  Thus it suffices to
produce an acyclic orientation of $G$ whose transitive closure has Hasse
diagram which is not a tree.

Since $G$ is biconnected but not complete, there must be two vertices $\{x,y\}$ that do not span
an edge in $E$, but which lie in some circuit $C$.  Traverse this circuit $C$ in some cyclic order,
starting at the node $x$,
passing through some nonempty set of vertices $V_1$ before passing through $y$, and then through
a nonempty set of vertices $V_2$ before returning to $x$.
One can then choose arbitrarily a total order on the node set $V$ so that these sets appear
as segments in this order:
$$
V_1, \quad x, \quad y, \quad V_2, \quad V-(V_1 \cup V_2 \cup \{x,y\}).
$$
It is then easily checked that if one orients the edges of $G$ consistently with this total order,
then the associated poset has a non-tree Hasse diagram:
for any $v_1 \in V_1$ and $v_2 \in V_2$, one has $v_1 < x,y < v_2$ with $x,y$ incomparable.
\end{proof}

\begin{corollary}
Let $Z_G$ be a simple graphic zonotope, with notation as in
Proposition~\ref{prop:simple-graphic-zonotopes}.

Then $Z_G$ is flag, and its $f$-polynomial, $h$-polynomials, $\gamma$-polynomials
are all equal to products for $j=1,2,\ldots,r$ of the $f$-, $h$-, or $\gamma$-polynomials
of $(|V_j|-1)$-dimensional permutohedra.
\end{corollary}

\begin{proof}
Use Proposition~\ref{prop:simple-graphic-zonotopes} along with the fact that
a Cartesian product of simple polytopes is flag if and only if each factor is
flag, and has $f$-, $h$- and $\gamma$-polynomial equal to the product of
the same polynomials for each factor.
\end{proof}

Note that the $h$-polynomial for an $(n-1)$-dimensional permutohedron is the
Eulerian polynomial $A_n(t)$ described in \eqref{eq:usual-Eulerian-polynomial}
above, and the $\gamma$-polynomial is given explicitly in Theorem~\ref{th:CDG} below.

\begin{example}
\label{ex:zonotopal-example}
Consider the graph $G=(V,E)$ with $V=[4]:=\{1,2,3,4\}$ and $E=\{12,13,23,14\}$,
whose biconnected components are the triangle $123$ and the edge $14$, which are
both complete subgraphs on node sets $V_1=\{1,2,3\}$ and $V_2=\{1,4\}$.
Hence the associated graphic zonotope $Z_G$ is simple and flag, equal to
the Cartesian product of a hexagon with a line segment, that is, $Z_G$ is a hexagonal
prism.

Its $f$-, $h$- and $\gamma$-polynomials are
$$
\begin{array}{lll}
f_{Z_G}(t)&=(2+t)(6+6t+t^2)                  &=12+18t+8t^2+t^3\\
h_{Z_G}(t)& =A_2(t) A_3(t)  =(1+t)(1+4t+t^2) &= 1+5t+5t^2+t^3\\
\gamma_{Z_G}(t)&=(1)(1+2t)                  &=1+2t.
\end{array}
$$

One can arrive at the same $h$-polynomial using Theorem~\ref{th:shelling-and-descents}.
One lists the tree-posets $Q_v$ for each of the $12$ vertices $v$ of the hexagonal prism $Z_G$,
coming in $5$ isomorphism types,
along with the number of descents for each:
$$
\begin{array}{clc}
\text{type} &\text{ poset }Q_v & \des \\[.1in]
\text{chain:} & 2<3<1<4      & 1 \\
&3<2<1<4                     & 2 \\
&4<1<2<3                     & 1 \\
&4<1<3<2                     & 2 \\[.1in]
\text{vee:} & 1<2<3 \text{ and } 1<4     & 0 \\
& 1<3<2 \text{ and } 1<4                & 1
\end{array} 
\qquad\quad
\begin{array}{clc}
\text{type} &\text{ poset }Q_v & \des \\[.1in]
\text{wedge:} & 2<3<1 \text{ and } 4<1    & 2 \\
& 3<2<1 \text{ and } 4<1                 & 3 \\[.1in]
\text{wye:} & 2<1<3 \text{ and } 1<4      & 1 \\
& 3<1<2 \text{ and } 1<4                 & 1 \\[.1in]
\text{lambda:} & 3<1<2 \text{ and }4<1    & 2 \\
&2<1<3 \text{ and }4<1                   & 2
\end{array}
$$
and finds that $\sum_{v} t^{\des(Q_v)} = 1+5t+5t^2+t^3$.

Lastly one can get this $h$-polynomial from Corollary~\ref{cor:h-as-average}, by
calculating directly that
$$
\begin{aligned}
A_\text{chain}(t)&=1+11t+11t^2+t^3 = A_4(t)\\
A_\text{vee}(t)&=3+10t+8t^2+3t^3\\
A_\text{wedge}(t)&=3+8t+10t^2+3t^3\\
A_\text{wye}(t)=A_\text{lambda}(t)&=2+10t+10t^2+2t^3\\
\end{aligned}
$$
and then the $h$-polynomial is
$$
\frac{1}{4!} \left[
4 A_\text{chain}(t) +
2 A_\text{vee}(t)+
2 A_\text{wedge}(t) +
2 A_\text{wye}(t) +
2 A_\text{lambda}(t)
\right]
= 1+5t+5t^2+t^3.
$$
\end{example}

\section{Building sets and nestohedra}
\label{sec:nested}

This section reviews some results from~\cite{FeichtnerSturmfels}, \cite{Post},
and~\cite{Zel} regarding the important special case of generalized permutohedra that
arise from building sets.  These generalized permutohedra turn out always to be simple. 
Their dual simplicial complexes, the nested set complexes, are defined,
and several tools are given for calculating their $f$- and $h$-vectors.
The notion of nested sets goes back to work of Fulton and MacPherson~\cite{FM},
and DeConcini and Procesi~\cite{DP} defined building sets
and nested set complexes.
However, our 
exposition mostly follows~\cite{Post} and~\cite{Zel}.

\subsection{Building sets, nestohedra, and nested set complexes}

\begin{definition}
\label{def:building-set-definition}
{\rm \cite[Definition~7.1]{Post}} \
Let us say that a collection $\B$ of nonempty subsets
of a finite set $S$ is a {\it building set}
if it satisfies the conditions:
\begin{itemize}
\item[(B1)]  If $I,J\in \B$ and $I\cap J\ne \emptyset$, then $I\cup J\in \B$.
\item[(B2)]  $\B$ contains all singletons $\{i\}$, for $i\in S$.
\end{itemize}
\end{definition}

For a building set $\B$ on $S$ and a subset $I\subseteq S$,
define the {\it restriction\/} of $\B$ to $I$ as
$\B|_I:=\{J\in \B\mid J\subseteq I\}$.
Let $\B_{\max}\subset \B$ denote the inclusion-maximal subsets of a building
$\B$.  Then elements of $\B_{\max}$ are pairwise disjoint subsets that
partition the set  $S$.  Call the restrictions $\B|_I$, for $I\in \B_{\max}$,
the {\it connected components\/} of $\B$.  Say that a building set is {\it connected\/} if
$\B_{\max}$ has only one element: $\B_{\max} =\{S\}$.

\begin{example}
\label{ex:graphical_building}
Let $G$ be a graph (with no loops nor multiple edges) on the node set $S$.
The {\it graphical building\/} $\B(G)$ is the set of nonempty subsets $J\subseteq S$
such that the induced graph $G|_J$ on node set $J$ is connected.  Then $\B(G)$ is indeed a building set.

The graphical building set $\B(G)$ is connected if and only if the graph $G$ is connected.
The connected components of the graphical $\B(G)$ building set correspond to connected
components of the graph $G$.  Also each restriction $\B(G)|_I$ is
the graphical building set $\B(G|_I)$ for the induced subgraph $G|_I$.
\end{example}


\begin{definition}
\label{def:nestohedron}
Let $\B$ be a building set on $[n]:=\{1,\dots,n\}$.
Faces of the standard coordinate simplex in $\R^n$
are the simplices $\Delta_I:=\mathrm{ConvexHull}(e_i\mid i\in I)$, for $I\subseteq [n]$,
where the $e_i$ are the endpoints of the coordinate vectors in $\R^n$.

Define the {\it nestohedron\/}\footnote{Called the {\it nested set polytope}
in \cite{Zel}.} $P_\B$
as the Minkowski sum of these simplices
\begin{equation}
\label{eq:nestohedron-definition}
P_\B :=\sum_{I\in \B} y_I \Delta_I,
\end{equation}
where $y_I$ are strictly positive real parameters;
see~\cite[Section~6]{Post}.
\end{definition}

Note that since each of the normal fans $\N(\Delta_I)$ is refined by the
braid arrangement fan,
the same holds for their Minkowski sum \cite[Prop. 7.12]{Ziegler}, and hence
the nestohedra $P_\B$ are generalized permutohedra by Proposition~\ref{prop:permutohedron-summand-characterization}.

  It turns out that $P_\B$ is always a simple polytope, whose
combinatorial structure (poset of faces) does not depend upon the choice of the
positive parameters $y_I$.   In describing this combinatorial structure,
it is convenient to instead describe the dual simplicial complex of $P_\B$.


\begin{definition}
\label{def:nested}
{\rm \cite[Definition~7.3]{Post}} \
For a building set $\B$,
let us say that a subset $N\subseteq \B\setminus \B_{\max}$ is a
{\it nested set\/} if it satisfies the conditions:
\begin{itemize}
\item[(N1)]  For any $I,J\in N$ one has either $I\subseteq J$,
$J\subseteq I$, or $I$ and $J$ are disjoint.
\item[(N2)]  For any collection of $k\geq 2$ disjoint subsets $J_1,\dots,J_k\in N$,
their union $J_1\cup\cdots\cup J_k$ is not in $\B$.
\end{itemize}

Define the {\it nested set complex\/} $\Delta_\B$ as the collection of
all nested sets for $\B$.
\end{definition}

It is immediate from the definition that the nested set complex $\Delta_\B$ is an
abstract simplicial complex on node set $\B$.
Note that this slightly modifies the definition of
a nested set from~\cite{Post}, following~\cite{Zel}, in that one does
not include elements of $\B_{\max}$ in nested sets.

\begin{theorem}
\label{th:nestohedra_dual_to_nested_complexes}
{\rm \cite[Theorem~7.4]{Post}, \cite[Theorem~3.14]{FeichtnerSturmfels}} \
Let $\B$ be a building set on $[n]$.
The nestohedron $P_\B$ is a simple polytope of dimension $n-|\B_{\max}|$.
Its dual simplicial complex is isomorphic to the nested set complex
$\Delta_\B$.
\end{theorem}

An explicit correspondence between faces of $P_\B$ and nested sets
in $\Delta_\B$ is described in~\cite[Proposition~7.5]{Post}.  The dimension of the
face of $P_\B$ associated with a nested set $N\in\Delta_\B$ equals
$n-|N|-|\B_{\max}|$.  Thus vertices of $P_\B$ correspond to inclusion-maximal
nested sets in $\Delta_\B$, and all maximal nested sets contain exactly
$n-|\B_{\max}|$ elements.

\begin{remark}
\label{rem:bistellar-construction}
For a building set $\B$ on $[n]$, it is known \cite[Theorem 4]{FeichtnerYuzvinsky}
that one can obtain the nested set complex
$\Delta_\B$ (resp., the nestohedron $P_\B$)
via the following stellar subdivision
(resp., shaving) construction, a common generalization of
\begin{enumerate}
\item[$\bullet$]
the barycentric subdivision of a simplex as the dual of the permutohedron,
\item[$\bullet$]
Lee's construction of the associahedron \cite[\S3]{Lee}.
\end{enumerate}
Start with an $(n-1)$-simplex whose vertices (resp., facets) have been labelled
by the singletons ${i}$ for $i \in [n]$, which are all in $\B$.  Then proceed
through each of the non-singleton sets $I$ in $\B$, in any order that reverses
inclusion (i.e., where larger sets come before smaller sets), performing a {\it
stellar subdivision} on the face with vertices (resp., {\it shave off\/} the
face which is the intersection of facets) indexed by the singletons in $I$.
\end{remark}

\begin{remark} \rm \
Note that if $\B_1,\dots,\B_k$ are the connected components of a building set $\B$, 
then $P_\B$ is
isomorphic to the direct product of polytopes $P_{\B_1}\times\cdots\times
P_{\B_k}$.  Thus it is enough to investigate generalized permutohedra $P_\B$ and
nested set complexes $\Delta_\B$ only for connected buildings.
\end{remark}

\begin{remark}
\label{rem:building-set-monotonicity}
The definition \eqref{eq:nestohedron-definition} of the nestohedron
$P_\B$ as a Minkowski sum should make it clear that whenever one
has two building sets $\B \subseteq \B'$, then $P_\B$ is a Minkowski
summand of $P_{\B'}$.  Hence Theorem~\ref{th:subdivision-monotonicity} implies
the $h$-vector of $P_{\B'}$ is componentwise weakly larger than that of
$P_\B$.
\end{remark}

\begin{remark}
\label{rem:graphical-monotonicity}
Nestohedra $P_{\B(G)}$ associated with graphical building sets
$\B(G)$
are called {\it graph-associahedra\/}, and have been studied
in~\cite{CarrDevadoss, Post, Tol, Zel}.  In \cite{CarrDevadoss},
the sets in $\B(G)$ are called {\it tubes}, and the nested
sets are called {\it tubings} of~$G$.

In particular, the $h$-vector monotonicity discussed 
in Remark~\ref{rem:building-set-monotonicity}
applies to graph-associahedra $P_{\B(G)}, P_{\B(G')}$ associated to graphs $G, G'$
where $G$ is an edge-subgraph of $G'$.
\end{remark}

\begin{example}  
\label{ex:complete_upper}
({\it Upper bound for nestohedra: the permutohedron})  
see~\cite[Sect.~8.1]{Post} \
For the complete graph $K_n$, the building set $\B(K_n) = 
2^{[n]}\setminus \{\emptyset\}$ consists 
of all nonempty subsets in $[n]$.  
Let us call it the {\it complete building set.}
The corresponding nestohedron (the graph-associahedron
of the complete graph) is the usual $(n-1)$-dimensional permutohedron in $\R^n$.
The $k$-th component $h_k$ of its $h$-vector is the Eulerian number,
that is the number of permutations in $\Sym_n$ with $k$ descents;
and its $h$-polynomial is the Eulerian polynomial $A_n(t)$;
see~\eqref{eq:usual-Eulerian-polynomial}.

This $h$-vector gives the componentwise upper bound
on $h$-vectors for all $(d-1)$-dimensional nestohedra.  
This also implies that the $f$-vector of the permutohedron
gives componentwise upper bound on $f$-vectors of nestohedra.
\end{example}

\begin{example}
\label{ex:simplex-example}
({\it Lower bound for nestohedra: the simplex})  \
The smallest possible connected building set
$\B=\{\{1\},\{2\},\ldots,\{n\}, [n] \}$ gives rise to the nestohedron $P_\B$ 
which is the $(n-1)$-simplex in $\R^n$.
In this case
$$
f(t)=\sum_{i=1}^n \binom{n}{i} t^{i-1}=\frac{(1+t)^n-1}{t}
\quad\text{and}\quad h(t)=1+t+t^2+\cdots+t^{n-1}
$$
give trivial componentwise lower bounds on the $f$-, $h$-vectors
of nestohedra.
\end{example}

\subsection{Two recurrences for $f$-polynomials of nestohedra}

It turns out that there are two useful recurrences for $f$-polynomials of nestohedra
and nested set complexes.

Let $f_\B(t)$ be the $f$-polynomial of the nestohedron $P_\B$:
$$
f_\B(t) := \sum f_i \, t^i = \sum_{N\in\Delta_\B} t^{|S|-|\B_{\max}|-|N|},
$$
where $f_i=f_i(P_\B)$ is the number of $i$-dimensional faces of $P_\B$.
As usual, it is related to the $h$-polynomial as $f_\B(t) = h_\B(t+1)$.

\begin{theorem}
\label{th:f_B_recurrence}
{\rm \cite[Theorem~7.11]{Post}} \
The $f$-polynomial $f_\B(t)$ is determined by the following recurrence relations:
\begin{enumerate}
\item If $\B$ consists of a single singleton, then $f_\B(t)=1$.
\item If $\B$ has connected components $\B_1,\dots, \B_k$, then
$$
f_\B(t) = f_{\B_1}(t)\cdots f_{\B_k}(t).
$$
\item  If $\B$ is a connected building set on $S$, then
$$
f_\B(t) = \sum_{I\subsetneq S} t^{|S|-|I|-1} f_{\B|_{I}}(t).
$$
\end{enumerate}
\end{theorem}

Another recurrence relation for $f$-polynomials was derived in~\cite{Zel},
and will be used in Section~\ref{sec:PDE} below.
It will be more convenient to work with the $f$-polynomial of nested set complexes
$$
\tilde f_\B(t) := \sum_{N\in\Delta_\B} t^{|N|}=
t^{|S|-|\B_{\max}|} f_\B(t^{-1}),
$$
where $\B$ is a building set on $S$.

For a building set $\B$ on $S$ and a subset $I\subset S$, recall that the
restriction of $\B$ to $I$ is defined as $\B|_I = \{J\in \B\mid J\subseteq I\}$.
Also define the {\it contraction\/} of $I$ from $\B$ as the
building set on $S\setminus I$ given by
$$
\B/I := \{J\in S\setminus I \mid J\in \B\textrm{ or } J\cup I\in \B\},
$$
see~\cite[Definition~3.1]{Zel}.
A link decomposition of nested set complexes was constructed in~\cite{Zel}.
It implies the following recurrence relation for the $f$-vector.

\begin{theorem}
\label{th:f_recurrence_zelevinsky}
{\rm \cite[Proposition~4.7]{Zel}} \
For a building set $\B$ on a nonempty set $S$, one has
$$
\frac{d}{d t}\, \tilde f_\B(t) = \sum_{I\in \B\setminus \B_{\max}}
\tilde f_{\B|_I}(t) \cdot \tilde f_{\B/I}(t)
\quad\text{and}\quad
\tilde f_\B(0) = 1.
$$
\end{theorem}

Let $G$ be a simple graph on $S$ and let $I\in \B(G)$, i.e., $I$ is a connected
subset of nodes of $G$.  It has already been mentioned that
$\B(G)|_I = \B(G|_I)$; see Example~\ref{ex:graphical_building}.
Let $G/I$ be the graph on the node set $S\setminus I$
such that two nodes $i, j\in S\setminus I$ are
connected by an edge in $G/I$ if and only if
\begin{enumerate}
\item
$i$ and $j$ are connected by an edge in $G$, or
\item
there are two edges $(i,k)$ and $(j,l)$ in $G$
with $k,l\in I$.
\end{enumerate}
 Then the contraction of $I$ from the graphical building set $\B(G)$ is
the graphical building set associated with the graph $G/I$, that is
$\B(G)/I = \B(G/I)$.

\section{Flag nestohedra}
\label{sec:FlagNestohedra}

This section characterizes the flag nested set complexes and nestohedra, and
then identifies those which are ``smallest''.

\subsection{When is the nested set complex flag?}
\label{ssec:when-is-flag}

For a graphical building set $\B(G)$ it has been observed
(\cite[\S8,4]{Post}, \cite[Corollary 7.4]{Zel}) that
one can replace condition (N2) in  Definition~\ref{def:nested}
with a weaker condition:

\begin{itemize}
\item[(N2')]  For a disjoint pair of subsets $I,J\in N$, one has $I\cup J\not\in \B$.
\end{itemize}

This implies that nested set complexes associated to graphical
buildings are flag complexes.  More generally, one has the following
characterization of the nested set complexes which are flag.

\begin{proposition}
\label{prop:flag-characterization}
For a building set $\B$, the following are equivalent.
\begin{enumerate}
\item[(i)] The nested set complex $\Delta_\B$ (or equivalently, the nestohedron $P_\B$) is
flag.
\item[(ii)] The nested sets for $\B$ are the subsets $N \subseteq \B\setminus \B_{\max}$ 
which satisfy conditions {\rm (N1)} and {\rm (N2').}
\item[(iii)] If $J_1,\ldots,J_\ell \in \B$ with $\ell \geq 2$ are pairwise disjoint 
and their union
$J_1 \cup \cdots \cup J_\ell$ is in $\B$, then one can reindex so that for some $k$ with
$1 \leq k \leq \ell-1$ one has both $J_1 \cup \cdots \cup J_k$ 
and $J_{k+1} \cup \cdots \cup J_\ell$ in $\B$.
\end{enumerate}
\end{proposition}

\begin{proof}
The equivalence of (i) and (ii) essentially follows from the definitions.  We will show
here the equivalence of (i) and (iii).

Assume that (iii) fails, and
let $J_1,\ldots,J_\ell$ provide such a failure with $\ell$  minimal.  Note that this means
$\ell \geq 3$, and minimality of $\ell$
forces $J_r \cup J_s \not\in \B$ for each $r\neq s$; otherwise one could replace the two sets $J_r, J_s$
on the list with the one set $J_r \cup J_s$ to obtain a counterexample of size $\ell-1$.  This means that all of the
pairs $\{J_r,J_s\}$ index edges of $\Delta_\B$, although $\{J_1,\ldots,J_\ell\}$ does not.
Hence $\Delta_\B$ is not flag, i.e., (i) fails.

  Now assume (i) fails, i.e., $\Delta_\B$ is not flag.  Let $J_1,\ldots,J_\ell$ be subsets in
$\B$, for which each pair $\{J_r,J_s\}$ with $r \neq s$ is a nested set,
but the whole collection $M:=\{J_1,\ldots,J_\ell\}$ is not, and assume that this violation has
$\ell$ minimal. Because $\{J_r,J_s\}$ are nested for $r \neq s$, it must be that
$M$ {\it does} satisfy condition (N1), and so $M$ must fail condition (N2).  By minimality of $\ell$,
it must be that the $J_1,\ldots,J_\ell$ are pairwise disjoint and
their union $J_1 \cup \cdots \cup J_\ell$ is in $\B$.
Bearing in mind that $J_r \cup J_s \not\in \B$ for $r \neq s$, it must be that $\ell \geq 3$.  But then
$M$ must give a violation of property (iii), else one could use property (iii) to produce
a violation of (i) either of size $k$ or of size $\ell - k$, which are both smaller than $\ell$.
\end{proof}

\begin{corollary}
For graphical buildings $\B(G)$, the graph-associahedron $P_{\B(G)}$ and nested set
complex $\Delta_{\B(G)}$ are flag.
\end{corollary}

\subsection{Stanley-Pitman polytopes and their relatives}
\label{sec:Stanley-Pitman-example}



One can now use Proposition~\ref{prop:flag-characterization} to characterize the inclusion-minimal
connected building sets $\B$ for which $\Delta_\B$ and $P_\B$ are flag.

\begin{figure}[ht]
\centering
\includegraphics[height=2.3in]{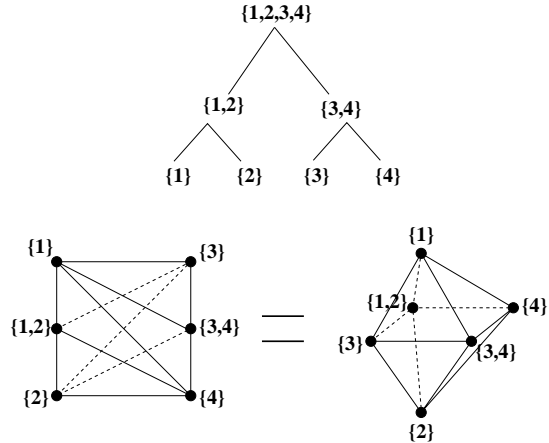}
\caption{A binary tree $\tau$ and building set $\B_\tau$, along with its
complex of nested sets $\Delta_{\B_\tau}$, drawn first as in the construction
of Remark~\ref{rem:bistellar-construction}, and then redrawn as the boundary of an octahedron.}
\label{Stanley-Pitman-relative}
\end{figure}

For any building set $\B$ on $[n]$ with  $\Delta_\B$  flag, one can apply
Proposition~\ref{prop:flag-characterization}(iii) with $\{J_1,\ldots,J_\ell\}$
equal to the collection of singletons $\{\{1\},\ldots,\{n\}\},$ since they are disjoint
and their union $[n]$ is also in $\B$.  Thus after reindexing, some initial segment
$[k]$ and some final segment $[n] \setminus [k]$ must also be in $\B$.  Iterating this,
one can assume after reindexing that there is a {\it plane binary tree} $\tau$ with these
properties
\begin{enumerate}
\item[$\bullet$] the singletons $\{\{1\},\ldots,\{n\}\}$ label the leaves of $\tau$,
\item[$\bullet$] each internal node of $\tau$ is labelled by the set $I$ which is the union of the
singletons labelling the leaves of the subtree below it (so $[n]$ labels
the root node), and
\item[$\bullet$] the building set $\B$ contains of all of the sets labelling nodes in this tree.
\end{enumerate}
It is not hard to see that these sets labelling the nodes of $\tau$
already comprise a building set $\B_\tau$ which satisfies
Proposition~\ref{prop:flag-characterization}(iii), and therefore give rise to a
nested set complex $\Delta_{\B_\tau}$ and nestohedron $P_{\B_\tau}$ which are flag.
See Figure~\ref{Stanley-Pitman-relative} for an example.

The previous discussion shows the following.
\begin{proposition}
The building sets $\B_\tau$ parametrized by plane binary trees $\tau$ are
exactly the inclusion-minimal building sets among those which are
connected and have the nested set complex and nestohedron flag.
\end{proposition}

As a special case, when $\tau$ is the plane binary tree having leaves
labelled by the singletons and internal nodes labelled by all initial segments
$[k]$, one obtains the building set $\B_\tau$ whose nestohedron 
$P_{\B_\tau}$ is {\it Stanley-Pitman polytope\/} from \cite{PS}; see \cite[\S8.5]{Post}.
The Stanley-Pitman polytope is shown there to
be combinatorially (but not affinely) isomorphic to an $(n-1)$-cube;  
the argument given there generalizes to prove the following.

\begin{proposition}
For any plane binary tree $\tau$ with $n$ leaves, the nested set complex $\Delta_{\B_\tau}$
is isomorphic to the boundary of a $(n-1)$-dimensional cross-polytope (hyperoctahedron),
and $P_{\B_\tau}$ is combinatorially isomorphic to an $(n-1)$-cube.
\end{proposition}

\begin{proof}
Note that the sets labelling the non-root nodes of $\tau$ can be grouped
into $n-1$ pairs $\{I_1,J_1\}, \ldots, \{I_{n-1},J_{n-1}\}$ of siblings, meaning that
$I_k,J_k$ are nodes with a common parent in $\tau$.  One then checks that
the nested sets for $\B_\tau$ are exactly the collections $N$ containing
{\it at most one set from each pair} $\{I_k,J_k\}$.  As a simplicial complex,
this is the boundary complex of an $(n-1)$-dimensional cross-polytope
in which each pair $\{I_k,J_k\}$ indexes an antipodal pair of vertices.
\end{proof}

Note that in this case,
$$
f_{\B_\tau}(t)=(2+t)^{n-1},
\quad
h_{\B_\tau}(t)=(1+t)^{n-1},
\quad
\gamma_{\B_\tau}(t)=1 = 1 + 0\cdot t + 0 \cdot t^2 + \cdots.
$$
which gives a lower bound for the $f$- and $h$-vectors of flag nestohedra by
Remark~\ref{rem:building-set-monotonicity}.  If one assumes
Conjecture~\ref{conj:Gal-conjecture}, then it would
also give a lower bound for $\gamma$-vectors of flag nestohedra
(and for flag simplicial polytopes in general).

Note that the permutohedron is a graph-associahedron (and hence a flag
nestohedra).  Therefore, Corollary~\ref{cor:upperbound} implies that the
permutohedron provides the upper bound on the $f$- and $h$-vectors among the
flag nestohedra.



\section{$\B$-trees and $\B$-permutations}
\label{sec:B-trees-permutations}

This section discusses $\B$-trees
and $\B$-permutations, which are two types of combinatorial objects 
associated with vertices of the nestohedron 
$P_\B$.  The $h$-polynomial of $P_\B$ equals the descent-generating
function for $\B$-trees.

\subsection{$\B$-trees and $h$-polynomials}
\label{ssect:B-trees}

This section gives a combinatorial interpretation of
the $h$-polynomials  of nestohedra.
Since nestohedra $P_\B$ are always simple, one should expect some
description of their vertex tree-posets $Q_v$
(see Corollaries~\ref{cor:fans_posets} and~\ref{cor:Map})
in terms of the building set $\B$.

\label{sssec:Bh}
Recall that a {\it rooted tree\/} is a tree with a distinguished node, called its {\it root}.
One can view a rooted tree $T$ as a partial order on its nodes
in which $i <_T j$ if $j$ lies on the unique path from $i$ to the root.
One can also view it as a directed graph in which all edges are directed towards the root;
both viewpoints will be employed here.

For a node $i$ in a rooted tree $T$, let $T_{\leq i}$ denote the set of all descendants
of $i$, that is $j\in T_{\leq i}$ if there is a directed path from the node $j$ to the node
$i$.  Note that $i \in T_{\leq i}$.
Nodes $i$ and $j$ in a rooted tree are called {\it incomparable\/}
if neither $i$ is a descendant of $j$, nor $j$ is a descendant of $i$.


\begin{definition}
\label{def:B-trees}
{\rm \cite[Definition~7.7]{Post}, cf.~\cite{FeichtnerSturmfels}} \
For a connected building set $\B$ on $[n]$, let us define a {\it $\B$-tree\/}
as a rooted tree $T$ on the node set $[n]$ such that
\begin{itemize}
\item[(T1)] For any $i\in[n]$, one has $T_{\leq i}\in \B$.
\item[(T2)] For $k\geq 2$ incomparable nodes $i_1,\dots,i_k\in[n]$,
one has $\bigcup_{j=1}^k T_{\leq i_j} \not\in \B$.
\end{itemize}
\end{definition}

Note that, when the nested set complex $\Delta_\B$ is flag,
that is when $\B$ satisfies any of the conditions of Proposition~\ref{prop:flag-characterization},
one can define a $\B$-tree by requiring condition (T2) only for $k=2$.

\begin{proposition}
\label{prop:B-trees}
{\rm \cite[Proposition~7.8]{Post}, \cite[Proposition~3.17]{FeichtnerSturmfels}} \
For a connected building set $\B$, the map sending a rooted
tree $T$ to the collection of sets 
$\{T_{\leq i} \mid $i$ \text{ is a nonroot vertex}\}\subset \B$
gives a bijection between $\B$-trees and maximal nested sets. 
(Recall that maximal nested sets correspond to the facets of 
the nested set complex $\Delta_\B$ and to the vertices of the 
nestohedron $P_\B$.)

Furthermore, if the $\B$-tree $T$ corresponds to the vertex $v$ of $P_\B$
then $T=Q_v$, that is, $T$ is the vertex tree-poset for $v$
in the notation of Corollary~\ref{cor:fans_posets}.
\end{proposition}


\begin{question}
Does a simple (indecomposable) generalized permutohedron $P$ come from a (connected) building set if
and only if every poset $Q_v$ is a rooted tree, i.e. has a unique maximal element?
\end{question}

Proposition~\ref{prop:B-trees} and Theorem~\ref{th:shelling-and-descents} yield the following corollary.

\begin{corollary}
\label{cor:h_des_T}
For a connected building set $\B$ on $[n]$,
the $h$-polynomial of the generalized permutohedron $P_\B$ is given by
$$
h_\B(t) = \sum_T  t^{\des(T)},
$$
where the sum is over $\B$-trees $T$.
\end{corollary}

The following recursive description of $\B$-trees is 
straightforward from the definition.

\begin{proposition}
\label{prop:B_trees_construction}
{\rm \cite[Section 7]{Post}} \
Let $\B$ be a connected building set on $S$ and let $i\in S$.
Let $\B_1,\dots,\B_r$ be the connected components of the restriction
$\B|_{S\setminus\{i\}}$.
Then all $\B$-trees with the root at $i$ are obtained by picking
a $\B_j$-tree $T_j$, for each component $\B_j$, $j=1,\dots,r$, and
connecting the roots of $T_1,\dots,T_r$ with the node $i$ by edges.
\end{proposition}

In other words, each $\B$-tree is obtained by picking a root $i\in S$,
splitting the restriction $\B|_{S\setminus \{i\}}$ into connected
components, then picking nodes in all connected components, splitting
corresponding restrictions into components, etc.

Recall Definition~\ref{def:generalized_permutohedron} of the surjection $\Psi_\B:= \Psi_{P_\B}$
$$
\Psi_{\B}: \Sym_n \longrightarrow \Vertices(P_\B) = \{\text{$\B$-trees}\},
$$
Here and below one identifies vertices of the nestohedron $P_\B$ with $\B$-trees
via Proposition~\ref{prop:B-trees}.
By Corollary~\ref{cor:Map}, for a $\B$-tree $T$, one has $\Psi_\B(w) = T$ 
if and only if $w$ is a linear extension of $T$.

Proposition~\ref{prop:B_trees_construction} leads to an explicit recursive 
description of the surjection $\Psi_\B$.

\begin{proposition}
\label{prop:B-map} 
Let $\B$ be a connected building set on $[n]$.
Given a permutation $w=(w(1), \dots , w(n)) \in \Sym_n$, one
recursively constructs a $\B$-tree $T=T(w)$, as follows.

The root of $T$ is the node $w(n)$.
Let $\B_1, \dots , \B_r$ be the connected components of
the restriction $\B|_{\{w(1), \dots , w(n-1)\}}$. Restricting $w$ to
each of the sets $\B_i$ gives a subword of $w$, to which one can
recursively apply the construction and obtain a $\B_i$-tree $T_i$.
Then attach these $T_1,\ldots,T_r$ as subtrees of the root node
$w(n)$ in $T$.  This association $w \mapsto T(w)$ is the map $\Psi_{\B}$.
\end{proposition}

\subsection{$\B$-permutations}

It is natural to ask for a nice section of the surjection $\Psi_\B$;
these are the $\B$-permutations defined next.

\begin{definition}
\label{def:characterization}
Let $\B$ be a building set on $[n]$.
Define the set $\Sym_n(\B)\subset \Sym_n$ of {\it $\B$-permutations\/}
as the set of permutations $w\in \Sym_n$ 
such that for any $i\in[n]$, the elements $w(i)$  
and $\max \{w(1), w(2), \dots , w(i) \}$ lie in the same connected component 
of the restricted building set $\B|_{\{w(1), \dots , w(i) \}}$.
\end{definition}

The following recursive construction of $\B$-permutations is 
immediate from the definition.

\begin{lemma}
\label{lem:B-permutation-backward}
A permutation $w\in \Sym_n$ is a $\B$-permutation if and only if
it can be constructed via the following procedure.

Pick $w(n)$ from the connected component of $\B$ that contains $n$;
then pick $w(n-1)$ from the connected component of $\B|_{[n]\setminus\{w(n)\}}$
that contains the {\it maximal} element of $[n] \setminus \{w(n)\}$; then 
pick $w(n-2)$ from the connected component of 
$\B|_{[n] \setminus \{w(n), w(n-1)\}}$ that contains the maximal element of 
$[n]\setminus \{w(n),w(n-1)\}$, etc.  Continue in this manner until $w(1)$ has been chosen.  
\end{lemma}

Let $T$ be a rooted tree on $[n]$ viewed as a tree-poset where the root is the unique maximal element.
The {\it lexicographically minimal linear extension\/} of $T$ is the permutation 
$w\in \Sym_n$ such that $w(1)$ is the minimal leaf of $T$
(in the usual order on $\Z$), $w(2)$ is the minimal leaf of $T-\{w(1)\}$ 
(the tree $T$ with the vertex $w(1)$ removed), 
$w(3)$ is the minimal leaf of $T - \{w(1),w(2)\}$, etc.
There is the following alternative ``backward'' construction for the lexicographically 
minimal linear extension of $T$.

\begin{lemma}
\label{lem:lex_min_backward}
Let $w$ be the lexicographically minimal linear extension of a rooted tree $T$ on $[n]$.
Then the permutation $w$ can be constructed from $T$, as follows: $w(n)$ is the root 
of $T$; $w(n-1)$ is the root of the connected component of $T - \{w(n)\}$ 
that contains the maximal vertex of this forest (in the usual order on $\Z$); 
$w(n-2)$ is the root of the connected component of $T - \{w(n),w(n-1)\}$ 
that contains the maximal vertex 
of this forest, etc.  

In general, $w(i)$ is the root of the connected component 
of the forest 
$$
T - \{w(n),\dots,w(i+1)\}
$$
that contains the vertex $\max(w(1),\dots,w(i))$.
\end{lemma}

\begin{proof}  The proof is by induction on the number of vertices in $T$.  Let $T'$ be the rooted tree obtained from $T$
by removing the minimal leaf $l$.  Then the lexicographically
minimal linear extension $w$ of $T$ is $w = (l,w')$, where
$w'$ is the lexicographically minimal linear extension of $T'$,
and both $w$ and $w'$ are written in list notation.
By induction, $w'$ can be constructed from $T'$ 
backwards.   When one performs the backward construction for $T$, the vertex $l$ can never be the 
root of the connected component of $T -  \{w(n),\dots,w(i+1)\}$ containing the maximal vertex, for $i>1$.
So the backward procedure for $T$ produces the same permutation $w = (l, w')$.
\end{proof}

The next claim gives a correspondence between $\B$-trees and $\B$-permutations.

\begin{proposition}
\label{prop:B-trees-B-permutations}
Let $\B$ be a connected building set on $[n]$.
The set $\Sym_n(\B)$ of $\B$-permutations is exactly the set of lexicographically minimal linear
extensions of the $\B$-trees. (Equivalently, $\Sym_n(\B)$ is the set of lexicographically minimal representatives
of fibers of the map $\Psi_\B$.)  

In particular, the map $\Psi_\B$ induces 
a bijection between $\B$-permutations and $\B$-trees, and $\Sym_n(\B)$ is a section of the map $\Psi_\B$.
\end{proposition}

\begin{proof}
Let $w\in \Sym_n$ be a permutation and let $T=T(w)$ be the corresponding $\B$-tree
constructed as in Proposition~\ref{prop:B-map}.
Note that, for $i=n-1,n-2,\dots,1$, the connected components of the forest $T|_{\{w(1),\dots,w(i)\}}
=T-\{w(n),\dots,w(i+1)\}$ correspond
to the connected components of the building set $\B|_{\{w(1),\dots,w(i)\}}$,
and corresponding components have the same vertex sets.
According to Lemma~\ref{lem:lex_min_backward}, the permutation $w$ is the lexicographically 
minimal linear extension of $T$ if and only if $w$ is a $\B$-permutation as described in 
Lemma~\ref{lem:B-permutation-backward}.
\end{proof}


\section{Chordal building sets and their nestohedra}
\label{sec:chordal}

This section describes an important class of building 
sets $\B$, for which the descent numbers of $\B$-trees are
equal to the descent numbers of $\B$-permutations.
In this case, the $h$-polynomial of the nestohedron $P_\B$
equals the descent-generating function of the corresponding $\B$-permutations.

\subsection{Descents in posets vs.\ descents in permutations}

Let us say
 that a {\it descent} of a permutation
$w\in \Sym_n$ is a pair\footnote{A more standard convention is say that 
a descent is an {\it index\/} $i$ such that $w(i)>w(i+1)$.}
 $(w(i), w(i+1))$ such that $w(i) > w(i+1)$.
Let $\Des(w)$ be the set of all descents in $w$.
Also recall that the descent set $\Des(Q)$ of a poset $Q$ 
is the set of pairs $(a,b)$ such that $a\lessdot_Q b$ 
and $a>_{\mathbb{Z}} b$; 
see Definition~\ref{def:des_tree_poset}.

\begin{lemma}
\label{lem:aux_Descent}
Let $Q$ be any poset on $[n]$, and let $w=w(Q)$ be the
lexicographically minimal linear extension of $Q$.
Then one has $\Des(w)\subseteq \Des(Q)$.
\end{lemma}

\begin{proof}
One must show that any descent $(a,b)$ (with $a>_{\Z}b$)
in $w$ must come from a covering relation $a\lessdot_Q b$
in the poset $Q$.
Indeed, if $a$ and $b$ are incomparable in $Q$, then
the permutation obtained from $w$ by transposing
$a$ and $b$ would be a linear extension of $P$ which is lexicographically
smaller than $w$.  On the other hand, if $a$ and $b$ are comparable
but not adjacent elements in $Q$, then they can never be adjacent
elements in a linear extension of $Q$.
\end{proof}

In particular, this lemma implies that, for a $\B$-tree $T$ and the
corresponding $\B$-permutation $w$ (i.e., $w$ is the lexicographically
minimal linear extension of $T$), one has $\Des(w)\subseteq \Des(T)$.
The rest of this section discusses a special class of building sets
for which one always has $\Des(w) = \Des(T)$.

\subsection{Chordal building sets}

\begin{definition}
\label{def:chordal_building}
Let us say that a building set $\B$ on $[n]$ is {\it chordal\/}
if it satisfies the following condition:  for any $I=\{i_1<\dots <i_r\}\in \B$
and $s=1,\dots,r$,
the subset $\{i_s,i_{s+1},\dots,i_r\}$ also belongs to $\B$.
\end{definition}

Recall that a graph is called {\it chordal\/} if it has no induced $k$-cycles for $k\geq 4$.
It is well known~\cite{FulkersonGross} that chordal graphs are exactly the graphs that admit
a {\it perfect elimination ordering,} which is an ordering of vertices such that, for each vertex $v$, 
the neighbors of $v$ that occur later than $v$ in the order form a clique.
Equivalently, a graph $G$ is chordal if its vertices can be labelled by numbers in $[n]$ so that
$G$ has no induced subgraph $G|_{\{i<j<k\}}$ with the edges $(i,j)$, $(i,k)$ but without the edge $(j,k)$.
Let us call such graphs on $[n]$ 
{\it perfectly labelled chordal graphs.}\footnote{We can also
call them {\it $312$-avoiding graphs\/} because they are exactly the graphs that have no induced
$3$-path $a$---$b$---$c$ with the relative order of the vertices $a,b,c$ as in
the permutation 312.  Note that, unlike the pattern avoidance in permutations,
a $312$-avoiding graph is the same thing as a $213$-avoiding graph.}

\begin{example}  
\label{ex:decreasing_trees}
Let us say that a tree on $[n]$ is {\it decreasing\/} 
if the labels decrease in the shortest path from the vertex $n$ (the root) 
to another vertex.  It is easy to see that decreasing trees are exactly the trees which
are perfectly labelled chordal graphs.  
Clearly, any unlabelled tree has such a decreasing labelling of vertices.
\end{example}

The following claim justifies the name ``chordal building set.''

\begin{proposition}
\label{prop:graphical_chordal}
A graphical building set $\B(G)$ is chordal if and only if
$G$ is a perfectly labelled chordal graph.
\end{proposition}

\begin{proof}  Suppose that $G$ contains an induced subgraph 
$G|_{\{i<j<k\}}$ with exactly two edges $(i,j)$, $(i,k)$.  Then $\{i,j,k\}\in\B(G)$
but $\{j,k\}\not\in \B(G)$.  Thus $\B(G)$ is not a chordal building set.

On the other hand, suppose that $\B(G)$ is not chordal.
Then one can find a connected subset $I=\{i_1<\cdots < i_r\}$ of vertices in $G$
such that $\{i_s,i_{s+1},\dots,i_r\}\not\in\B(G)$, for some $s$.  In other words,
the induced graph $G'= G|_{\{i_s,\dots,i_k\}}$ is disconnected.  Let us pick a shortest
path $P$ in $G|_{\{i_1,\dots,i_r\}}$ that connects two different components of $G'$.
Let $i$ be the minimal vertex in $P$ and let $j$ and $k$ be the two vertices adjacent of $i$
in the path $P$.  Clearly, $j>i$ and $k>i$.  It is also clear that $(i,j)$ is not 
an edge of $G$.  Otherwise there is a shorter path obtained from $P$ by replacing
the edges $(i,j)$ and $(i,k)$ with the edge $(j,k)$.  So one has found a forbidden
induced subgraph $G|_{\{i,j,k\}}$.  Thus $G$ is not a perfectly labelled chordal graph.
\end{proof}

\begin{proposition}
\label{prop:des=des-chordal}
Let $\B$ be a connected chordal building set.
Then, for any $\B$-tree $T$ and the corresponding $\B$-permutation $w$,
one has $\Des(w)=\Des(T)$.
\end{proposition}

\begin{proof}
Let $T$ be a $\B$-tree and let $w$ be the corresponding $\B$-permutation,
which can be constructed backward from $T$ as described in Lemma~\ref{lem:lex_min_backward}.
Let us fix $i\in \{n-1,n-2,\dots,1\}$.  
Let $T_1,\dots,T_r,T_1',\dots,T_s'$ be the connected components of the forest
$T - \{w(n),w(n-1),\dots,w(i+1)\}$, where $T_1,\dots,T_r$ are the subtrees whose roots are
the children of the vertex $w(i+1)$, and $T_1',\dots,T_s'$ are the remaining subtrees.
Let $I=T_{\leq w(i+1)} \subset [n]$ be the set of all descendants of $w(i+1)$ in $T$.
By Definition~\ref{def:B-trees}(T1), one has $I\in \B$. 

Suppose that the vertex $m=\max(w(1),\dots,w(i))$ appears in one of the subtrees $T_1,\dots,T_r$,
say, in the tree $T_1$.  Then, by  Lemma~\ref{lem:lex_min_backward}, $w(i)$ should be the root of $T_1$. 
We claim that all vertices in $T_2,\dots,T_r$ are less than $w(i+1)$.
Indeed, this is clear if $w(i+1)$ is the maximal element in $I$.
Otherwise, the set $I'=I\cap\{w(i+1)+1,\dots, n-1,n\}$ is nonempty, $I'\in \B$ because $\B$ is chordal,
and $I'$ contains the maximal vertex $m$.  Since the vertex set $J$ of $T_1$ should be an element of $\B$,
it follows that $I'\subseteq J$.  So all vertices of $T_2,\dots,T_r$ are less than $w(i+1)$.

Thus none of the edges of $T$ joining the vertex $w(i+1)$ with the roots
of $T_2,T_3,\dots,T_r$ can be a descent edge.  The only potential descent
edge is the edge $(w(i),w(i+1))$ that attaches the subtree $T_1$ to $w(i+1)$.
This edge will be a descent edge in $T$ if and only if $w(i)>w(i+1)$, i.e.,
exactly when $(w(i),w(i+1))$ is a descent in the permutation $w$.

Now suppose that the maximal vertex $m=\max(w(1),\dots,w(i))$ appears in one of the remaining subtrees
$T_1',\dots,T_s'$, which are not attached to the vertex $w(i+1)$, say, in $T_1'$.
In this case $w(i+1)$ should be greater than all $w(1),\dots,w(i)$.
(Otherwise, if $w(i+1)<m$, then at the previous step of the backward construction for $w$,
$T_1'$ is the connected component of $T-\{w(n),\dots,w(i+1)\}$ that contains 
the vertex $\max(w(1),\dots,w(i+1)) =m$.  So $w(i+1)$ should have been the root of $T_1'$.)
In this case, none of the edges joining the vertex $w(i+1)$ with the components $T_1,\dots,T_r$ can be
a descent edge and $(w(i),w(i+1))$ cannot be a descent in $w$.

This proves that descent edges of $T$ are in bijection with descents in $w$.
\end{proof}

Corollary~\ref{cor:h_des_T} and Proposition~\ref{prop:des=des-chordal} imply 
the following formula.

\begin{corollary}
\label{cor:h_des_chordal}
For a connected chordal building set $\B$, the $h$-polynomial of the nestohedron 
$P_\B$ equals
$$
h_\B(t) = \sum_{w\in \Sym_n(\B)}  t^{\des(w)},
$$
where $\des(w)$ is the usual descent number of a permutation $w\in \Sym_n(\B)$.
\end{corollary}

Let us give an additional nice property of nestohedra for chordal building sets.

\begin{proposition}
\label{prop:chordal_are_flag}
 For a chordal building set $\B$, the nestohedron $P_\B$ is a flag simple polytope.
\end{proposition}

\begin{proof}
Let us check that a chordal building set $\B$ satisfies
the condition in Proposition~\ref{prop:flag-characterization}(iii).
Using the notation of that proposition, 
let $J_1 \cup \dots \cup J_\ell = \{i_1< \cdots < i_r\}$.
Let $U_s$ be the union of those subsets $J_1,...,J_\ell$ that 
have a nonempty intersection with $\{i_s, i_{s+1},\dots,i_r\}$.
Since $\{i_s, i_{s+1},\dots,i_n\}$ is in $\B$ (because $\B$ is chordal), 
the subset $U_s$ should also be in $\B$ (by 
Definition~\ref{def:building-set-definition}(B1)).  
Clearly, $U_1$ is the union of all $J_i$'s
and $U_r$ consists of a single $J_i$.
It is also clear that $U_{j+1}$ either equals $U_j$ or is obtained
from $U_j$ by removing a single subset $J_i$.
It follows that there exists an index $s$ such that 
$U_s = (J_1 \cup \dots \cup J_\ell)\setminus J_i$.
This gives an index $i$ such that  
$(J_1 \cup \dots \cup J_\ell)\setminus J_i$ and $J_i$ are both in $\B$,
as needed.
\end{proof}

\section{Examples of nestohedra}
\label{sec:three-graph-ass-examples}

Let us give several examples which illustrate 
Corollary~\ref{cor:h_des_T} and Corollary~\ref{cor:h_des_chordal}.
The $f$- and $h$-numbers for the permutohedron and associahedron are well-known.

\subsection{The permutohedron}
\label{ssect:complete-graph-example-1}

For the complete building set $\B=\B(K_n)$ 
the nestohedron $P_\B$ is the usual permutohedron;
see Example~\ref{ex:complete_upper} and~\cite[Sect.~8.1]{Post}.   
In this case $\B$-trees are linear orders on $[n]$ and
$\B$-permutations are all permutations 
$\Sym_n(\B) = \Sym_n$.
Thus, as noted before in Example~\ref{ex:complete_upper},
the $h$-polynomial is the usual Eulerian polynomial $A_n(t)$,
and the $h$-numbers are the 
Eulerian numbers $h_k(P_{\B}) = A(n,k):=\#\{w\in\Sym_n\mid \des(w)=k\}$.

\subsection{The associahedron} 
\label{ssect:path-graph-example-1}

Let $G=\A_n$ denote the graph which is a path having $n$ nodes labelled consecutively
$1,\dots,n$.  The graphical building set $\B=\B(\A_n)$ consists of
all intervals $[i,j]$, for $1\leq i\leq j\leq n$.
The corresponding nestohedron
$P_{\B(\A_n)}$ is the usual Stasheff associahedron; 
see \cite{CarrDevadoss, Post}.

In this case, the $\B$-trees correspond to unlabelled plane binary trees on $n$ nodes,
as follows;  see~\cite[Sect.~8.2]{Post} for more details.
A {\it plane binary tree\/} is a rooted tree with two types of edges
(left and right) such that every node has at most one left and
at most one right edge descending from it.  From Proposition~\ref{prop:B_trees_construction},
one can see that a $\B$-tree is a binary tree with $n$ nodes labelled
$1,2,\ldots,n$ so that, for any node, all nodes in its left (resp., right)
branch have smaller (resp., bigger) labels.
Conversely, given an unlabelled plane binary tree, there is a unique
way to label its nodes $1,2,\dots,n$ to create a
$\B$-tree, namely in the order of traversal of a {\it depth-first search}.
Furthermore, note that descent edges correspond to right edges.

It is well-known that the
number of unlabelled binary trees on $n$ nodes is equal to
the {\it Catalan number\/} $C_n = \frac{1}{n+1}\binom{2n}{n}$,
and the number of binary trees on $n$
nodes with $k-1$ right edges is the {\it Narayana number}
$N(n,k) = \frac{1}{n} \binom{n}{k} \binom{n}{k-1}$;
see~\cite[Exer.~6.19c and Exer. 6.36]{EC2}.
Therefore, the $h$-numbers of the associahedron $P_{\B(\A_n)}$
are the Narayana numbers:
$h_k(P_{\B(\A_n})) = N(n,k+1)$, for $k=0,\dots,n-1$.

%
%

It is also well-known that the $f$-numbers of the associahedron
are $f_k(P_{\B(\A_n})) = \frac{1}{n+1} \binom{n-1}{k} \binom{2n-k}{n}$.
This follows from a classical Kirkman-Cayley formula \cite{Cayley} for
the number of ways to draw $k$ noncrossing diagonals in an $n$-gon.

In this case, the $\B$-permutations are exactly $312$-avoiding
permutations $w\in \Sym_n$.  Recall that a permutation $w$ is
{\it $312$-avoiding} if there is no triple of indices
$i<j<k$ such that $w(j) < w(k) < w(i)$.
Thus Corollary~\ref{cor:h_des_chordal} says that
the $h$-polynomial of the associahedron $P_{\B(\A_n)}$ is
$\sum_{w} t^{\des(w)}$
where the sum runs over all $312$-avoiding permutations in $\Sym_n$.
This is consistent with the known fact that the Narayana numbers
count $312$-avoiding permutations according to their number of descents;
see Simion \cite[Theorem 5.4]{Simion-NCstats} for a stronger
statement.

\subsection{The cyclohedron}\label{ssect:cyclo}

If $G=\mathrm{Cycle}_n$ is the $n$-cycle, then the nestohedron 
$P_{\B(\mathrm{Cycle}_n)}$ is the {\it cyclohedron\/} also introduced by
Stasheff; see~\cite{CarrDevadoss, Post}.
The $h$-polynomial of the cyclohedron was computed by
Simion \cite[Corollary 1]{Simion}:
\begin{equation}
\label{eq:h_cyclohedron}
h_{\B(\cycle_n)}(t)=\sum_{k=0}^{n} \binom{n}{k}^2 t^k.
\end{equation}

Note that the $n$-cycle (for $n>3$) is not a chordal graph, so
Corollary~\ref{cor:h_des_chordal} does not apply to this case.

\subsection{The stellohedron} 
\label{ssect:star-graph-example-1}

Let $m=n-1$.  Let  $G=K_{1,m}$ be the $m$-star graph with the central node
$m+1$ connected to the nodes $1,\dots,m$.  
Let us call the associated polytope $P_{\B(K_{1,m})}$ the
{\it stellohedron.}

From Proposition~\ref{prop:B_trees_construction} 
one sees that $\B(K_{1,m})$-trees are in  bijection 
with {\it partial permutations\/} of $[m]$,
which are ordered sequences $u=(u_1,\dots,u_r)$ of distinct numbers in
$[m]$, where $r=0,\dots,m$.  The tree $T$ associated to a partial permutation 
$u=(u_1,\dots,u_r)$ has the edges
$$
(u_r,u_{r-1}), \dots, (u_2,u_1), (u_1,m+1), (m+1,i_1), \ldots, (m+1,i_{m-r})
$$
where $i_1,\dots,i_{m-r}$ are the elements of 
$[m]\setminus \{u_1,\dots,u_r\}$.
The root of $T$ is $u_r$ if $r\geq 1$, or $m+1$ if $r=0$.
For $r\geq 1$, one has $\des(T) = \des(u) + 1$,
where the descent number of a partial permutation is
$$
\des(u) :=\#\{i=1,\dots,r-1\mid u_{i}>u_{i+1}\}.
$$
Also for the tree $T$ associated with the empty partial permutation
(for $r=0$) one has $\des(T)=0$.
Corollary~\ref{cor:h_des_T} then says that
\begin{equation}
h_{\B(K_{1,m})}(t) =  1 + \sum_{u} t^{\des(u) + 1} =
1 + \sum_{r=1}^m \binom{m}{r} \sum_{k=1}^{r} A(r,k) \,t^{k},
\end{equation}
where the first sum is over nonempty partial permutations $w$ of $[m]$.
In particular, the total number of vertices of the stellohedron
$P_{\B(K_{1,m})}$ equals
$$
f_0(P_{\B(K_{1,m})}) =
\sum _{r=0}^m \binom{m}{r} \cdot r! =
\sum _{r=0}^m \frac{m!}{r!}.
$$
This sequence appears in Sloan's On-Line Encyclopedia of Integer 
Sequences\footnote{{\tt http://akpublic.research.att.com/\~{}njas/sequences/}} 
as A000522.  

In this case, $\B(K_{1,m})$-permutations are
permutations $w\in \Sym_{m+1}$ such that $m+1$
appears before the first descent.
Such permutations $w$ are in bijection with partial permutations $u$ of $[m]$.
Indeed, $u$ is the part of $w$ after the entry $m+1$.
Since our labelling of $K_{1,m}$ (with the central node
labelled $m+1$) is decreasing (see Example~\ref{ex:decreasing_trees}),
Corollary~\ref{cor:h_des_chordal} implies that
the $h$-polynomial of the stellohedron $P_{B(K_{1,m})}$ is
$h(t)=\sum_{w} t^{\des(w)}$,
where the sum runs over all such permutations $w\in \Sym_{m+1}$.
This agrees with the above expression in terms of partial permutations.

\subsection{The Stanley-Pitman polytope}

Let $\B_{\PS}=\{[i,n], \{i\}\mid i=1,\dots,n\}$ 
(the collection of all intervals $[i,n]$ and singletons $\{i\}$).
This (non-graphical) building set is chordal.
According to \cite[\S8.5]{Post}, the corresponding nestohedron 
$P_{\B_{\PS}}$ is the Stanley-Pitman polytope from \cite{PS}.

By Proposition~\ref{prop:B_trees_construction},
$\B_{\PS}$-trees have the following form $T(I)$.  For an increasing sequence 
$I$ of positive integers $i_1<i_2<\cdots <i_k=n$, construct the
tree $T(I)$ on $[n]$ with the root at $i_1$ and the 
chain of edges $(i_1,i_2),(i_2,i_3),\dots, (i_{k-1},i_k)$;
also, for each $j\in[n]\setminus I$,
one has the edge $(i_l,j)$ where $i_l$ is the minimal element of $I$
such that $i_l>j$. 

In this case, $\B_{\PS}$-permutations are permutations $w\in \Sym_n$
such that $w(1)<w(2)<\cdots < w(k)>w(k+1)> \cdots > w(n)$, for some
$k=1,\dots,n$.

Using $\B_{\PS}$-trees or $\B_{\PS}$-permutations one can easily deduce
that the $h$-polynomial of the Stanley-Pitman polytope
is $h_{\B_{\PS}}(t) = (1+t)^{n-1}$.  This is not surprising since 
$P_{\B_{\PS}}$ is combinatorially isomorphic to the $(n-1)$-dimensional cube.


\section{$\gamma$-vectors of nestohedra}
\label{sec:gamma_nestohedra}

Recall that the $\gamma$-vector $(\gamma_0,\gamma_1,\dots,\gamma_{\lfloor d/2 \rfloor})$
of a $d$-dimensional simple polytope is defined via its $h$-polynomial as
$h(t) = \sum \gamma_i \,t^i (1+t)^{d-2i}$; and the $\gamma$-polynomial
is $\gamma(t) = \sum \gamma_i \,t^i$; see Section~\ref{ssec:flag_simple_gamma}.

The main result of this section is a formula for the $\gamma$-polynomial of a
chordal nestohedron as a descent-generating function 
(or peak-generating function) for some set of permutations.
This implies that Gal's conjecture (Conjecture~\ref{conj:Gal-conjecture}) holds
for this class of flag simple polytopes.

\subsection{A warm up: $\gamma$-vector for the permutohedron}
\label{ssec:permutohedron-gamma}

We review here the beautiful construction of Shapiro, Woan, and Getu \cite{ShapiroWoanGetu} 
that leads to a nonnegative formula for the $\gamma$-vector 
of the usual permutohedron.
This subsection also serves as a warm-up for a more general construction
in the following subsection.

Some notation is necessary.  
Recall that a {\it descent} in a permutation $w\in\Sym_n$ is
a pair $(w(i),w(i+1))$ such that $w(i)>w(i+1)$, where $i\in[n-1]$.
A {\it final descent} is when $w(n-1) > w(n)$, and
a {\it double descent} is a pair of consecutive descents,
i.e.\ a triple $w(i) > w(i+1) > w(i+2)$.


Additionally,
define a {\it peak} of $w$ to be an entry 
$w(i)$ for $1\leq i\leq n$ such that $w(i-1)<w(i)>w(i+1)$.
Here (and below) set $w(0)=w(n+1)=0$ and so a peak can occur
in positions $1$ or $n$.  On the other hand, a
{\it valley} of $w$ is an entry $w(i)$ for $1 < i < n$ such that
$w(i-1)>w(i)<w(i+1)$.
The {\it peak-valley sequence\/} of $w$ is the subsequence 
in $w$ formed by all peaks and valleys.

Let $\widehat{\Sym}_n$ denote the set of permutations in $\Sym_n$
which do not contain any final descents or double descents.
Let $\peak(w)$ denote the number of peaks in a permutation $w$.
It is clear that $\peak(w)-1 = \des(w)$,
for permutations $w\in \widehat{\Sym}_n$ (and only for these permutations).


\begin{theorem}
\label{th:CDG}
\rm{(cf.~ \cite[Proposition 4]{ShapiroWoanGetu})} \
The $\gamma$-polynomial of the usual permutohedron
$P_{\B(K_n)}$ is
\begin{equation*}
\sum_{w \in \widehat{\Sym}_n} t^{\peak(w)-1}
=\sum_{w \in \widehat{\Sym}_n} t^{\des(w)}.
\end{equation*}
\end{theorem}

\begin{example}
Let us calculate the $\gamma$-polynomial of the
two dimensional permutohedron $P_{\B(K_3)}$.
One has $\widehat{\Sym}_3 = \{(1,2,3), (2,1,3), (3,1,2)\}$.
Of these, $(1,2,3)$ has one peak (and no descents),
and $(2,1,3)$ and $(3,1,2)$ have two peaks (and one descent).
Therefore, the $\gamma$-polynomial is $1+2t$.
\end{example}

Say that an entry $w(i)$ of $w$ is an {\it intermediary\/} entry
if $w(i)$ is not a peak or a valley.  Say that $w(i)$ is an {\it
ascent-intermediary entry\/} if $w(i-1) < w(i) < w(i+1)$ and that it is
a {\it descent-intermediary entry\/} if $w(i-1)> w(i) > w(i+1)$. 
(Here again one should assume that $w(0)=w(n+1)=0$.) 
Note that the set $\widehat{\Sym}_n$ is exactly the set of permutations 
in $\Sym_n$ without descent-intermediary entries.

\begin{figure}[ht]
\centering
\includegraphics[height=1.6in]{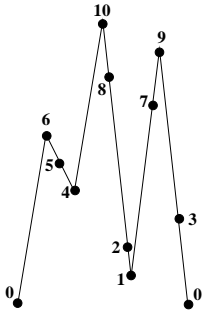}
\caption{Mountain range $M_{w}$ for $w = (6,5,4,10,8,2,1,7,9,3)$}
\label{fig:perm-graph}
\end{figure}

It is convenient to graphically represent a permutation $w\in\Sym_n$
by a piecewise linear ``mountain range'' $M_w$ obtained by connecting 
the points $(x_0,0)$, $(x_1,w(1))$, $(x_2,w(2))$, \dots, $(x_n,w(n))$, $(x_{n+1},0)$
on $\R^2$ by straight line intervals, for some $x_0<x_1<\dots<x_{n+1}$;
see Figure~\ref{fig:perm-graph}.  
Then peaks in $w$ correspond to local maxima
of $M_w$, valleys correspond to local minima of $M_w$, ascent-intermediary
entries correspond to nodes on ascending slopes of $M_w$,
and descent-intermediary entries correspond to nodes on descending slopes of $M_w$.
For example, the permutation $w = (6,5,4,10,8,2,1,7,9,3)$ shown in Figure~\ref{fig:perm-graph}
has three peaks $6,10,9$, two valleys $4,1$,
one ascent-intermediary entry $7$, and four descent-intermediary 
entries $5,8,2,3$.  Its peak-valley sequence is $(6,4,10,1,9)$.

As noted in Section~\ref{ssec:descents_of_tree_posets}, the $h$-polynomial of
the permutohedron is the descent-generating function for permutations in
$\Sym_n$ (the Eulerian polynomial).  In order to prove Theorem~\ref{th:CDG}, one
constructs an appropriate partitioning of $\Sym_n$ into equivalence classes
(cf.\ Remark~\ref{Boolean-classes-remark}), where each equivalence classes has
exactly one element from $\widehat{\Sym}_n$.  To describe the equivalence
classes of permutations, one must introduce some operations on permutations.

\begin{definition}
\label{def:leap}
Let us define the {\it leap\/} operations $L_a$ and $L_a^{-1}$ that act on permutations.  
Informally, the permutation $L_a(w)$ is obtained from $w$ by moving an
intermediary node $a$ on the mountain range $M_w$ directly to the right until
it hits the next slope of $M_w$.  The permutation $L_a^{-1}(w)$ is obtained
from $w$ by moving $a$ directly to the left until it hits the next slope of
$M_w$.

More formally, for an intermediary entry $a=w(i)$ in $w$, 
the permutation $L_a(w)$ is obtained from $w$ by removing $a$ from the $i$-th position
and inserting $a$ in the position between $w(j)$ and $w(j+1)$,
where  $j$ is the minimal index such that
$j>i$ and $a$ is between $w(j)$ and $w(j+1)$, i.e., $w(j)<a<w(j+1)$ or
$w(j)>a>w(j+1)$.  The leap operation $L_a$ is not defined if all entries following $a$
in $w$ are less than $a$.

Similarly, the inverse operation $L_a^{-1}(w)$ is given by removing $a$ from the
$i$-th position in $w$ and inserting $a$ between $w(k)$ and $w(k+1)$, where $k$ is
the 
maximum index such that $k<i$ and $a$ is between $w(k)$ and $w(k+1)$.  The
operation $L_a^{-1}$ is is not defined if all entries preceding $a$ in $w$
are less than $a$.
\end{definition}

For example, for the permutation $w$ 
shown on Figure~\ref{fig:perm-graph}, one has
$L_2(w) = (6,5,4,10,8,1,2,7,9,3)$
and $L_2^{-1}(w) = (2,6,5,4,10,8,1,7,9,3)$.

Clearly, if $a$ is an ascent-intermediary entry in $w$ then $a$ is
a descent-intermediary entry in $L_a^{\pm 1}(w)$, and vise versa.  Note that if
$a$ is an ascent-intermediary entry in $w$, then $L_a(w)$ is always defined, and
if $a$ is a decent-intermediary entry, then $L_{a}^{-1}(w)$ is always defined.

\begin{definition}  Let us also define the {\it hop\/} operations $H_a$ on permutations.
For an ascent-intermediary entry $a$ in $w$, define $H_a(w) = L_a(w)$;
and, for a descent-intermediary entry $a$ in $w$, define $H_a(w) = L_a^{-1}(w)$.
\end{definition}

For example, for the permutation $w$ shown on Figure~\ref{fig:perm-graph},  
the permutation $H_2(w) = (2,6,5,4,10,8,1,7,9,3)$ 
is obtained by moving the descent-intermediary 
entry $2$ to the left to the first ascending slope,
and $H_7(w) =(6,5,4,10,8,2,1,9,7,3)$ is obtained by moving 
the ascent-intermediary entry $7$ to the right to the last descending slope.

Note that leaps and hops never change the shape of the mountain range
$M_w$, that is, they never change the peak-valley sequence of $w$.
They just move intermediary nodes from one slope of $M_w$ to another.
It is quite clear from the definition that all leap and hop operations pairwise commute 
with each other.  It is also clear that two hops $H_a$ get us back to the original
permutation.

\begin{lemma}
\label{lem:hops_permutohedron}
For intermediary entries $a$ and $b$ in $w$,
one has $(H_a)^2(w) = w$ and $H_a (H_b(w)) = H_b (H_a(w))$.
\end{lemma}

Thus the hop operations $H_a$ generate the action of the group 
$(\Z/2\Z)^m$ on the set of permutations with a given peak-valley
sequence, where $m$ is the number of intermediary entries in
such permutations.

Let us say that two permutations are {\it hop-equivalent\/}
if they can be obtained from each other by the hop operations
$H_a$ for various $a$'s.
The partitioning of $\Sym_n$ into hop-equivalence classes
allows us to prove Theorem~\ref{th:CDG}.

\begin{proof}[Proof of Theorem~\ref{th:CDG}]
The number $\des(w)$ of descents in $w$ equals
the number of peaks in $w$ plus the number of descent-intermediary
entries in $w$ minus $1$ (because the last entry is either a peak
or a descent-intermediary entry, but it does not contribute a descent).
Notice that if $a$ is an ascent-intermediary (resp., descent-intermediary)
entry in $w$ then the number of descent-intermediary entries in $H_a(w)$
increases (resp., decreases) by 1 and the number of peaks does not change.

If $w\in\Sym_n$ has $p=\peak(w)$ peaks then it has $p-1$ valleys and $n-2p+1$ 
intermediary entries.  Lemma~\ref{lem:hops_permutohedron} implies that 
the hop-equivalence class $C$ of $w$ involves $2^{n-2p+1}$ permutations.
Moreover, the descent-generating function for these permutations
is $\sum_{u\in C} t^{\des(u)} = t^p (t+1)^{n-2p+1}$.
Each hop-equivalence class has exactly one representative $u$
without descent-intermediary entries, that is $u\in\widehat{\Sym}_n$.
Thus, summing the contributions of hop-equivalence classes, one can 
write the $h$-polynomial of the permutohedron as
$$
h(t)= \sum_{w\in\Sym_n} t^{\des(w)} = 
\sum_{w \in \widehat{\Sym}_n} t^{\peak(w)-1} (t+1)^{n+1-2\,\peak(w)}.
$$
Comparing this to the definition of the $\gamma$-polynomial,
one derives the theorem.
\end{proof}


\subsection{$\gamma$-vectors of chordal nestohedra}
\label{ssec:gamma-chordal}

According to Proposition~\ref{prop:chordal_are_flag}, nestohedra for chordal building sets 
are flag simple polytopes. 
Thus Gal's conjecture (Conjecture~\ref{conj:Gal-conjecture}) applies.
This section proves this conjecture and
present a nonnegative combinatorial formula for 
$\gamma$-polynomials of such  nestohedra as peak-generating functions
for some subsets of permutations.

Let $\B$ be a connected chordal building set on $[n]$.
Recall that $\Sym_n(\B)$ is the set of $\B$-permutations;
see Definition~\ref{def:characterization}.  Let
$\widehat{\Sym}_n(\B):=\Sym_n(\B)\cap \widehat{\Sym}_n$ be the subset
of $\B$-permutations which have no final descent or double descent.

The following theorem is the main result of this section.

\begin{theorem}
\label{th:gamma-chordal}
For a connected chordal building $\B$ on $[n]$,
the $\gamma$-polynomial of the nestohedron $P_\B$
is the peak-generating function for the permutations in $\widehat{\Sym}_n(\B)$:
\begin{equation*}
  \gamma_{\B}(t) = \sum_{w \in \widehat{\Sym}_n(\B)}
  t^{\peak(w)-1}
   = \sum_{w \in \widehat{\Sym}_n(\B)}
  t^{\des(w)}.
\end{equation*}
\end{theorem}

As noted earlier, $\peak(w)-1 = \des(w)$ for $w\in \widehat{\Sym}_n$.

The proof of Theorem~\ref{th:gamma-chordal} will be an extension of the proof
given for the $\gamma$-vector of the permutohedron in
Section~\ref{ssec:permutohedron-gamma}.
Recall that Corollary~\ref{cor:h_des_chordal} interprets
the $h$-polynomial of $P_\B$ as the descent-generating function
for $\B$-permutations $w\in\Sym_n(\B)$.
Theorem~\ref{th:gamma-chordal} will be proven by constructing
an appropriate partitioning of the set $\Sym_n(\B)$ into equivalence classes,
where each equivalence class has exactly one representative
from $\widehat{\Sym}_n(\B)$.
As before, one uses (suitably generalized) hop operations
to describe equivalence classes of elements of $\Sym_n(\B)$.

One needs powers of the leap operations $L_a^{r} := (L_a)^{r}$, for $r\geq 0$, and 
$L_a^{r} := (L_a^{-1})^{-r}$, for $r\leq 0$; see Definition~\ref{def:leap}.
In other words, for $r>0$, $L_a^{r}(w)$ is obtained from $w$ by moving 
the intermediary entry $a$ to the right until it hits the $r$-th slope from its 
original location; and, for $r<0$, by moving $a$ to the left until it hits
$(-r)$-th slope from its original location.
Clearly, $L_a^{r}(w)$ is defined whenever $r$ is in a certain integer 
interval $r\in[r_{\min},r_{\max}]$.
It is also clear that, if $a$ is an ascent-intermediary entry in $w$,
then $a$ is ascent-intermediary in $L_a^r(w)$ for even $r$
and $a$ is descent-intermediary in $L_a^r(w)$ for odd $r$,
and vice versa if $a$ is descent-intermediary in $w$.

Note that for a $\B$-permutation $w\in \Sym_n(\B)$, the permutations
$L_a^r(w)$ may no longer be $\B$-permutations.  The next lemma 
ensures that at least some of them will be $\B$-permutations.

\begin{lemma} 
\label{lem:leaps}
Let $\B$ be a chordal building on $[n]$.
Suppose that $w\in \Sym_n(\B)$ is a $\B$-permutation.

{\rm (1)} 
If $a$ is an ascent-intermediary letter in $w$,
then there exists an odd positive integer $r>0$ such that 
$L_a^r(w)\in \Sym_n(\B)$ and $L_a^s(w)\not \in \Sym_n(\B)$,
for all $0<s<r$.

{\rm (2)}
If $a$ is a descent-intermediary letter in $w$,
then there exists an odd negative integer $r<0$ such that 
$L_a^r(w)\in \Sym_n(\B)$ and $L_a^s(w)\not \in \Sym_n(\B)$,
for all $0>s>r$.
\end{lemma}

The proof of Lemma \ref{lem:leaps} will require some preparatory notation
and observations.

For a permutation $w\in \Sym_n$ and $a\in[n]$ such that $w(i)=a$, 
let 
$$
\{w\nw a\}:= \{w(j)\mid j\leq i,\, w(j)\geq a\}
$$
be the set of all entries in $w$ which are located to the left of $a$
and are greater than or equal to $a$ (including the entry $a$ itself).
The arrow in this notation refers to our graphical representation of a permutation
as a mountain range $M_w$: the set
$\{w \nw a\}$ is the set of entries in $w$ located to the North-West of the entry $a$.

According to Definition~\ref{def:characterization}, 
the set $\Sym_n(\B)$ is the set 
of permutations $w$ such that, for $i=1,\dots,n$, there exists
$I\in\B$ such that both $w(i)$ and $\max(w(1),\dots,w(i))$ are in $I$
and $I\subset \{w(1),\dots,w(i)\}$.
If $\B$ is chordal, then $I' := I\cap [w(i),\infty]$ also belongs to $\B$
(see Definition~\ref{def:chordal_building}) and satisfies the same properties.
Clearly $\max(w(1),\dots,w(i)) = \max \{w \nw w(i)\}$.
Thus, for a chordal building set, one can reformulate 
Definition~\ref{def:characterization} of $\B$-permutations as follows.

\begin{lemma}
\label{lem:B-permutations-chordal}
Let $\B$ be a chordal building set.
Then $\Sym_n(\B)$ is the set of permutations $w\in \Sym_n$ such that
for any $a\in[n]$, the elements $a$ and $\max \{w \nw a\}$ are in the same connected
component of $\B|_{\{w \nw a\}}$.  Equivalently, there exists $I\in \B$ such that 
$a\in I$,  $\max \{w \nw a\}\in I$, and $I \subset \{w \nw a\}$.
\end{lemma}

Let us now return to the setup of Lemma~\ref{lem:leaps}.
There are 2 possible reasons why the permutation  $u=L_a^{r}(w)$ may no longer 
be a $\B$-permutation, that is, fail to satisfy the conditions
in Lemma~\ref{lem:B-permutations-chordal}:
\begin{itemize}
\item[(A)]
It is possible that the entry $a$ 
and the entry $\max \{u \nw a\}$ are in different connected
components of $\B|_{\{u \nw a\}}$.
\item[(B)] It is also possible that another entry $b\ne a$
in $u$ and $\max \{u \nw b\}$ are in different connected
components of $\B|_{\{u \nw b\}}$.
\end{itemize}
Let us call these two types of failure 
{\it A-failure\/} and {\it B-failure.}
The following auxiliary result is needed.

\begin{lemma}  
\label{lem:leaps_aux}
Let us use the notation of Lemma~\ref{lem:leaps}.

{\rm (1)}
For left leaps $u=L_a^{r}(w)$, $r<0$, one  can 
never have a B-failure.  

{\rm (2)}
For the maximal left leap $u = L_a^{r_{\min}}(w)$,
where the entry $a$ goes all the way to the left,
one cannot have an A-failure. 

{\rm (3)}
For the maximal right leap $u = L_a^{r_{\max}}(w)$,
where the entry $a$ goes all the way to the right,
one cannot have an A-failure. 

{\rm (4)}
Let $u=L_a^r(w)$ and $u'=L_a^{r+1}(w)$, for $r\in\Z$, be two 
adjacent leaps such that $a$ is descent-intermediary in $u$ 
(and, thus, $a$ is ascent-intermediary in $u'$). 
Then there is an A-failure in $u$ if and and only 
if there is an A-failure in $u'$.
\end{lemma}

\begin{proof}
(1)
Since $w\in\Sym_n(\B)$, there is
a subset  $I\in \B$ that contains both $b$ and $\max \{w \nw b\}$ 
and such that $I\subset \{w \nw b\}$.
The same subset $I$ works for $u$ because $\{u \nw b\} = \{w \nw b\}$
or $\{u \nw b\}=\{w \nw b\}\cup \{a\}$.

(2) 
In this case, $a$ is greater than all preceding entries in $u$, so
$a = \max \{u \nw a\}$.

(3)  In this case, $a$ is greater than all following entries in $u$.
The interval $I=[a,n]$ contains both $a$ and $\max \{u \nw a\}$,
$I\subset \{u \nw a\}$, and $I\in \B$ because $\B$ is chordal.

(4)
In this case, all entries between the position of $a$ in $u$
and the position of $a$ in $u'$ are less than $a$.
Thus $\{u \nw a\} = \{u' \nw a\}$.
So $u$ has an A-failure if and only if $u'$ has an A-failure.
\end{proof}

\begin{proof}[Proof of Lemma~\ref{lem:leaps}]
It is easier to prove the second part of the lemma.

(2)  By parts (1) and (2) of Lemma~\ref{lem:leaps_aux},
there exists a negative $r$ such that $L_a^r(w)\in \Sym_n(\B)$.
Let us pick such an $r$ with minimal possible absolute value.
Then $r$ should be odd, by part (4) of Lemma~\ref{lem:leaps_aux}, 
which proves (2).

(1)
Suppose that there is an entry $b\ne a$ in the permutation $w$ such that
$b$ and $m=\max \{w \nw b\}$ are in different connected
components of $\B|_{\{w \nw b\} \setminus \{a\}}$.
In this case, $a \in \{w \nw b\}$,
that is $b<a$ and $b$ is located to the right of $a$ in $w$.
(Otherwise, $b$ and $m$ are in different connected components of $\B|_{\{w\nw b\}}$,
which is impossible because $w$ is a $\B$-permutation.)
Let us pick the leftmost entry $b$ in $w$ that satisfies this condition.
Then the permutation $u=L_a^r(w)$ has a B-failure if the letter $a$
moves to the right of this entry $b$; and $u$ has no B-failure if $a$
stays to the left of $b$.  By our assumptions, $a$ stays to the left of $b$ in
$L_a^1(w)$, so such a $u$ exists.

Let $u=L_a^r(w)$ be the maximal right leap (i.e., with maximal 
$r>0$) such that the entry $a$ stays to the left 
of $b$.  Then all entries in $u$ between the positions of $a$ and $b$
should be less than $a$.  Thus 
$m=\max \{u \nw a\} = \max   \{w \nw b\}$.
Since $w\in \Sym_n(\B)$, there is an $I\in \B$ such that $b,m\in I$ and $I\subset \{w \nw b\}$.
This subset $I$ should also contain the entry $a$.
(Otherwise, $b$ and $m$ would be in the same connected component  
$\B|_{\{w \nw b\}\setminus \{a\}}$, contrary to our choice of $b$.)
Thus $I':=I\cap[a,+\infty]\in\B$ contains both $a$ and $m$ 
and $I'\subset \{u \nw a\}$.  This means that there is no A-failure in $u$.  
Thus $u\in\Sym_n(\B)$. 

If there is no entry $b$ in $w$ as above, then none of the permutations
$L_a^r(w)$ has a B-failure.  In this case $L_a^{r_{\max}}(w)\in \Sym_n(\B)$
by  part (3) of Lemma~\ref{lem:leaps_aux}.

In all cases, there exists a positive $r$ such that $L_a^r(w)\in\Sym_n(\B)$
and only A-failures are possible in $L_a^s(w)$, for $0<s<r$.
Let us pick the minimal such $r$.  Then $r$ should be odd by part (4) 
of Lemma~\ref{lem:leaps_aux},
as needed.
\end{proof}

\begin{definition}  Let us define the {\it $\B$-hop\/} operations $\B H_a$.
For a $\B$-permutation $w$ with an ascent-intermediary (resp., descent-intermediary)
entry $a$, the permutation $\B H_a(w)$ is the right leap $u = L_a^r(w)$,
$r>0$ (resp., the left leap $u = L_a^r(w)$, $r<0$) with minimal possible $|r|$ 
such that $u$ is a $\B$-permutation.

Informally, $\B H_a(w)$ is obtained from $w$ by moving the 
node $a$ on its mountain range $M_w$  
directly to the right if $a$ is ascent-intermediary in $w$,
or directly left if $a$ is descent-intermediary in $w$
(possibly passing through several slopes) until one hits a slope
and obtain a $\B$-permutation.
\end{definition}
%

Lemma~\ref{lem:leaps} says that the $\B$-hop $\B H_a(w)$ is 
well-defined for any intermediary entry $a$ in $w$.
It also says that if $a$ is ascent-intermediary in $w$
then $a$ is descent-intermediary in $\B H_a(w)$, and vice versa.
Moreover, according to that lemma, $(\B H_a)^2(w) = w$.

\begin{example}
Let $G$ be the decreasing tree shown on Figure~\ref{fig:behop}. 
Then the graphical building $\B=\B(G)$ is chordal;
see Example~\ref{ex:decreasing_trees}.
Figure~\ref{fig:behop} shows several $\B$-hops 
of the $\B$-permutation $w=(1,10,8,3,6,9,7,4,12,11,5,2)$:
$$
\begin{array}{l}
\B H_1(w) = L_1(w) = (10,8,3,6,9,7,4,12,11,5,2,1),\\
\B H_5(w) = (L_5)^{-5}(w) = (1,5,10,8,3,6,9,7,4,12,11,2),\\
\B H_6(w) = L_6(w) = (1,10,8,3,9,7,6,4,12,11,5,2).
\end{array}
$$
\begin{figure}[ht]
\centering
\includegraphics[height=1.6in]{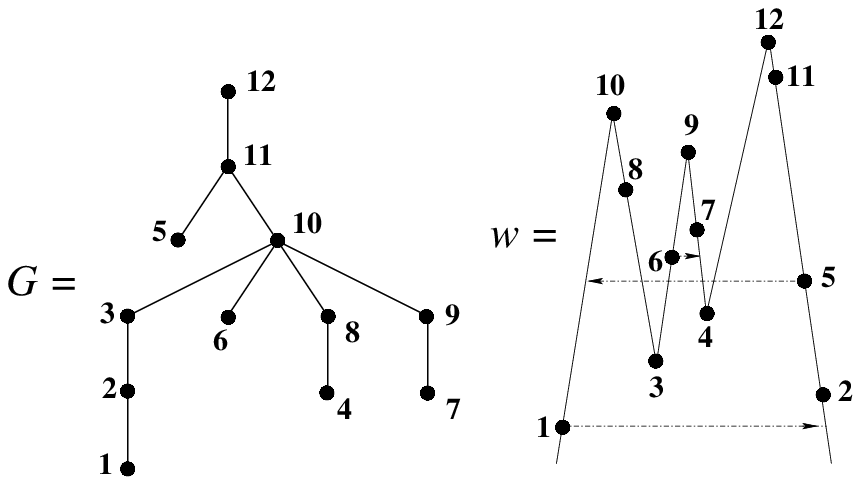}
\caption{A $\B(G)$-permutation $w$ and some $\B$-hops}
\label{fig:behop}
\end{figure}
\end{example}


Let us now show that the $\B$-hop operations pairwise commute with each other.

\begin{lemma}  Let $a$ and $b$ be two intermediary entries in 
a $\B$-permutation $w$.   Then $\B H_a(\B H_b(w))= \B H_b(\B H_a(w))$.
\end{lemma}

\begin{proof}
Let us first assume that both $a$ and $b$ are descent-intermediary 
entries in $w$.  Without loss of generality assume that $a>b$.
In this case $\B H_a(w) = L_a^{r}(w)$
and $\B H_b(w) = L_b^{s}(w)$ for some {\it negative\/} odd $r$ and $s$,
that is the entries $a$ and $b$ of $w$ are moved to the left. 
According to Lemma~\ref{lem:leaps_aux}(1), in this case one 
does not need to worry about B-failures.  In other words,
$\B H_a(w)$ is the first left leap $L_a^r(w)$ (i.e., with minimal $-r>0$) 
that has no A-failure.  Similarly, $\B H_b(w)$ is the first 
left leap $L_b^s(w)$ without A-failures (where A-failures concern the entry $b$).

Since A-failures for permutations $u=L_a^t(w)$, $t<0$, are described in terms of 
the set $\{u\nw a\}\subset [a,\infty]$, moving the entry $b<a$ in $w$
will have no effect on these A-failures.  Thus, for the permutation $w'=\B H_b(w)$,
one has $\B H_a(w') = L_a^r(w')$ with exactly the same $r$ as  
in $\B H_a(w) = L_a^r(w)$.

However, for permutations $u=L_b^t(w)$, $t<0$, the sets $\{u\nw b\}$ might change
if one first performs the operation $\B H_a$ to $w$.  Namely, 
let $\tilde w = \B H_a(w)$  and $\tilde u = L_b^t(\tilde w) = L_b^t(L_a^r(w)).$
Then $\{\tilde u\nw b\} = \{u\nw b\} \cup \{a\}$ if $a$ is located to the left of $b$ 
in $\tilde u$ and $a$ is located to the right of $b$ in $u$
(and $\{\tilde u\nw b\} = \{u\nw b\}$ otherwise).
Notice that one always has $m=\max \{u\nw b\} = \max \{\tilde u \nw b\}$,
since this maximum is the maximal peak preceding $b$ in $u$ (or in $\tilde u$),
and leaps and hops have no affect on the peaks.

If $b$ and $m$ are in the same connected component of 
$\B|_{\{u\nw b\}}$ then they are also in the same connected component of 
$\B|_{\{\tilde u\nw b\}}$, that is if there is no A-failure for $u$
then there is no A-failure for $\tilde u$.

Suppose that there is no A-failure for $\tilde u$ but there is an A-failure
for $u$.  Then the sets $\{u\nw b\}$ and $\{\tilde u\nw b\}$
have to be different.  That means that $a$ is located to the left of $b$
in $\tilde u$ and  $a$ is located to the right of $b$ in $u$.
Let $I$ be the element $I\in \B$ such that $b,m\in I$ and $I\subset \{\tilde u\nw b\}$.
Then $I$ should contain the entry $a$.  (Otherwise, $I \subset \{u\nw b\}$
and there would be no A-failure for $u$.)

Let $\hat w = L_a^{\hat t}(w)$ be the 
left leap with maximal possible $-\hat t\geq 0$ such that the position of $a$ in $\hat w$ 
is located to the right of the position of $b$ in $\tilde u$.
Since $\tilde u = L_b^t(L_a^r(w))$, it follows that 
$|\hat t| < |r|$.
In other words, if one starts moving to the right from the node $b$ along the mountain
range $M_{\tilde u}$, the (ascending) slope that first crosses the level $a$ 
is the place where the entry $a$ is located in $\hat w$.  
Note that 
$\hat t$ is odd because $a$ should be an ascent-intermediary entry in $\hat w$;
in particular $\hat t <0$.

Since all entries in $\hat w$ located between the position of $b$ in $\tilde u$
and the position of $a$ in $\hat w$ are less than $a$, one deduces that 
$\{\tilde u \nw b\}\cap [a,\infty ] = \{\hat w \nw a\}$.
Thus the subset $\hat I = I \cap [a,\infty]$ has 
three important properties: it lies in $\B$ (because $\B$ is chordal); it
contains both $a$ and $m=\max\{\hat w \nw a\}$; 
and it is a subset of $\{\hat w\nw a\}$.
It follows that there is no A-failure in $\hat w$.  
This contradicts the fact that $L_a^r(w)\ne L_a^{\hat t}(w)$ is the first left leap that 
has no A-failure.

Thus $u$ has an A-failure if and only if $\tilde u$ has an A-failure.
It follows that $\B H_b (\tilde w) = L_b^s(\tilde w)$ with exactly the same $s$
as in $\B H_b (w) = L_b^s(w)$.

This proves that $\B H_a (\B H_b(w)) = L_a^r (L_b^s(w)) = L_b^s (L_a^r(w)) = 
\B H_b (\B H_a (w))$, in the case when both $a$ and $b$ are descent-intermediary in $w$.

Let us now show that the general case easily follows.
Suppose that, say, $a$ is ascent-intermediary and $b$ is descent-intermediary
in $w$.  Then, for $w''= \B H_a(w)$ both $a$ and $b$ are descent-intermediary.
One has $\B H_a (\B H_b(w'')) =  \B H_b (\B H_a (w''))$.
Thus $\B H_a (\B H_b( \B H_a (w))) =  \B H_b (\B H_a ( \B H_a(w))) = \B H_b(w)$.
Applying $\B H_a$ to both sides, one deduces 
$\B H_b( \B H_a (w)) =  \B H_a(\B H_b(w))$.  The other cases are similar.
\end{proof}

Thus the $\B$-hop operations $\B H_a$ generate the action of the group 
$(\Z/2\Z)^m$ on the set of $\B$-permutations with a given peak-valley
sequence, where $m$ is the number of intermediary entries in
such permutations.

Let us say that two $\B$-permutation are {\it $\B$-hop-equivalent\/} if they can be
obtained from each other by the $\B$-hop operations $\B H_a$ for various $a$'s.
This gives the partitioning of the set of $\B$-permutations into
$\B$-hop-equivalence classes.

One can now prove Theorem~\ref{th:gamma-chordal} by literally
repeating the argument in the proof of Theorem~\ref{th:CDG}.

\begin{proof}[Proof of Theorem~\ref{th:gamma-chordal}]
For a $\B$-permutation $w\in \Sym_n(\B)$ 
with $p=\peak(w)$,
the descent-generating function of the $\B$-hop-equivalence class
$C$ of $w$ is $\sum_{u\in C} t^{\des(u)} = t^p (t+1)^{n-2p+1}$.
Each $\B$-hop-equivalence class has exactly one representative
without descent-intermediary entries, that is, in the set $\widehat{\Sym}_n(\B)$.
Thus the $h$-polynomial of the nestohedron $P_{\B}$ (see
 Corollary~\ref{cor:h_des_chordal}) is 
$$
h_{P_\B}(t)= \sum_{w\in\Sym_n(\B)} t^{\des(w)} = 
\sum_{w \in \widehat{\Sym}_n(\B)} t^{\peak(w)-1} (t+1)^{n+1-2\,\peak(w)}.
$$
Comparing this to the definition of the $\gamma$-polynomial,
one derives the theorem.
\end{proof}

\begin{corollary}
Gal's conjecture holds for all graph-associahedra corresponding to chordal 
graphs.
\end{corollary}

\subsection{$\gamma$-vectors for the associahedron and cyclohedron}

\begin{proposition}
\label{prop:gamma_assoc}
The $\gamma$-polynomial of the associahedron $P_{\B(\A_n)}$ is
$$
\gamma(t)=\sum_{r=0}^{\lfloor \frac{n-1}{2} \rfloor} C_r  \binom{n-1}{2r}\,t^r,
$$
where $C_r=\frac{1}{r+1}\binom{2r}{r}$ is the $r$-th Catalan number.
\end{proposition}

\begin{proposition}
The $\gamma$-polynomial of the cyclohedron $P_{\B(\cycle_n)}$ is
$$
\gamma(t)=\sum_{r=0}^{\lfloor \frac{n}{2}\rfloor} \binom{n}{r,r,n-2r}\, t^r,
$$
\end{proposition}

These two formulas can be derived from
the expressions 
for the corresponding $h$-polynomials
(see Sections \ref{ssect:path-graph-example-1} and \ref{ssect:cyclo})
using standard quadratic transformations of hypergeometric series;
e.g., see  \cite[Lemma 4.1]{RStantonWelker}.

On the other hand, let us mention the following three combinatorial
interpretations of the $\gamma$-vector 
for the associahedron $P_{\B(\A_n)}$.

\begin{proof}[First proof of Proposition~\ref{prop:gamma_assoc}]
It is known that the Narayana polynomial which is
the $h$-polynomial of $P_{\B(\A_n)}$
is also the rank generating function for the well-studied 
{\it lattice of noncrossing partitions}
$NC(n)$.  An explicit {\it symmetric chain decomposition} for $NC(n)$ was given
by Simion and Ullman \cite{SimionUllman}, who actually produced
a much stronger decomposition of $NC(n)$  into disjoint Boolean
intervals placed symmetrically about the middle rank(s) of $NC(n)$.
Their decomposition contains
exactly $C_r\, \binom{n-1}{2r}$ such Boolean intervals of rank $n-(2r+1)$ for each $r=0,1,\ldots,\frac{n-1}{2}$,
which immediately implies the formula for the $\gamma$-polynomial;
see \cite[Corollary 3.2]{SimionUllman}.
\end{proof}

\begin{proof}[Second proof of Proposition~\ref{prop:gamma_assoc}]
By Section~\ref{ssect:path-graph-example-1}, the $h$-polynomial of
$P_{\B(\A_n)}$ counts plane binary trees on $n$ nodes according to their number
of right edges.
There is a natural map from binary trees to {\it full\/} binary trees, i.e.,
those in which each node has zero or two children: if a node has a unique
child, contract this edge from the node to its child.  If the original binary
tree $T$ has $n$ nodes, then the resulting full binary tree $T'$ will have
$2r+1$ nodes, $2r$ edges and $r$ right edges for some $r=0,1,\dots,
\lfloor(n-1)/2\rfloor$.
There are $C_r$ such full binary trees for each $r$.  Given such a full binary
tree $T'$, one can produce all of the binary trees in its preimage by inserting
$n-(2r+1)$ more nodes and deciding if they create left or right edges.
One chooses the locations of these nodes from $2r+1$ choices, either an edge
of the full binary tree they will subdivide or located above the root, giving
$\binom{n-(2r+1)+ (2r+1)-1}{n-(2r+1)} = \binom{n-1}{2r}$ possible locations.
Thus the generating function with respect to the number of right edges
for the preimage of $T'$ is $\binom{n-1}{2r}\,t^r (t+1)^{n-(2r+1)}$,
where the term $t^r (t+1)^{n-(2r+1)}$ comes from choosing whether each of
the new nodes creates a left or a right edge.
It follows that the generating function for all binary trees on $n$ nodes
is $h_{\A_n}(t) = \sum_r C_r\,\binom{n-1}{2r}\,t^r (t+1)^{n-(2r+1)}$, where
$C_r$ counts full binary trees.  This implies the needed expression
for the $\gamma$-vector of the associahedron $P_{\B(\A_n)}$.

Equivalently, one can describe the subdivision of all binary trees into
classes where two binary trees are in the same class if they can be obtained
from each other by switches of left and right edges
coming from single child nodes.
Then one gets exactly $C_r \binom{n-1}{2r}$ classes
having $t^r (t+1)^{n-(2r+1)}$ as its generating function counting
number of right edges, for each $r=0,1,\dots, \lfloor(n-1)/2\rfloor$.
\end{proof}

\begin{proof}[Third proof of Proposition~\ref{prop:gamma_assoc}]
This proof is based on our general approach to $\gamma$-vectors of chordal nestohedra.
According to Section~\ref{ssect:path-graph-example-1},
$\B$-permutations for the associahedron are $312$-avoiding permutations
and $h$-polynomial is equal to the sum
$h_{P_{\B(\A_n)}}(t) = \sum_{w} q^{\peak(w)-1}$ over all $312$-avoiding
permutations  $w\in \Sym_n$.  
By Theorem~\ref{th:gamma-chordal}, $\gamma_r(P_{\B(\A_n)})$ equals the number of
$312$-avoiding permutations with no descent-intermediary elements
and $r+1$ peaks.
The (flattenings of) peak-valley sequences of such permutations are exactly
$312$-avoiding {\it alternating permutations\/} in $\Sym_{2r+1}$,
that is 312-avoiding permutations $w'$ such that
$w'_1>w'_2<w'_3>\cdots < w'_{2r+1}$.
It is known that the number of such permutations equals
the Catalan number $C_r$; see~\cite[Theorem~2.2]{Mansour}.
Then there are $\binom{n-1}{2r}$ ways to insert
the remaining $n-(2r+1)$ descent-intermediary elements.
\end{proof}

\section{Graph-associahedra for single branched trees}
\label{sec:one-branching}

Our goal in this section is to
compute a generating function that computes the $h$-polynomials of all
graph-associahedra in which the graph is a tree having at most one
{\it branched} vertex (i.e., a vertex of valence $3$ or more).

\subsection{Associahedra and Narayana polynomials}

First recall (see Section~\ref{ssect:path-graph-example-1})
that the $h$-numbers of the associahedron $P_{\B(\A_n)}$
are the Narayana numbers
$h_k(P_{\B(\A_n)}) = N(n,k) := \frac{1}{n} \binom{n}{k} \binom{n}{k-1}$,
and the $h$-polynomial of the associahedron is the 
{\it Narayana polynomial}:
\begin{equation}
\label{eq:Narayanas-as-h-vector}
h_{\B(\A_n)}(t) = C_n(t) := \sum_{k=1}^n N(n,k)\,t^{k-1}.
\end{equation}

Recall the well-known recurrence relation
and the generating function for 
the Narayana polynomials $C_n(t)$.
The recurrence for the $f$-polynomials $f_{\B(\A_n)}(t) =
h_{\B(\A_n)}(t+1)= C_n(t+1)$ given by Theorem~\ref{th:f_B_recurrence}
can be written as follows.
When one removes $k$ vertices from the $n$-path,
it splits into $k+1$ (possibly empty) paths.
Thus one obtains
\begin{equation}
\label{eq:assoc_recurr}
C_n(t) = \sum_{k\geq 1} (t-1)^{k-1} \sum_{m_1+\cdots + m_{k+1} = n-k}
C_{m_1}(t) \cdots C_{m_{k+1}}(t),
\quad\text{for }n\geq 1,
\end{equation}
where the sum is over $m_1,\dots,m_{k+1}\geq 0$
such that $\sum m_i = n-k$.  Here one assumes that $C_0(t)=1$.

Let $C(t,x)$ be the generating function for the
Narayana polynomials:
\begin{eqnarray}
\label{eq:C_formula}
C(t,x) := \sum_{n\geq 1} C_n(t)\, x^n =
x + (1+t)x^2 + (1+3t+t^2)x^3 + \cdots &&\\
\nonumber
= \frac{1-x-tx-\sqrt{(1-x-tx)^2-4tx^2}}{2tx}.&&
\end{eqnarray}
The recurrence relation~\eqref{eq:assoc_recurr}
is equivalent to the following well-known
functional equation:
\begin{equation}
\label{eq:C}
C = tx\, C^2 + (1+t)\,x\, C + x,
\end{equation}
see \cite[Exer.~6.36b]{EC2}.

\subsection{Generating function for single branched trees}

Trees with at most one branched vertex have the following form.
For $a_1,\dots,a_k\geq 0$, let $T_{a_1,\dots,a_k}$ be the graph obtained
by attaching $k$ chains of lengths $a_1,\dots,a_k$ to one central node.
For example, $T_{0,\dots,0}$ is the graph with a single node
and $T_{1,\dots,1}$ is the $k$-star graph $K_{1,k}$.

\begin{theorem}
\label{th:T_formula}
One has the following generating function for the $h$-polynomials of graph-associahedra
$P_{\B(T_{a_1,\dots,a_k})}$ for the graphs $T_{a_1,\dots,a_k}$:
$$
\begin{aligned}
T(t,x_1,\dots,x_k) &:= \sum_{a_1,\dots,a_k\geq 0} h_{T_{a_1,\dots,a_k}}(t)\, x_1^{a_1+1}\cdots x_k^{a_k+1}\\
 &= \frac{ (t-1)\,\phi_1\cdots \phi_k }{t- \prod_{i=1}^k(1+(t-1)\,\phi_i)}
\end{aligned}
$$
where $\phi_i = x_i(1 + t \, C(t,x_i))$,
and $C(t,x)$ is the generating function
for the Narayana polynomials from \eqref{eq:C_formula}.
\end{theorem}

This theorem immediately implies the following formula from~\cite{Post}.

\begin{corollary} {\rm \cite[Proposition~8.7]{Post}} \
The generating function for the number of vertices
in the graph-associahedron $P_{\B(T_{a_1,\dots,a_k})}$ is
$$
\sum_{a_1,\dots,a_k} f_0(P_{\B(T_{a_1,\dots,a_k})})\, x_1^{a_1}\cdots x_k^{a_k}
=
\frac{\bar C(x_1)\cdots \bar C(x_k)}{1-x_1\,\bar C(x_1)-\cdots -
x_k\,\bar C(x_k)},
$$
where $\bar C(x) = \sum_{n\geq 0} C_n\,x^n = \frac{1-\sqrt{1-4x}}{2x}$
is the generating function for the Catalan numbers.
\end{corollary}

\begin{proof}
The claim is obtained from  Theorem~\ref{th:T_formula}
in the limit $t\to 1$.
Note however that one needs to use l'H\^opital's rule before plugging in $t=1$.
\end{proof}

The first proof of Theorem~\ref{th:T_formula}
is fairly direct, using Corollary~\ref{cor:h_des_T}
and the solution to Simon Newcomb's problem.  The second uses Theorem~\ref{th:f_recurrence_zelevinsky}
to set up a system of PDE's and solve them;  it has the advantage of producing a generating function for
the $h$-polynomials of one further family of graph-associahedra.

\subsection{Theorem~\ref{th:T_formula} via Simon Newcomb's problem}
\label{eq:Simon-Newcomb-section}

Let us first review Simon Newcomb's problem and its solution.

Let $w = (w(1), \dots, w_m)$ be a permutation of the multiset
$\{1^{c_1},\dots,k^{c_k}\}$, that is,
each $i$ appears in $w$ exactly $c_i$ times, for $i=1,\dots,k$.
A {\it descent\/} in $w$ is an index $i$ such that $w(i)>w(i+1)$.
Let $\des(w)$ denote the number of descents in $w$.
{\it Simon Newcomb's Problem\/} is the problem
of counting permutations of a multiset
with a given number of descents,
see~\cite[Sec.~IV, Ch.~IV]{MacMahon} and~\cite[Sec.~4.2.13]{GJ}.
Let us define the {\it multiset Eulerian polynomial\/} as
$$
A_{c_1,\dots,c_k}(t) := \sum_{w} t^{\des(w)},
$$
where the sum is over all permutations $w$ of the
multiset $\{1^{c_1},\dots,k^{c_k}\}$.
By convention, set $A_{0,\dots,0}(t) = 1$.

In particular, the polynomial $A_{1,\dots,1}(t)$
is the usual Eulerian polynomial.
It is clear that
$A_{c_1,\dots,c_k}(1)=
\binom{m}{c_1,\dots,c_k}$,
the total number of multiset permutations.
A solution to Simon Newcomb's problem can be expressed
by the following generating function for
the $A_{c_1,\dots,c_k}(t)$.

\begin{proposition}
\label{prop:Simon_Newcomb}
{\rm \cite[Sec.~4.2.13]{GJ}} \
One has
$$
\sum_{c_1,\dots,c_k\geq 0}  A_{c_1,\dots,c_k}(t)\,y_1^{c_1}\cdots y_k^{c_k}
= \frac{t-1} {t-\prod_{i=1}^k (1+(t-1)\,y_i)}.
$$
\end{proposition}

Theorem~\ref{th:T_formula} then immediately follows
from Proposition~\ref{prop:Simon_Newcomb} and the following
proposition.

\begin{proposition}
\label{prop:multi_Eulerian}
The generating function for the $h$-polynomials of the polytopes
$P_{\B(T_{a_1,\dots,a_k})}$ equals
$$
T(t,x_1,\dots,x_k)=\sum_{c_1,\dots,c_k\geq 0} A_{c_1,\dots,c_k}(t)\,
\phi_1^{c_1+1}\cdots \phi_{k}^{c_k+1}.
$$
\end{proposition}

\begin{proof}
Let us label nodes of the graph $T_{a_1,\dots,a_k}$ by integers in $[n]$, where
$n=a_1+\cdots + a_k +1$, so that the first chain is labelled by $1,\dots,a_1$,
the second chain is labelled by $a_1+1,\dots,a_1+a_2$, etc., with all labels
increasing towards the central node, and finally the central
node has the maximal label $n$.

Let $T$ be a $T_{a_1,\dots,a_k}$-tree.  Suppose that the root $r$ of $T$
belongs to the $w(1)$-st chain of the graph $T_{a_1,\dots,a_k}$.
If one removes the node $r$ from
the graph $T_{a_1,\dots,a_k}$, then the graph decomposes into 2 connected
components, one of which is a chain $\A_{b_1}$ and the other is
$T_{a_1,\dots,a_{w(1)}',\dots,a_k}$, where $a_{w(1)}' = a_{w(1)} - b_1 -1$
and all other indices are the same as before.  (The first component is
empty if $b_1=0$.)
According to Proposition~\ref{prop:B_trees_construction}, the tree $T$ is
obtained by attaching a $\A_{b_1}$-tree $T_1$ and a
$T_{a_1,\dots,a_{w(1)}',\dots,a_k}$-tree $T'$ to the root $r$.
(Here one assumes that there is one empty $\A_{0}$-tree $T_1$, for $b_1=0$.)
Let us repeat the same procedure with the tree $T'$.  Assume that
its root belongs to the $w(2)$-nd chain and split it into a $\A_{b_2}$-tree
$T_2$ and a tree $T''$.  Then repeat this procedure with $T''$, etc.
Keep on doing this until one gets a tree $T^{'\cdots'}$ with the root at
the central node $n$.  Finally, if one removes the central node $n$ from
$T^{'\cdots'}$, then it splits into $k$ trees $\tilde T_1,\dots,\tilde T_k$
such that $\tilde T_j$ is a $\A_{d_j}$-tree, for $j=1,\dots,k$.

So each $T_{a_1,\dots,a_k}$-tree $T$ gives us the following data:
\begin{enumerate}
\item a sequence $(w(1),\dots,w_m)\in [k]^m$;
\item a $\A_{b_i}$-tree $T_i$, for $i=1,\dots,m$;
\item a $\A_{d_j}$-tree $\tilde T_j$, for $j=1,\dots,k$.
\end{enumerate}
This data satisfies the following conditions:
$$
\begin{aligned}
m,
b_1,\dots,b_m,
d_1,\dots,d_k  &\geq 0, \text{ and }\\
(b_1+1) e_{w(1)} + \cdots + (b_m+1) e_{w_m}
+ (d_1,\dots,d_k) &= (a_1,\dots,a_k),
\end{aligned}
$$
where $e_1,\dots,e_k$ are the standard basis vectors in $\R^k$.  Conversely,  data of this form
gives us a unique $T_{a_1,\dots,a_k}$-tree $T$.
The number of descents in the tree $T$ is
$$
\des(T)= \sum_{i=1}^m \des(T_i)
+\sum_{j=1}^k \des (\tilde T_j) + l + \des(w),
$$
where $l$ is the number of
nonempty trees among $T_1$, \dots, $T_m$, $\tilde T_1$, \dots, $\tilde T_k$.
Indeed, all descents in trees $T_i$ and $\tilde T_j$ correspond to descents in $T$,
each nonempty tree $T_i$ or $\tilde T_j$ gives an additional descent
for the edge that attaches this tree, and descents in $w$ correspond
to descent edges that attach trees $T', T'',\dots$.

Let us fix a sequence $w=w(1),\dots,w(m)$.  For $i\in[k]$, let $c_i$ be the number
of times the integer $i$ appears in $w$.  In other
words, $w$ is a permutation of the multiset $\{1^{c_1},\dots,k^{c_k}\}$.
Then the total contribution to the generating function $T(t,x_1,\dots,x_k)$
of trees $T$ whose data involve  $w$ is equal to
$t^{\des(w)}\, \phi_1^{c_1+1}\cdots \phi_k^{c_k+1}$.
Indeed, the term 1 in $\phi_i = x_i(1+ t\cdot C(t,x_i))$ corresponds to an empty tree,
and the term $t\cdot C(t,x_i)$ corresponds to nonempty trees, which contribute
one additional
descent.  The term $\phi_i^{c_i}$ comes from the $c_i$ trees
$T_{j_1},\dots,T_{j_{c_i}}$, where $w_{j_1},\dots,w_{j_{c_i}}$ are all occurrences
of $i$ in $w$.  Finally, additional $1$'s in the exponents of $\phi_i$'s come
from the trees $\tilde T_1,\dots,\tilde T_k$.
Summing this expression over all permutations $w$ of the multiset
$\{1^{c_1},\dots,k^{c_k}\}$ and then over all $c_1,\dots,c_k\geq 0$, one obtains
the needed expression for the generating function $T(t,x_1,\dots,x_k)$.

\end{proof}

\begin{remark}
One can dualize all definitions, statements, and arguments in this section,
as follows.  An equivalent dual formulation to Theorem~\ref{th:T_formula} says
$$
T(t,x_1,\dots,x_k) =\frac{ (1-t)\,\psi_1\cdots \psi_k }{ 1- t \,\prod_{i=1}^k (1+(1-t)\,\psi_i)}
$$
where $\psi_i = x_i(1 + C(t,x_i))$.
The equivalence to Theorem~\ref{th:T_formula} follows from the
relation $\phi_i\cdot \psi_i = (t-1)(\phi_i-\psi_i)$,
which is a reformulation of the functional equation~\eqref{eq:C}.

The {\it dual multiset Eulerian polynomial} is
$\bar A_{c_1,\dots,c_k}(t) :=
\sum_{w} t^{\wdes(w)+1}$,
where the sum is over permutations $w$ of the
multiset $M=\{1^{c_1},\dots,k^{c_k}\}$,
$m = c_1+\cdots + c_k$,  and
$\wdes(w)$ is the
number of {\it weak descents\/} in the multiset permutation $w$,
that is, the number of indices $i$ for which $w(i)\geq w(i+1)$.
The bijection which reverses the word $w$ shows that
$\bar A_{c_1, \dots , c_k}(t)= t^m \, A_{c_1,\dots,c_k}(t^{-1})$
and consequently one has an equivalent formulation of the solution to
Simon Newcomb's problem:
$$
\sum_{c_1,\dots,c_k\geq 0}
\bar A_{c_1,\dots,c_k}(t)\,y_1^{c_1}\cdots y_k^{c_k}
= \frac{1-t} {1-t\prod_{i=1}^k (1+(1-t)\,y_i)}.
$$
Then one can modify the proof of Proposition~\ref{prop:multi_Eulerian},
by switching the labels $i \leftrightarrow n+1-i$ in the graph
$T_{a_1,\dots,a_k}$, and applying  a similar argument to show
$$
T(t,x_1,\dots,x_k)
= \sum_{c_1,\dots,c_k\geq 0} \bar A_{c_1,\dots,c_k}(t)\,
\psi_1^{c_1+1}\cdots \psi_{k}^{c_k+1}.
$$
\end{remark}

\subsection{Proof of Theorem~\ref{th:T_formula} via PDE}
\label{sec:PDE}

This section rederives Theorem~\ref{th:T_formula} using
Theorem~\ref{th:f_recurrence_zelevinsky}.  It also calculates 
the generating function for $f$-polynomials of graph-associahedra
corresponding to another class of graphs, the {\it hedgehog graphs} defined
below.


Recall that $\A_n$ is the path with $n$ nodes, and
$T_{a_1,\dots,a_k}$ is the graph obtained by attaching the paths
$\A_{a_{1}}$, \dots, $\A_{a_k}$ to a central
node.  Let us also define the {\it hedgehog graph\/}
$H_{a_1,\dots,a_k}$ as the graph obtained from the disjoint
union of the chains $\A_{a_1}$, \dots, $\A_{a_k}$ by adding edges of the
complete graph between the first vertices of all chains.  For example,
$H_{0,\dots,0}$ is the empty graph, $H_{1,\dots,1} = K_k$,
and $H_{2,\dots,2}$ is a graph with $2k$ vertices obtained from the complete graph $K_k$ by adding a
``leaf'' edge hanging from each of the $k$ original nodes. By convention,
for the empty graph, one has $\tilde f_{H_{0,\dots,0}}(t) = 0$.

Theorem~\ref{th:f_recurrence_zelevinsky} gives the following recurrence
relation for $f$-polynomials of path graphs:
$$
\frac{d}{d t} \tilde f_{\A_n}(t) =
\sum_{r=1}^{n-1} (n-r+1)\cdot \tilde f_{\A_r}(t) \cdot \tilde f_{\A_{n-r}}(t).
$$
Indeed, there are $n-r+1$ connected $r$-element subsets $I$ of
nodes of $\A_n$, the deletion ${\A_n}|_I$ is isomorphic to $\A_r$,
and the contraction $\A_n/I$ is isomorphic to $\A_{n-r}$.

For graphs $T_{a_1,\dots,a_k}$, Theorem~\ref{th:f_recurrence_zelevinsky} gives
the following recurrence relation
$$
\begin{array}{l}
\displaystyle
\frac{d}{dt} \tilde f_{T_{a_1,\dots,a_k}}(t) =
\sum_{i=1}^k \sum_{r=1}^{a_i} \tilde f_{\A_r}(t)\cdot
\tilde f_{T_{a_1,\dots,a_i-r,\dots, a_k}}(t)\cdot (a_i-r+1),\\[.2in]
\displaystyle
\qquad \qquad \qquad \qquad
\qquad \qquad \qquad \qquad
 + \sum \tilde f_{T_{b_1,\dots,b_k}}(t)\cdot \tilde f_{H_{a_1-b_1,\dots,a_k-b_k}}(t),
\end{array}
$$
where the second sum is over $b_1,\dots,b_k$ such that $0\leq b_i\leq a_i$, for $i=1,\dots,k$.
Indeed, a connected subset $I$ of vertices of $G=T_{a_1,\dots,a_k}$
either belongs to one of the chains $\A_{a_i}$,
or contains the central node.
In the first case, the restriction is $G|_I = \A_{r}$
and the contraction is $G/I = T_{a_1,\dots,a_i-r,\dots,a_k}$,
where $r=|I|$.
In the second case, the restriction $G|_I$ has the form
$T_{b_1,\dots,b_k}$ and the contraction is
$G/I = H_{a_1-b_1,\dots,a_k-b_k}$.
Similarly, for hedgehog graphs $H_{a_1,\dots,a_k}$,
one obtains the recurrence relation
$$
\begin{array}{l}
\displaystyle
\frac{d}{dt} \tilde f_{H_{a_1,\dots,a_k}}(t) =
\sum_{i=1}^k \sum_{r=1}^{a_i} \tilde f_{\A_r}(t)\cdot
\tilde f_{H_{a_1,\dots,a_i-r,\dots, a_k}}(t)\cdot (a_i-r),\\[.2in]
\displaystyle
\qquad \qquad \qquad \qquad
\qquad \qquad \qquad \qquad
 + \sum \tilde f_{H_{b_1,\dots,b_k}}(t)\cdot \tilde f_{H_{a_1-b_1,\dots,a_k-b_k}}(t),
\end{array}
$$
where the second sum is over $b_1,\dots,b_k$ such that $0\leq b_i\leq a_i$, for $i=1,\dots,k$.
In all cases one has the initial conditions
$\tilde f_{\A_n}(0) = \tilde f_{T_{a_1,\dots,a_k}}(0)
=\tilde f_{H_{a_1,\dots,a_k}}(0)=1$, except $\tilde f_{\A_0}(t) =
\tilde f_{H_{0,\dots,0}}(t) = 0$.

The above recurrence relations can be written in a more compact form
using these generating functions:
\begin{eqnarray*}
&& F_A(t,x):=\sum_{n\geq 1} \tilde f_{\A_n}(t) \,x^{n+1}
= x^2 + (1+2t)\,x^3 + (1+5t+5t^2)\,x^4+\cdots,\\
&& F_T(t,x_1,\dots,x_k)
:=\sum_{a_1,\dots,a_k\geq 0} \tilde f_{T_{a_1,\dots,a_k}}(t) \,
x_1^{a_1+1}\cdots x_k^{a_k+1}, \\
&& F_H(t,x_1,\dots,x_k):=
\sum_{a_1,\dots,a_k\geq 0} \tilde f_{H_{a_1,\dots,a_k}}(t) \,
x_1^{a_1}\cdots x_k^{a_k}.
\end{eqnarray*}
Note that $F_A$ and $F_T$ are related to generating functions from
Section \ref{sec:one-branching}:
\begin{eqnarray*}
F_A(t,x) &=& t^{-1} x\,C(t^{-1}+1,tx)\\
F_T(t,x_1,\dots,x_k) &=& t^{-k}\, T(t^{-1}+1,tx_1,\dots,t x_k).
\end{eqnarray*}

The above recurrence relations can be expressed as the following
partial differential equations with initial conditions at $t=0$:
\begin{eqnarray}
\label{eq:F_A_diff_eqn}
&&\frac{\partial F_{A}}{\partial t} = F_A \cdot \frac{\partial F_A}{\partial x},
\quad F_A|_{t=0} = \frac {x^2}{1-x},\\
\label{eq:F_T_diff_eqn}
&&\frac{\partial F_T}{\partial t} =
\sum_{i=1}^k F_A(t,x_i) \, \frac{\partial F_T}{\partial x_i} +
F_T\cdot F_H,
\quad  F_T|_{t=0} = \frac{x_1\cdots x_k }{\prod_{i=1}^k(1-x_i)},\\
\label{eq:F_H_diff_eqn}
&&\frac{\partial F_H}{\partial t} =
\sum_{i=1}^k F_A(t,x_i)\,\frac{\partial F_H}{\partial x_i} +  (F_H)^2,
\quad F_H|_{t=0} = \frac{1-\prod_{i=1}^k(1-x_i)}{\prod_{i=1}^k (1-x_i)}.
\end{eqnarray}

One can actually solve these partial differential equations for
arbitrary initial conditions, as follows.

\begin{proposition}
\label{eq:diffeq_FGHR}
The solutions $F(t,x)$,
$G(t,x_1,\dots,x_k)$,
$H(t,x_1, \dots,x_k)$,
and $R(t,x_1,\dots,x_k)$ to
the following system of partial differential equations
with initial conditions
\begin{eqnarray}
\label{eq:dFdq}
&& \frac{\partial F}{\partial t} = F \cdot \frac{\partial F}{\partial x},
\quad F|_{t=0} = f_0(x),\\
\label{eq:dGdq}
&& \frac{\partial G}{\partial t} =
\sum_{i=1}^k F(t,x_i)\,\frac{\partial G}{\partial x_i},
\quad G|_{t=0} = g_0(x_1,\dots,x_k),\\
\label{eq:dHdq}
&& \frac{\partial H}{\partial t} =
\sum_{i=1}^k F(t,x_i)\,\frac{\partial H}{\partial x_i} + H^2,
\quad H|_{t=0} = h_0(x_1,\dots,x_k),\\
\label{eq:dRdq}
&& \frac{\partial R}{\partial t} =
\sum_{i=1}^k F(t,x_i)\,\frac{\partial R}{\partial x_i} + R\cdot H,
\quad R|_{t=0} = r_0(x_1,\dots,x_k)
\end{eqnarray}
are given by
\begin{eqnarray*}
&&f_0(x + t\cdot F) = F \text{ (implicit form)}\\
&&G=g_0(\xi_1,\dots,\xi_k)\\
&&H = - (t +  (h_0(\xi_1,\dots,\xi_k))^{-1})^{-1}\\
&&R = - r_0(\xi_1,\dots,\xi_k)\cdot (1 + t\cdot h_0(\xi_1,\cdots,\xi_k))^{-1}
\end{eqnarray*}
where  $\xi_i = x_i + t\cdot F(t,x_i)$, for $i=1,\dots,k$.
\end{proposition}

\begin{proof}
Let us first solve~(\ref{eq:dFdq}).
For a constant $C$, consider the function $x(t)$ given implicitly as
$F(t,x)= C$, i.e., the graph of $x(t)$ is a level curve for $F(t,x)$.
The tangent vector to the graph of $x(t)$ at some point
$(t_0,x_0)$ such that
$F(t_0,x_0)=C$ is $(1,\frac{dx(t_0)}{dt})$.
The derivative of
the function $F(t,x)$ at the point $(t_0,x_0)$ in
the direction of this vector should be 0, i.e.,
$1\cdot \frac {\partial F(t_0,x_0)}{\partial t} +
\frac {dx(t_0)}{dt} \cdot \frac{\partial F(t_0,x_0)}{\partial x} =0$.
This equation, together with the differential equation (\ref{eq:dFdq})
for $F$, implies that
$\frac{d}{dt}\,x(t) = -C$.
Solving this trivial differential equation for $x(t)$ one deduces that
$x(t) = - C\cdot t + B(C)$, where $B$ is a function that depends only on
the constant $C$.
Since $C$ can be an arbitrary constant, one deduces that
$$
x = - F(t,x) \cdot t + B (F(t,x)),
\textrm{ or, equivalently, }
B^{\left<-1\right>}(x+t\cdot F(t,x)) = F(t,x).
$$
Plugging the initial condition $F|_{t=0}= f_0(x)$ in the last expression,
one gets
$$
B^{\left<-1\right>}(x) = f_0(x).
$$
Thus the solution $F(t,x)$ is given by
$f_0(x+t\cdot F) = F$, as needed.

Direct verification shows that the function
$G=R(F(t,x_1),\dots,F(t,x_k))$
satisfies the differential equation~(\ref{eq:dGdq}),
for an arbitrary $R(y_1,\dots,y_k)$.
The initial condition for $t=0$ gives
$R(f_0(x_1),\dots,f_0(x_k)) = g_0(x_1,\dots,x_k)$.
Thus $R(y_1,\dots,y_k) = g_0(B(y_1),\dots,B(y_k))$,
where $B = f_0^{\left<-1\right>}$, as above.
Since $B(F(t,x)) = x + t\cdot F(t,x)$,
one deduces that
$G= g_0(B(F(t,x_1)),\dots, B(F(t,x_k))) =
g_0(\xi_1,\dots,\xi_k)$, as needed.

Making the substitution $H = - (t + G(t,x_1,\dots,x_k))^{-1}$
in differential equation (\ref{eq:dHdq}) for $H$, one obtains
equation (\ref{eq:dHdq}) for $G$ with $g_0 = - (h_0)^{-1}$.
By the previous calculation, one has $G = - (h_0(\xi_1,\dots,\xi_k))^{-1}$.
Thus the solution for (\ref{eq:dHdq}) is
$H = - (t +  (h_0(\xi_1,\dots,\xi_k))^{-1})^{-1}$.

Making the substitution $R = H\cdot G$
in equation~(\ref{eq:dHdq}) for $R$,
where $H$ is the solution to (\ref{eq:dHdq}),
one obtains equation (\ref{eq:dGdq}) for $G$
with $g_0 = r_0/h_0$.
By the above calculation, one has
$G = r_0(\xi_1,\dots,x_k)/h_0(\xi_1,\dots,\xi_k)$.  Thus,
$$
R= - \frac{1}{t +  (h_0(\xi_1,\dots,\xi_k))^{-1}}\cdot
\frac{r_0(\xi_1,\dots,x_k)}{h_0(\xi_1,\dots,\xi_k)}
= - \frac{ r_0(\xi_1,\dots,x_k)}{1 + t\cdot h_0(\xi_1,\cdots,\xi_k)},
$$
as needed.
\end{proof}

Applying Proposition~\ref{eq:diffeq_FGHR} to differential
equation~(\ref{eq:F_A_diff_eqn}) for $F_A(t,x)$,
one obtains the implicit solution:
$$
\frac{(x+t\cdot F_A)^2}{1-x- t\cdot F_A}  = F_A.
$$
This is equivalent to the quadratic equation~(\ref{eq:C})
for $C(t,x)$.
Explicitly, one gets
\begin{equation}
\label{eq:expression_for_FA}
F_A(t,x) = \frac{(1-x-2tx) - \sqrt{(1-x-2tx)^2 - 4 t(t+1)x^2}}{2 t(t+1)}.
\end{equation}

Applying Proposition~\ref{eq:diffeq_FGHR} to differential
equations (\ref{eq:F_T_diff_eqn}) and (\ref{eq:F_H_diff_eqn}) for
the generating functions $F_T$ and $F_H$, one obtains
the following result.

\begin{theorem}
\label{th:FH_FT}
The generating functions $F_T(t,x_1,\dots,x_k)$ and
$F_H(t,x_1,\dots,x_k)$ are given by the following expressions
\begin{eqnarray*}
&&F_T(t,x_1,\dots,x_k) =
\frac{\xi_1\cdots \xi_n} {(t+1)(1-\xi_1)\cdots (1-\xi_n) - t},\\
&&F_H(t,x_1,\dots,x_k) =
\frac{1-(1-\xi_1)\cdots (1-\xi_k)} {(t+1)(1-\xi_1)\cdots (1-\xi_n) - t},
\end{eqnarray*}
where $\xi_i = x_i + t\cdot F_A(t,x_i)$
and $F_A$ is given by~{\rm(\ref{eq:expression_for_FA}).}
\end{theorem}

Note that the above expression for $F_T$ is equivalent to
Theorem~\ref{th:T_formula}, using~\eqref{eq:f-h-relation}.

\section{Graph-associahedra for path-like graphs}
\label{sec:almost-paths}

The goal of this section is to use Theorem~\ref{th:f_B_recurrence} to compute
the $h$-polynomials of the graph-associahedra
of a fairly general infinite family of graphs,
including all Dynkin diagrams of finite and affine Coxeter groups.

Let $A$ and $B$ be two connected graphs with a marked vertex in each,
and let $n_0$ be the total number of unmarked vertices in $A$ and $B$.
For $n> n_0$, let $G_n = G_n(A,B)$ be the graph obtained by
connecting the marked vertices in $A$ and $B$ by the path $\A_{n-n_0}$
so that the total number of vertices in $G_n$ is $n$.
Call graphs of the form $G_n$ {\it path-like graphs\/}
because, for large $n$, they look like paths with some ``small'' graphs
attached to their ends.

The following claim shows that the $h$-polynomials
of the graph-associahedra $P_{\B(G_n)}$ can be expressed as
linear combinations (with polynomial coefficients)
of the $h$-polynomials $C_n(t)$ of usual associahedra;
see~\eqref{eq:Narayanas-as-h-vector}.

\begin{theorem}
\label{th:almost_path_combination}
  There exist unique polynomials
$g_0(t),g_1(t),\dots,g_{n_0}(t)\in\Z[t]$
of degrees $\deg g_i(t) \leq i$
such that, for any $n> n_0$, one has
$$
h_{G_n}(t) = g_0(t) \, C_n(t) +
g_1(t) \, C_{n-1}(t) +\cdots + g_{n_0}(t)\, C_{n-n_0}(t).
$$
The polynomial $g_i(t)$ is a palindromic polynomial,
that is $g_i(t) = t^i\,g_i(t^{-1})$, for $i=0,\dots,n_0$.
\end{theorem}

Similarly, one can express the $f$-polynomials of $P_{\B(G_n)}$
as a linear combination of the $f$-polynomials of usual associahedra,
because $f_G(t) = h_G(t+1)$.

One can rewrite this theorem in terms the generating function
$C(t,x)$ for the Narayana numbers; see~\eqref{eq:C_formula}.

\begin{corollary}
\label{cor:path_like}
There exists a unique polynomial $g(t,x)\in\Z[t,x]$
such that, for any $n>n_0$, one has
$$
h_{G_n}(t) = \text{the coefficient of $x^n$ in } g(t,x) \,C(t,x).
$$
The polynomial $g(t,x)$ has degree $\leq n_0$ with respect
to the variable $x$.  It satisfies the equation $g(t,x) = g(t^{-1},tx)$.
\end{corollary}

\begin{proof}
This follows from Theorem~\ref{th:almost_path_combination}, by
setting $g(t,x) = g_0(t) + g_1(t)\,x + \cdots +g_{n_0}(t)\, x^{n_0}$.
\end{proof}

\begin{proof}[Proof of Theorem~\ref{th:almost_path_combination}]
Let us first prove the existence of the linear expansion (and later prove
its uniqueness and the palindromic property of the coefficients $g_i(t)$).  The recurrence from
Theorem~\ref{th:f_B_recurrence} will be used to prove this clam by induction
on the total number of vertices in $A$ and $B$.
For this argument one should drop the assumption that $A$ and $B$ are connected.
Suppose that $A$ or $B$ is disconnected, say, $A$ is a disjoint union of
graphs $A_1$ and $A_2$ where
$A_1$ contains the marked vertex.  Let $\tilde G_n := G_n(A_1,B)$ and let $r$
be the number of vertices in $A_2$.  Then $h_{G_n}(t) =
h_{\tilde G_{n-r}}(t)\, h_{A_2}(t)$, where $\deg h_{A_2}(t) \leq r-1$.
By induction, $h_{\tilde G_{n-r}}(t)$
can be expressed as a linear combination of $C_{n-r}(t), C_{n-r-1}(t),
\dots, C_{n-n_0}(t)$, which produces the needed expression for $h_{G_n}(t)$.

Now assume that both $A$ and $B$ are connected graphs.
Theorem~\ref{th:f_B_recurrence}(3) gives the  expression
for the $h$-polynomial as the sum
$h_{G_n}(t) = \sum_{L} (t-1)^{|L|-1} h_{G_n\setminus L}(t)$
over nonempty subsets $L$ of vertices of $G_n$, where
$G_n\setminus L$ denotes the graph $G_n$ with removed vertices in $L$.
(Here one has shifted $t$ by $-1$ to transform $f$-polynomials into
$h$-polynomials.)
Let us write $L$ as a disjoint union $L=I\cup J\cup K$,
where $I$ is a subset of unmarked vertices of $A$,
$J$ is a subset of unmarked vertices of $B$,
and $K$ is a subset of vertices in the path connecting the marked vertices.
The contribution of the terms with $K=\emptyset$ to the above sum is
$\sum_{I,J} (t-1)^{|I|+|J|-1} h_{G_n \setminus (I\cup J)}(t)$.
Note that $G_n \setminus (I\cup J) = G_{n-r}(A\setminus I,B\setminus J)$,
where $r=|I|+|J|$.  By induction, one can express each term
$h_{G_n \setminus (I\cup J)}(t)$ as a combination of
$C_{n-r}(t),\dots, C_{n-n_0}(t)$.

The remaining terms involve a nonempty subset $K$ of vertices
in the path $\A_{n-n_0}$.  When one removes these $k=|K|$ vertices from the path,
it splits into $k+1$ smaller paths $\A_{m_1},\dots,\A_{m_{k+1}}$
with $m_i\geq 0$; cf.\ paragraph before~\eqref{eq:assoc_recurr}.
Thus the remaining contribution to $h_{G_n}(t)$ can be written as
$$
\sum_{I,J} \sum_{m_1,\dots,m_{k+1}\geq 0}
(t-1)^{|I|+|J|+k-1}\,
h_{G_p(A\setminus I,\circ)}(t)\,
C_{m_2}(t) \cdots C_{m_k}(t) \, h_{G_q(\circ,B\setminus J)}(t),
$$
where $\circ$ is the graph with a single vertex,
$$
\begin{aligned}
p &= m_1 + |A\setminus I| -1,\\
q &= m_{k+1} + |B\setminus J| -1,\text{ and }\\
k+\sum m_i &= n-n_0.
\end{aligned}
$$
By induction, one can express
$h_{G_p(A\setminus I,\circ)}(t)$
and $h_{G_q(\circ, B\setminus J)}(t)$ as linear combinations of
the $C_{m'}(t)$.  So the remaining contribution to $h_{G_n}(t)$
can be written as a sum of several expressions of the form
$$
s(t)\ \sum_{k\geq 1}
\sum_{m_1',m_2,\dots,m_k,m_{k+1}'}
(t-1)^{k-1}
C_{m_1'}(t)\, C_{m_2}(t) \cdots C_{m_k}(t)\, C_{m_{k+1}'}(t),
$$
where the sum is over $m_1',m_2,\dots,m_k,m_{k+1}'$ such that
$m_1'\geq a$, $m_2,\dots,m_k\geq 0$, $m_{k+1}'\geq b$,
$m_1'+m_2+\cdots + m_k + m_{k+1}' + k = n - c$.
This expression depends on nonnegative integers $a,b,c$ such
that $a+b+c = n_0$ and a polynomial $s(t)$
of degree $\deg s(t) \leq c$.
If one extends the summation to all $m_1',m_{k+1}'\geq 0$,
one obtains the expansion \eqref{eq:assoc_recurr} for $C_{n-c}(t)$
times $s(t)$.
Applying the inclusion-exclusion principle
and equation~\eqref{eq:assoc_recurr}, one deduces that the previous sum
equals
$$
s(t)\,\left(C_{n-c}(t)- \sum_{m_1'=0}^{a-1}
t\,C_{m_1'}(t)\,C_{n-c-m_1'-1}(t)-\cdots
\right),
$$
which is a combination of $C_n(t),\dots,C_{n-n_0}(t)$
as needed.

It remains to show the uniqueness of the linear expansion
and show that the $g_i(t)$ are palindromic
polynomials.  (Here one assumes that the graphs $A$ and $B$ are connected.)
According to  Corollary~\ref{cor:path_like},
the polynomial $H(t,x):=\sum_{n>n_0} h_{G_n}(t)\, x^n$
can be written as $H(t,x) = g(t,x) \,C(t,x) + r(t,x)$,
where $g(t,x), r(t,x)\in\Z[t,x]$.
If $H(t,x) = \tilde g(t,x) \,C(t,x) + \tilde r(t,x)$
with $\tilde g(t,x) \ne g(t,x)$, then this would imply
that $C(t,x)$ is a rational function, which contradicts the 
formula~\eqref{eq:C_formula} involving a square root.
This proves the uniqueness claim.
One has $H(t,x) = H(t^{-1},tx)/t$ and $C(t,x) = C(t^{-1},tx)/t$
because $h$-polynomials are palindromic.
Thus
$$
H(t,x) = g(t,x)\,C(t,x) + r(t,x) =
g(t^{-1},tx) \,C(t,x) + \frac{1}{t}r(t^{-1},tx).
$$
This implies that $g(t,x) = g(t^{-1},tx)$.  Otherwise, $C(t,x)$
would be a rational function.  The equation $g(t,x) = g(t^{-1},tx)$
says that the coefficients $g_i(t)$ of $g(t,x) = \sum_i g_i(t)\,x^i$
are palindromic.
\end{proof}

Let us illustrate Theorem~\ref{th:almost_path_combination}
by several examples.
For a series of path-like graphs $G_n$, let $g\{G_n\}$ denote the polynomial
$g(t,x)$ that appears in the generating functions
$\sum_{n\geq n_0} h_{G_n}(t) \,x^n = g(t,x)\,C(t,x) + r(t,x)$.
For instance, the expression $g\{D_n\} = 2 - (t+1)x -tx^2$ (see the example
below) is equivalent to the expression
$h_{D_n}(t) = 2C_n(t) - (t+1)\,C_{n-1}(t) - t\,C_{n-2}(t)$,
for $n>2$.

\begin{examples}
Define {\it daisy graphs\/} as $\mathrm{Daisy}_{n,k} :=
T_{n-k-1,1^k}$;
see Section~\ref{sec:one-branching}.  (Here $1^k$ means a sequence
of $k$ ones.) They include {\it type $D$ Dynkin diagrams\/} $D_n :=
\mathrm{Daisy}_{n,2}$. For fixed $k$, the $\mathrm{Daisy}_{n,k}$
form the series of path-like graphs for $A = K_{1,k}$ (the $k$-star
with marked central vertex) and $B=\circ$ (the graph with a single
vertex). Also define {\it kite graphs\/} as $\mathrm{Kite}_{n,k} :=
H_{n-k+1,1^{k-1}}$; see Section~\ref{sec:PDE}. They are path-like
graphs for $A=K_k$ and $B=\circ$. The {\it affine Dynkin diagram of
type $\widetilde{D}_{n-1}$} is the $n$th path-like graph in the case
when both $A$ and $B$ are $3$-paths with marked central vertices.

Here are the polynomials $g\{G_n\}$ for
several series of such graphs:
\begin{eqnarray*}
&&g\{D_n\} =  2-(t+1)\,x - t\,x^2, \\
&&g\{\widetilde{D}_{n-1}\} = 4 - 4(t+1)\,x + (t-1)^2\,x^2 +
2\, t(t+1)\,x^3 + t^2\,x^4,\\
&&g\{\Kite_{n,3}\} = 2 - (t+1)\,x,\\
&&g\{\mathrm{Daisy}_{n,3}\}=6-6(t+1)\,x+(1-5\,t+t^2)\,x^2-t(t+1)\,x^3,\\
&& g\{\mathrm{Daisy}_{n,4}\} =24-36(t+1)\,x+(14-16\,t+14\,t^2)\,x^2 + \\
&& \qquad\qquad\qquad {}+ (-1+3\,t+3\,t^2-t^3)\,x^3 - (t+t^2+t^3)\,x^4.
\end{eqnarray*}
The formulas for $D_n$, $\Kite_{n,k}$, and $\mathrm{Daisy}_{n,k}$
were derived from Theorems~\ref{th:T_formula} and~\ref{th:FH_FT}
(using the {\sc Maple} package). The formula for the affine Dynkin
diagram $\widetilde{D}_{n-1}$ was obtained using the inductive
procedure given in the proof of
Theorem~\ref{th:almost_path_combination}.
\end{examples}

%
%
%
%
%
%

\begin{remark}
The induction from the proof of
Theorem~\ref{th:almost_path_combination} is quite involved.  It is
very difficult to calculate by hand other examples for bigger graphs
$A$ and $B$.  It would be interesting to find a simpler procedure
for finding the polynomials $g\{G_n\}$. Also it would be interesting
to find explicit formulas for the polynomials $g\{G_n\}$ for all
daisy graphs, kite graphs, and other ``natural'' series of path-like
graphs.
\end{remark}

\section{Bounds and monotonicity for face numbers of  graph-associahedra}
\label{sec:flossing}

Section~\ref{sec:Stanley-Pitman-example} showed that the $f$- and $h$-vectors of flag nestohedra
coming from connected building sets are bounded below by those of hypercubes and
bounded above by those of permutohedra.
It is natural to ask for the bounds within the subclass of graph-associahedra corresponding to connected
graphs, and to ask for bounds on their $\gamma$-vectors.

Permutohedra are graph-associahedra corresponding to complete graphs, and so
they still provide the upper bound for the $f$- and $h$-vectors.
For lower bounds on $f$- and $h$-vectors, the monotonicity discussed
in Remark~\ref{rem:graphical-monotonicity} implies that
the $f$- and $h$-vector of any graph-associahedron $P_{\B(G)}$ for a 
connected graph $G$ is bounded
below by the graph-associahedron for any spanning tree inside $G$.  
Thus it is of interest to
look at lower (and upper) bounds for $f$-, $h$- and $\gamma$-vectors of 
graph-associahedra for trees.

\begin{figure}[ht]
\centering
\includegraphics[height=2.5in]{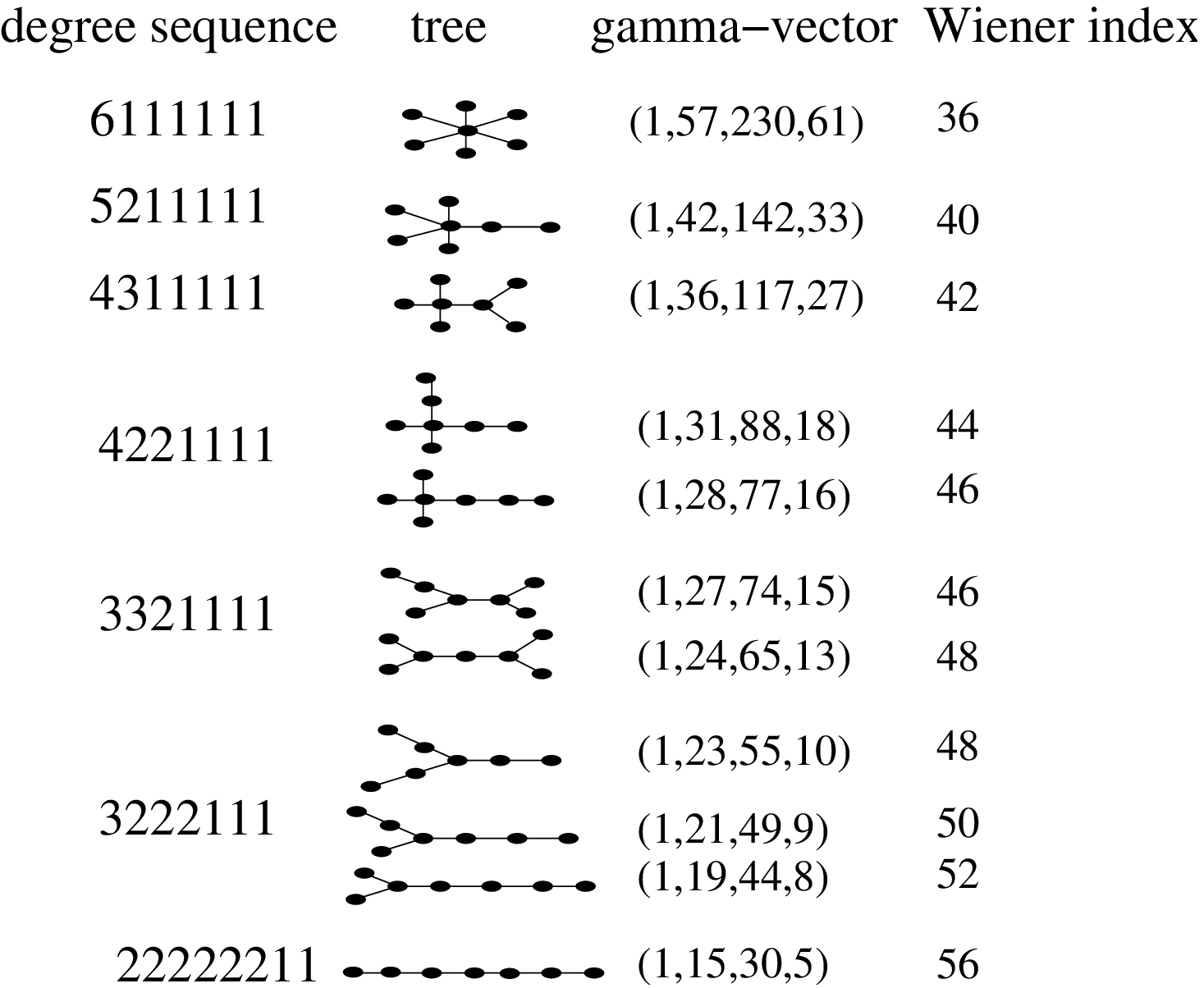}
\caption{The $\gamma$-vectors $(\gamma_0,\gamma_1,\gamma_2,\gamma_3)$ for graph-associahedra
of all trees on $7$ vertices, grouped according to their degree sequences.}
\label{fig:trees-on-seven-gammas}
\end{figure}

A glance at Figure~\ref{fig:trees-on-seven-gammas} suggests that, roughly speaking,
trees which are more branched and forked (that is, farther from being a path)
tend to have higher entries componentwise in their $\gamma$-vectors, and hence also
in their $f$- and $h$-vectors.
In fact, in that figure, which shows all trees on $7$ vertices
grouped by their degree sequences, one sees several (perhaps misleading) features:

\begin{enumerate}
\item[(i)] The degree sequences are ordered linearly under the {\it majorization} (or {\it dominance})
partial ordering on partitions of $2(n-1)$  ($= 2 \cdot (7-1) =12$ here).
\item[(ii)]  The $\gamma$-vectors of these trees are linearly ordered under the componentwise
partial order.
\item[(iii)]  Trees whose degree sequence are lower in the majorization order have
componentwise smaller $\gamma$-vectors.
\item[(iv)]  The trees are distinguished up to isomorphism by their $\gamma$-vectors.
\end{enumerate}

Additionally, it
seems that the {\it Wiener index} \cite{Wiener}
for graphs has some correlation with
the $\gamma$-vector.
The {\it Wiener index\/} $W(G)$ of a graph $G$ is defined as the sum of
distances $d(i,j)$ over unordered pairs $i,j$ of vertices in the graph $G$,
where $d(i,j)$ is the number of edges in the shortest path from $i$ to $j$.
The Wiener index $W(T)$ of a tree is equal to the number of forbidden
$312$ patterns (see the remarks following Definition 
\ref{def:chordal_building})
provided by the tree $T$ (plus the constant $\binom{n}{2}$).
Thus, for two trees on $n$ vertices, if one has $W(T) < W(T')$, then
roughly speaking one might expect that the generalized permutohedron $P_{\B(T)}$
has a larger gamma vector than $P_{\B(T')}$.

This is exactly the case for trees on $7$ vertices, as shown in
Figure~\ref{fig:trees-on-seven-gammas}.
It shows that as the $\gamma$-vectors decrease, the Wiener indices (weakly)
increase.
Note that in this case, the Wiener index together with the degree sequence
completely
distinguish all equivalence classes of trees.

For trees $T$ on $n$ vertices, the maximum and minimum values of the Wiener index
are, respectively,
$\sum_{i=1}^{n-1} i(n-i) = n(n^2-1)/6$ for $\A_n$, and
$(n-1)^2$ for $K_{1,n-1}$.

\begin{figure}[ht]
\centering
\includegraphics[height=2in]{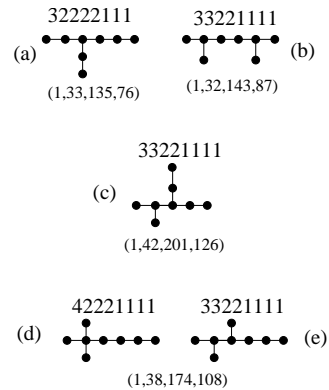}
\caption{Some trees on $8$ nodes and the $\gamma$-vectors of their
graph-associahedra.}
\label{fig:eight-node-ceg-gammas}
\end{figure}

None of the properties
(i)-(iv) persist for all trees. For example, when looking at trees
on $n=8$ nodes, one finds that
\begin{enumerate}
\item[(i)]  the degree sequences are only partially ordered by
the majorization order on partitions of $2(n-1)=14$:
$$
\begin{aligned}
22222211 &< 32222111 \\
         &< 33221111 < 33311111, 42221111 \\
         & < 43211111 < 44111111, 52211111 \\
         & < 53111111 < 62111111 < 71111111
\end{aligned}
$$
\item[(ii)] there are trees, such as the two shown in Figure~\ref{fig:eight-node-ceg-gammas}(a) and (b),
whose $\gamma$-vectors are incomparable componentwise,
\item[(iii)] there are trees, such as the two shown in Figure~\ref{fig:eight-node-ceg-gammas}(d) and (c),
where the degree sequence of one strictly majorizes that of the other, but its $\gamma$-vector is
strictly smaller, and
\item[(iv)] there are nonisomorphic trees, such as the two shown in
Figure~\ref{fig:eight-node-ceg-gammas}(d) and (e), having the same $\gamma$-vector.
\end{enumerate}

Nevertheless, we do make some monotonicity conjectures about the face
numbers for graph-associahedra.

\begin{conjecture}
\label{conj:flossing-order-conjecture}
There exists a partial order $\prec$ on the set of (unlabelled, isomorphism classes of) trees with $n$ nodes,
having these properties:
\begin{enumerate}
\item[$\bullet$]
$\A_n$ is the unique $\prec$-minimum element,
\item[$\bullet$]
$K_{1,n-1}$ is the unique $\prec$-maximum element, and
\item[$\bullet$]
$T \prec T'$ implies $\gamma_{P_{\B(T)}} \leq \gamma_{P_{\B(T')}}$ componentwise.
\end{enumerate}
\end{conjecture}

\noindent
We suspect that such a partial order $\prec$ can be defined so that
$T,T'$ will, in particular, be comparable whenever $T, T'$ are related by one of the {\it flossing moves}
considered in \cite[\S4.2]{BabsonR}.

Assuming Conjecture~\ref{conj:flossing-order-conjecture}, the $\gamma$-vectors
(and hence also the $f$-, $h$-vectors) of graph-associahedra for trees on $n$ nodes
would have the associahedron $P_{\B(\A_n)}$ and the stellohedron $P_{\B(K_{1,n-1})}$
giving their lower and upper bounds.  This would also imply that the $f$-, $h$-vectors of graph-associahedra
for connected graphs on $n$ nodes would have associahedra and permutohedra giving their lower and upper bounds.
To make a similar assertion for $\gamma$-vectors it would be nice to have
the following analogue of Stanley's monotonicity result (Theorem~\ref{th:subdivision-monotonicity}).

\begin{conjecture}
\label{conj:gamma-monotonicity}
When $\Delta, \Delta'$ are two flag simplicial complexes and $\Delta'$ is a
geometric subdivision of $\Delta$, the $\gamma$-vector of $\Delta'$ is
componentwise weakly larger than that of $\Delta$.

In particular, when $\B, \B'$ are building sets giving rise to flag nestohedra
and $\B \subset \B'$, (such as graphical buildings $\B(G) \subset \B(G')$ for
an edge-subgraph $G \subset G'$) then the $\gamma$-vector of $P_{\B'}$ is
componentwise weakly larger than that of $P_{\B}$.
\end{conjecture}

We close with a question suggested by the sets of permutations 
$\Sym_n(\B)$ and $\widehat{\Sym}_n(\B)$ for a chordal building set $\B$ which appeared in
Corollary~\ref{cor:h_des_chordal} and Theorem~\ref{th:gamma-chordal}.

\begin{question}
Given a (non-chordal) building set $\B$,
is there a way to define two sets of permutations $\Sym_n'(\B)$
and $\widehat{\Sym}_n'(\B)$ such that:
\begin{enumerate}
\item[$\bullet$]
the $h$-polynomial for the nestohedron $P_\B$ is given by
the descent generating function for $\Sym_n'(\B)$, and
\item[$\bullet$]
the $\gamma$-polynomial is given by the peak generating function
for $\widehat{\Sym}_n'(\B)$?
\end{enumerate}
\end{question}

\section{Appendix:  Deformations of a simple polytope}
\label{sec:deformations}

The goal of this section is to clarify the equivalence between
various definitions of the deformations of a simple polytope,
either by
\begin{enumerate}
\item[$\bullet$] deforming vertex positions, keeping edges in the same parallelism class, or
\item[$\bullet$] deforming edge lengths, keeping them nonnegative, or
\item[$\bullet$] altering level sets of facet inequalities, or
but not allowing facets to move past any vertices.
\end{enumerate}
There will be defined below three cones of such
deformations which are all linearly isomorphic.
This discussion is essentially implicit in
\cite[Definition 6.1 and \S19]{Post}, but we hope the explication
here clarifies this relationship.

Let $P$ be a simple $d$-dimensional polytope in $\R^d$.  Let $V$ be its set of
vertices.
Let $E  \subseteq
\binom{V}{2}$ be its set of edge pairs.  Let $F$ be an indexing set for its
facets, so that $P$ is defined by facet inequalities $h_f(x) \leq \alpha_f$ for
$f \in F$, in which each $h_f$ is a linear functional in $(\R^d)^*$, and
$(\alpha_f)_{f \in F} \in \R^F$.

\begin{definition}
\label{def:three_cones}
(1) The {\it vertex deformation cone\/} $D^V_P$ of $P$ is the set
of points $(x_v)_{v\in V} \in(\R^d)^V$ such that
\begin{equation}
\label{eq:nonnegative-scaling_uv}
x_u - x_v \in \R_{\geq 0} (u-v),\quad
\text{for every edge $uv \in E$.}
\end{equation}

(2) The {\it edge length deformation cone\/} $D^E_P$ of $P$
is the set of points $(y_e)_{e\in E}\in \R^E$ such that all $y_e\geq 0$,
and, for any 2-dimensional face of $P$ with edges
$e_1 = v_1 v_2$, $e_2 = v_2 v_3$, \dots, $e_{k-1} = v_{k-1} v_k$,
$e_k = v_k v_1$, one has
$$
y_{e_1}(v_1-v_2) + y_{e_2} (v_2 - v_3) + \cdots + y_{e_k} (v_k - v_1) = 0.
$$

(3) For $\beta=(\beta_f)_{f\in F}\in\R^F$, let
$P_\beta :=\{x\in\R^d\mid h_f(x)\leq \beta_f,\text{ for } f\in F\}$
be the polytope obtained from $P$ by parallel translations of the facets.
In particular, $P_\alpha=P$.
The {\it open facet deformation cone\/}\footnote{This is linearly isomorphic to the
{\it type-cone} of $P$ described by McMullen in \cite[\S2, p. 88]{McMullen}.}
 $D^{F,open}_P$ for $P$ is the
set of $\beta \in \R^F$ for which the polytopes $P_\beta$ and $P$
have the same normal fan $\N(P_\beta)=\N(P)$.  (Equivalently,
$P_\beta$ and $P$ have the same combinatorial structure.)  The {\it (closed) facet deformation cone\/}
is the closure $D^F_P$ of $D^{F,open}_P$ inside $\R^F$.
\end{definition}

It is clear that the definitions of $D^V_P$
and $D^E_P$  translate into linear equations and
weak linear inequalities.  Thus $D^V_P$ and
$D^E_P$ are (closed) polyhedral cones in the spaces
$(\R^d)^V$ and $\R^E$, correspondingly.
The following lemma shows that $D^F_P$ is also
a polyhedral cone.


\begin{lemma}
\label{lem:def_cone_F}
For a simple polytope $P$, the facet deformation cone $D^{F,open}_P$ is a full
$|F|$-dimensional open polyhedral cone inside $\R^F$,
that is a nonempty subset in $\R^F$ given by strict linear inequalities.
Thus $D^{F}_P$ is the closed polyhedral cone in $\R^{F}$ given by replacing
the strict inequalities with the corresponding weak inequalities.
\end{lemma}

\begin{proof}
Since every polytope $P_\beta$ has facet normals in directions which
are a subset of those for $P$, the rays (=1-dimensional normal cones) in
$\N(P_\beta)$ are a subset of those in $\N(P)$.  Therefore, one will have
$\N(P_\beta) = \N(P)$ if and only if $P_\beta, P$ have the same face lattices,
or equivalently, the same collection of vertex-facet incidences $(v,f)$.
This means that one can define the set $D^{F,open}_P$ inside $\R^F$
by a collection of strict linear inequalities on the coordinates $\beta=(\beta_f)_{f \in F}$.
It is next explained how to produce one such inequality for each pair
$(v_0,f_0)$ of a vertex $v_0$ and facet $f_0$ of $P$ such that $v_0\not\in f_0$.

If $v_0$ lies on the $d$ facets $f_1,\ldots,f_d$ in $P$, then
$v_0$ is the unique solution to the linear system $h_{f_j}(x)=\alpha_{f_j}$ for $j=1,...,d$.
Its corresponding vertex $x_0$ in $P_\beta$ is then the unique solution to
$h_{f_j}(x)=\beta_{f_j}$ for $j=1,...,d$.  Note that this implies $x_0$ has
coordinates given by linear expressions in the coefficients $\beta$.
Then the inequality corresponding to the vertex-facet pair $(v_0,f_0)$
asserts that $h_{f_0}(x_0) < \beta_{f_0}$.

Lastly, note that this system of strict linear inequalities has at least one solution, namely
the $\alpha$ for which $P_\alpha=P$.  Hence this open polyhedral cone is nonempty.
\end{proof}

The following theorem gives several different ways to describe deformations
of a simple polytope.

\begin{theorem}
\label{th:deformation}
Let $P$ be a $d$-dimensional simple polytope in $\R^d$, with notation as above.
Then the following are equivalent for a polytope $P'$ in $\R^d$:
\begin{enumerate}
\item[(i)] The normal fan $\N(P)$ refines the normal fan $\N(P')$.
\item[(ii)] The vertices of $P'$ can be (possibly redundantly)
labelled $x_v$, $v \in V$,
so that $(x_v)_{v\in V}$ is a point in the vertex deformation cone $D_P^V$,
i.e., the $x_v$ satisfy conditions~\eqref{eq:nonnegative-scaling_uv}.
\item[(iii)]
The polytope $P'$ is the convex hull of points
$x_v$, $v\in V$, such that $(x_v)_{v\in V}$ is in the vertex
deformation cone $D_P^V$.
\item[(iv)] $P'=P_\beta$ for some $\beta$ in the closed
facet deformation cone $D^F_P$.
\item[(v)] $P'$ is a Minkowski summand of a dilated polytope $rP$,
that is there exist a polytope $Q \subset \R^d$
and a real number $r>0$ such that
$P'+Q = rP$.
\end{enumerate}
\end{theorem}

\begin{proof}
One proceeds by proving the following implications
(i)
$\Rightarrow$
(ii)
$\Rightarrow$
(iii)
$\Rightarrow$
(i)
$\Rightarrow$
(iv)
$\Rightarrow$
(iii),
(iv)
$\Rightarrow$
(v)
$\Rightarrow$
(i).

{\sf (i) implies (ii).}
The refinement of normal fans gives rise to the redundant labelling of vertices
$(x_v)_{v \in P}$ as follows:  given a vertex $x$ of $P'$, label it by $x_v$
for every vertex $v$ in $P$ having its normal cone $\N_v(P)$ contained in the
normal cone $\N_{x}(P')$.  There are then two possibilities for an edge $uv \in
E$ of $P$:  either $x_u = x_v$, in which case~\eqref{eq:nonnegative-scaling_uv}
is trivially satisfied, or else $x_u \neq x_v$ so that $\N_u(P), \N_v(P)$ lie
in different normal cones $\N_{x_u}(P') \neq \N_{x_v}(P')$.  But then since
$\N(P)$ refines $\N(P')$, these latter two cones must share a codimension one
subcone lying in the same hyperplane that separates $\N_u(P)$
and $\N_v(P)$.  As this
hyperplane has normal vector $u-v$, this forces $x_u - x_v$ to be a positive
multiple of this vector, as desired.

{\sf (ii) implies (iii).}  Trivial.

{\sf (iii) implies (i).}
Let $P'$ be the convex hull of the points $x_v$.
Fix a vertex $u\in V$.  Let
$\lambda\in (\R^d)^*$ be a generic linear functional that belongs to
the normal cone $\N_u(P)$ of $P$ at the vertex $u$.
Then the maximum of $\lambda$ on $P$
is achieved at the point $u$ and nowhere else.  Let us orient the $1$-skeleton
of $P$ so that $\lambda$ increases on each directed edge.   This connected graph
has a unique vertex of outdegree 0, namely the vertex $u$.  Thus, for any other
vertex $v\in V$, there is a directed path $(v_1,\dots,v_l)$ from $v_1=v$ to
$v_l=u$ in this directed graph.  According to the conditions of the lemma, one has
have $\lambda(x_{v_1}) \leq \lambda(x_{v_2}) \leq \cdots \leq
\lambda(x_{v_l})$, so $\lambda(x_v) \leq \lambda(x_u)$.  Thus the maximum of
$\lambda$ on the polytope $P'$ is achieved at the point $x_u$.  This implies
that $x_u$ is a vertex of $P'$
and the normal cone $\N_{x_{u}}(P')$ of $P'$ at
this point contains the normal cone $\N_u(P)$ of $P$ at $u$.  Since the same
statement is true for any vertex of $P$, one deduces that $\N(P)$ refines
$\N(P')$.

{\sf (i) implies (iv).}  First, note that if $\N(P')=\N(P)$
then $P'=P_\beta$ for some $\beta$ in the open facet deformation cone
$D^{F,open}_P$.   Indeed, the facets of $P'$ are orthogonal to
the 1-dimensional cones in $\N(P')$, thus they should be parallel
to the corresponding facets of $P$.

Now assume that $\N(P)$ refines $\N(P')$.
Recall the standard fact \cite[Prop. 7.12]{Ziegler}
that the normal fan $\N(Q_1+ Q_2)$ of a Minkowski sum $Q_1+Q_2$
is the common refinement of the normal fans $\N(Q_1)$ and $\N(Q_2)$.
Thus, for any $\epsilon>0$,
the normal fan of the Minkowski sum
$P'+\epsilon P$ coincides with $\N(P)$.
By the previous observation, one should have
$P'+\epsilon P = P_{\beta(\epsilon)}$ for some $\beta(\epsilon)
\in D^{F,open}_P$.  Since all coordinates of $\beta(\epsilon)$
linearly depend on $\epsilon$, one obtain $P'= P_\beta$ for $\beta =
\lim_{\epsilon\to 0}\beta(\epsilon)\in D^{F}_P$.

{\sf (iv) implies (iii).}
Given $\beta \in D^F_P$, it is the limit point for some family
$\{\beta(\epsilon)\} \subset D^{F,open}_P$.
One may assume that $\beta(\epsilon)$ linearly depends on
$\epsilon>0$ and $\lim_{\epsilon\to 0}\beta(\epsilon) = \beta$.
Hence $P' = P_\beta$ is the limit of the polytopes
$P_{\beta(\epsilon)}$, which each have
$\N(P_{\beta(\epsilon)})=\N(P)$.
In particular, the vertices of $P_{\beta(\epsilon)}$ can be
labelled by $x_v(\epsilon)$, $v\in V$.
These vertices linearly depend on $\epsilon$ and satisfy
$x_u(\epsilon)-x_v(\epsilon) = \R_{\geq 0} (u-v)$ for any edge $uv\in E$.
Taking the limit $\epsilon\to 0$, one obtains that $P'$ is the convex
hull of points $x_v = \lim_{\epsilon\to 0} x_v(\epsilon)$
satisfying~\eqref{eq:nonnegative-scaling_uv}.

{\sf (iv) implies (v).}  Note that $P_{\gamma} + P_{\delta} =
P_{\gamma+\delta}$, for $\gamma,\delta\in D_P^F$.
Let $P' = P_\beta$ for $\beta\in D_P^F$.
The point $\alpha$ (such that $P= P_\alpha$) belongs to the {\it open\/}
cone $D^{F,open}_P$.  Thus, for sufficiently large $r$, the point
$\gamma = r\alpha - \beta$ also belongs to the open cone $D^{F,open}_P$.
Let $Q = P_\gamma$.  Then one has $P'+Q = P_\beta + P_{r\alpha-\beta} =
P_{r\alpha} = rP$, as needed.

{\sf (v) implies (i).}
This follows from the standard
fact \cite[Prop. 7.12]{Ziegler} on normal fans of Minkowski sums
mentioned above.
\end{proof}

\begin{remark}
We are being somewhat careful here, since Theorem~\ref{th:deformation} can
fail when one allows a broader interpretation for a simple polytope $P$
to deform into a polytope $P'$
by parallel translations of its facets, e.g. if one allows facets to
translate {\it past\/} vertices.
For example,  letting $P'$ be a regular tetrahedron in $\R^3$, and $P$
the result of ``shaving off an edge'' from $P'$ with a generically tilted plane in $\R^3$,
one finds that $\N(P)$ does not refine $\N(P')$.
\end{remark}


Let us now describe the relationship between the three deformation
cones $D_P^V$, $D_P^E$, and $D_P^F$.
Let $H$ be the linear subspace in $(\R^d)^V$
given by
$$
H:=\{ (x_v)_{v\in V} \in (\R^d)^V\mid x_u-x_v \in \R(u-v)
\text{ for any edge }uv\in E \}.
$$
Clearly, the cone $D_P^V$ belongs to the subspace $H$.
Let us define two linear maps
$$
\phi: H \to \R^E\quad\textrm{and}\quad
\psi:\R^F\to H.
$$
The map $\phi$ sends $(x_v)_{v\in V} \in H$ to $(y_e)_{e\in E}\in \R^E$,
where $x_u-x_v = y_e(u-v)$, for any edge $e=uv\in E$.
The map $\psi$ sends $\beta = (\beta_f)_{f\in F}$ to $(x_v)_{v\in V}$,
where, for each vertex $v$ of $P$ given as the intersection
of the facets of $P$ indexed $f_1,\dots,f_d$, the point $x_v\in \R^d$
is the unique solution
of the linear system $\{h_{f_j}(x) = \beta_{f_j}\mid j=1,\dots,d\}$.
For $\beta \in D_P^{F,open}$, $\psi(\beta)= (x_v)_{v\in V}$,
where the $x_v$ are the vertices of the polytope $P_\beta$.
One can easily check that $\psi(\beta)\in H$.
Indeed, this is clear for $\beta\in D_P^{F,open}$ and thus this
extends to all $\beta\in \R^F$ by linearity.

Note that the kernel of the map $\phi$ is the subspace
$\Delta(\R^d)\simeq\R^d$ embedded diagonally into $(\R^d)^V$.
This follows from the fact that the $1$-skeleton of $P$ is connected.
The vertex deformation cone $D_P^V$ can be reduced
modulo the subspace $\Delta(\R^d)$ of parallel translations of polytopes.
Similarly, the facet deformation cone can be reduced modulo
the subspace $\Delta'(\R^d) = \psi^{-1}(\Delta(\R^d))\simeq\R^d$,
where $\Delta'(x) := (h_f(x))_{f\in F}$ for $x\in \R^d$.

\begin{theorem}  The map $\psi$ gives a linear isomorphism between
the cones $D_P^F$ and $D_P^V$.  The map $\phi$ induces a linear isomorphism
between the cones $D_P^V/\Delta(\R^d)$ and $D_P^E$.
Thus one has
$$
D_P^E \simeq D_P^V/\Delta(\R^d) \simeq D_P^F/
\Delta'(\R^d).
$$
In particular, $\dim D_P^E = \dim D_P^V - d = \dim D_P^F - d = |F| -d$.
\end{theorem}

\begin{proof}  The claim about the map $\psi$ follows immediately from
Theorem~\ref{th:deformation}.

Let us prove the claim about $\phi$.
Note that, for $(x_v)_{v\in V}\in D_P^V$,
the point $(y_e)_{e\in E} =
\phi((x_v)_{v\in V})$ satisfies the condition
of Definition~\ref{def:three_cones}(2)
because
$$
\begin{aligned}
&y_{e_1} (v_1-v_2) + y_{e_2} (v_2-v_3) + \cdots + y_{e_k}(v_k-v_1) \\
&= (x_{v_1}-x_{v_2})  + (x_{v_2}-x_{v_3})  + \cdots + (x_{v_k} - x_{v_1}) \\
&= 0.
\end{aligned}
$$

It remains to show that for any $(y_e)_{e\in E}\in D_P^E$
there exists a unique (modulo diagonal translations)
element $(x_v)_{v\in V}\in H $ such that $x_u-x_v = y_e (u-v)$
for any edge $e=uv\in E$.
Let us construct the points $x_v\in \R^d$, as follows.
Pick a vertex $v_0\in V$ and pick any point $x_{v_0}\in \R^d$.
For any other $v\in V$, find a path $(v_0, v_1,\dots,v_l)$
from $v_0$ to $v_l=v$ in the $1$-skeleton of $P$ and define
$x_v = x_{v_0} - y_{v_0v_1}(v_0-v_1) - y_{v_1 v_2}(v_1-v_2)-\cdots
- y_{v_{l-1} v_l}(v_{l-1}-v_l)$.
This point $x_v$ does not depend on a choice of path from $v_0$ to $v$,
because any other path in the $1$-skeleton can be obtained by switches
along 2-dimensional faces of $P$.  These $\left( x_v \right)_{v\in V}$
satisfy the needed conditions.

Finally, note that $\dim D_P^F = |F|$ because $D_P^F$ is a full-dimensional
cone.
\end{proof}

\end{document}